\pdfoutput=1
\documentclass[12pt]{amsart}
\usepackage{latexsym}
\usepackage{amsmath}
\usepackage{amssymb}
\usepackage{color}
\usepackage{amssymb, amscd}
\usepackage{ifthen}
\usepackage{amsthm}
\input epsf

\newtheorem{dfn}{Definition}[section]
\newtheorem{defi}[dfn]{Definition}

\newtheorem{rem}[dfn]{Remark}

\newtheorem{rems}[dfn]{Remarks}
\newtheorem{thm}{Theorem}

\newtheorem{theor}{Theorem}[section]

\newtheorem{lem}[dfn]{Lemma}

\newtheorem{sublem}[dfn]{Sublemma}
\newtheorem{prop}[dfn]{Proposition}
\newtheorem{cor}[dfn]{Corollary}

\newtheorem{conv}[dfn]{Convention}

\def\proof{\par\medskip\noindent{\it Proof: }}

\def\R{{\Bbb R}}
\def\H{{\Bbb H}}
\def\Z{{\Bbb Z}}

\def\A{{\mathcal A}}

\def\N{{\Bbb N}}
\def\D{{\Delta}}
\def\F{{{\mathcal F}}}
\def\T{{{\mathcal T}}}
\def\P{{{\mathcal P}}}

\def\al{\alpha}
\def\be{\beta}
\def\ga{\gamma}
\def\Th{\Theta}
\def\g{\gamma}
\def\Ga{{\mathcal G}}
\def\G{{\Gamma}}
\def\de{\delta}
\def\d{\delta}
\def\da{d_A}

\def\ve{\varepsilon}

\def\De{\Delta}

\def\L{{\mathcal L}}

\def\Ca{{\mathcal C}a}
\def\la{\lambda}
\def\La{\Lambda}

\def\d{\delta}

\def\ti{\widetilde}

\def\om{\omega}
\def\D{\partial}

\def\iso4{{\rm Isom}{\H^4}}
\def\hn{{\H^n}}
\def\h3{{{\Bbb H}^3}}
\def\h4{{\H^4}}
\def\sp{\S^{n-1}_\infty}
\def\sp3{{\S}^{3}_\infty}

\def\bx{$\hfill\square$}
\def\ent{{\rm Ent}~T}

\def\de{\Delta_{\bf e}}
\def\d{\delta}
\def\be{{\bf e}}
\def\ba{{\bf a}}
\def\bb{{\bf b}}

\def\smc{{\rm Small}(\bc)}

\def\bp{{\bf p}}
\def\bc{{\bf c}}

\def\bd{{\bf d}}
\def\bz{{\bf z}}
\def\by{{\bf y}}
\def\bw{{\bf w}}
\def\zx{{\bf x}}
\def\bg{{\bf g}}
\def\bh{{\bf h}}
\def\bu{{\bf u}}
\def\bv{{\bf v}}

\def\lba{{\mathsf L}\ba}
\def\ld{{\rm L}\bd}

\def\db{\Delta_{\bf b}}
\def\dc{\Delta_{\bf c}}
\def\od{{\overline\delta}}
\def\Sh{{\rm Sh}}
\def\sh{{\rm sh}}

\def\BF{\partial_fG}

\def\OFG{{\overline \G}_f}
\def\OFH{{\overline H}_f}

\def\1GF{\G_1\cup \partial_f\G_1}
\def\shab{\sh_{\ba}{\bb}}
\def\Shab{\Sh_{\ba}{\bb}}

\def\shba{\sh_{\bb}{\ba}}

\def\shbc{\sh_{\bb}{\bc}}
\def\shdb{\sh_{\bd}{\bb}}
\def\shdc{\sh_{\bd}{\bc}}
\def\shcb{\sh_{\bc}{\bb}}
\def\shda{\sh_{\bd}{\ba}}
\def\shca{\sh_{\bc}{\ba}}

\def\Shcb{\Sh_{\bc}{\bb}}
\def\Shca{\Sh_{\bc}{\ba}}

\def\shbue1{\sh^{\ve_1}_{\bb}{\bu}}

\def\sm0{{\rm Small}_0^{\ve}}
\def\sm1{{\rm Small}_1^{\ve}}

\def\act{\curvearrowright}
\def\GAT{G\act T}

\def\DFGa{\D_f\G}

\def\iso3{{\rm Isom}{\H^3}}
\def\hn{{\H^n}}
\def\lfv{L_{f,v}}

\def\bsn{{\bf S}^nT}
\def\bs2{{\bf S}^2T}
\def\bden{{\bf \Delta}^nT}
\def\bde2{{\bf \Delta}^2T}
\def\odg{\overline{\d}_g}

\def\tupl#1{\langle #1\rangle}

\def\OG{\overline {\G}_f}
\def\T{\mathcal T}
\def\rh{\mathsf{RH_{32}}}

\subjclass[2010]{Primary 20F65, 20F67; Secondary 57M07, 22D05}

\title{Non-finitely generated relatively hyperbolic groups and Floyd quasiconvexity}

\author{Victor Gerasimov}
\address{Victor Guerassimov, Departamento de Matem\'atica,
Universidade Federal de Minas Gerais,
Av. Antonio Carlos, 6627/ Caixa Postal 702  31270-901 Belo Horizonte, MG,
Brasil}
\email{victor@mat.ufmg.br}
\author{Leonid Potyagailo}
\address{ Leonid Potyagailo, UFR de Math\'ematiques, Universit\'e
de Lille 1, 59655 Villeneuve d'Ascq cedex, France}
\email{potyag@math.univ-lille1.fr}

\begin{document}


\def\pdFile#1#2#3{\ifnum\pdfoutput>0%
\pdfximage width #2pt height
#3pt{#1.pdf}\pdfrefximage\pdflastximage\else\kern#2pt\vbox to
#3pt{\vss}\fi}

\def\Vnote#1{\underline{\{V:\color{red}\small #1\}}}
\def\Excl#1{{\{\color{green} #1\}}}

\date{\today}
\maketitle

\begin{abstract}

\noindent We regard a relatively hyperbolic group as a group
acting non-trivially by homeomorphisms on a compactum $T$
discontinuously on the set of distinct triples and cocompactly on
the set of distinct pairs of points  of $T$.

In the first part of the paper we prove that such a group $G$
admits a graph of groups decomposition  given by a star graph
whose central vertex group is finitely generated relatively
hyperbolic with respect to the edge groups, and the other vertex
groups are stabilizers of non-equivalent parabolic points. It
follows from this result that every relatively hyperbolic group is
relatively finitely generated with respect to the parabolic
subgroups. Another corollary is that the definition of the
relative hyperbolicity which we are using is equivalent to those
of Bowditch and Osin (taken with respect to finitely many
peripheral subgroups) and they are all equivalent to the existence
of the above star graph of groups decomposition.

The second part of the paper uses the method of the first part.
Considering the induced action of $G$ on the space of distinct
pairs of $T$ we construct a connected graph on which $G$ acts
properly and cofinitely on edges. Equipping the graph with   Floyd
metrics we
  prove that the
 quasigeodesics in this metric are   close somewhere to the geodesics in the
 word metric. This allows us to prove that the parabolic subgroups of $G$ are quasiconvex
with respect to the Floyd metrics. As a corollary we prove that
the preimage of a parabolic point by the Floyd map is the Floyd
boundary of its stabilizer.
\end{abstract}

\section{Introduction}

 \noindent {\bf Part I of the Paper.} Let $T$ be a compact Hausdorff space
(compactum)
 containing at least $3$ points. The action of a  discrete  group $G$
  by homeomorphisms  of $T$ is called \it convergence action \rm
 if the induced action on the space $\Th^3T$ of
subsets of cardinality $3$  is discontinuous. We   say in this
case that the action is 3\it-discontinuous\rm.

The action of $G$ on $T$ is called 2\it-cocompact \rm
if the
 action on $\Th^2T$
  is cocompact. An action is called \it parabolic \rm
 if
$G$ is infinite
there is a unique fixed point.

If $G$ admits a non-parabolic action on $T$ which is
3-discontinuous and 2-cocompact
  then the action is
\it geometrically finite\rm, i.e.  every point of  $T$ is either
conical  or bounded parabolic or isolated \cite{Ge1}.
 Conversely if a group $G$
admits a minimal geometrically finite action on a {\bf metrisable}
space $T$ without isolated points then the action is $2$-cocompact
\cite{Tu3}. If $G$ is finitely generated then the existence of a
geometrically finite action of $G$ is equivalent (see \cite{Bo1}
and \cite{Ya}) to the ``classical'' relative hyperbolicity in the
sense of Farb \cite{Fa} and Gromov \cite[8.6]{Gr}.

These facts
justify the following ``dynamical'' definition.

\begin{dfn}\label{ourdfn} \cite{Ge1} A group $G$ is called
\it relatively hyperbolic \rm if it admits a non-parabolic
$3$-discontinuous and 2-cocompact action  ($\rh$-action) on a compactum $T$.
\end{dfn}

We point out that we do not impose any restriction on the
cardinality of $G$. We also do not require  the metrisability of
$T.$

 Our first result shows  that any
  relatively hyperbolic group can be ``nicely'' approximated by
  finitely generated relatively hyperbolic group.

\begin{thm}
\label{graph}
  Let $G$ be a relatively hyperbolic group  with respect to a
 collection of
parabolic subgroups $\{P_1,..., P_n\}.$ Then   $G$  is the
fundamental group of the following finite ``star graph''

\vskip-5pt
$$\vcenter{\pdFile{grph0}{200}{180}}
\eqno(1)$$
\vskip5pt
 \noindent  whose central
  vertex group  $G_0$ is   finitely generated relatively hyperbolic  with respect to
  those edge groups $Q_i=P_i\cap G_0$ which are infinite,  all other vertex groups
of the graph   are $P_i\ (i=1,...,n)$.

Moreover for every finite set $K\subset G$ the subgroup
$G_0$ can be chosen to contain $K$.
 \end{thm}

\medskip

Theorem A yields generalization of several known results omitting
the assumption of finite generatedness.

A group $G$ is said to be \it finitely generated with respect to \rm a collection
$\mathcal H$ of subgroups if there exists a finite set $S{\subset}G$ such that
$S{\cup}({\cup}\mathcal H)$ is a generating set for $G$.

\bigskip

\noindent  {\bf Corollary} (Corollaries \ref{fingen},
\ref{noparab}).
\it
 Let a group $G$ admit a $3$-discontinuous $2$-cocompact non-parabolic action
on a compactum $T$.
Then $G$ is finitely generated with respect to
a finite collection of the stabilizers
of
parabolic points.
In particular, if $G$ acts without parabolics then
$G$ is finitely generated\rm.

\bigskip

In \cite[Appendix]{GePo1} we gave a short proof
 of Bowditch's
theorem that the existence of a 3-discontinuous and 3-cocompact action of a
finitely generated group implies that the group is hyperbolic.
 The
above Corollary omits the
assumption of finite generatedness.

In most papers about relatively hyperbolic groups the authors assume
that the group is finitely generated.
Besides definition \ref{ourdfn}
there are two more definitions which do not require the finite generatedness.
The first is due to
B.~Bowditch \cite{Bo1} and the second is due to D.~Osin \cite{Os}.
We recall them now.

A graph $\Gamma$ is called \it fine \rm
if for any two vertices the set of arcs of fixed length
joining them is finite.
Bowditch calls a group $G$ relatively hyperbolic
if there is an action of $G$ on a fine hyperbolic graph $\Gamma$
such that the action on edges is \it proper \rm(i.e. the edge
stabilizers are finite), \it cofinite \rm(the
set of edge orbits is finite i.e. $\vert\Gamma^1/G\vert
<\infty$) and \it non-parabolic \rm(there is no vertex fixed by $G$).

We  use a "finite" version of  Osin's definition of relative
hyperbolicity according to which
  a group $G$ is relatively hyperbolic with respect to a
finite collection $\mathfrak P$ of subgroups of infinite index, if
 it is relatively finitely presented with respect to $\mathfrak P$
 and satisfies
linear isoperimetric inequality relative to this system (see
Definition \ref{osdef} in Section 3.5). We note that the original
Osin's definition does not require the finiteness of the
peripheral system $\mathfrak P$ whereas all other definitions
imply this property. So to relate   Osin's definition with all
other definitions we will always assume that the system $\mathfrak
P$ is finite. The assumption  that every element of $\mathfrak P$
is a subgroup of infinite index is needed to exclude the trivial
case of the relative hyperbolicity with respect to a subgroup of
finite index.
The existence of the star graph decomposition (1) directly follows
from  Osin's definition  (see \cite[Theorem 2.44]{Os}). Our
Theorem A is a different result as it uses another definition of
the relative hyperbolicity. On its turn since the existence of
such a graph of groups decomposition is a common point for both
these approaches it gives rise to the following
 equivalence of all known definitions of the relative
hyperbolicity  valid for  a group without any restriction on its
cardinality.

\bigskip

\noindent {\bf Corollary} (Theorem \ref{equiv}).   The following
conditions of the relative hyperbolicity for a group $G$ are
equivalent:

 \begin{itemize}

   \item[\sf 1)] (the above Definition \ref{ourdfn}) The group $G$ admits a $\rh$-action on a
  compactum $T$ containing at least $3$ points.

 \item[\sf 2)] (Bowditch's definition) The group $G$ acts non-parabolically on a
 connected fine hyperbolic graph $\Gamma$ properly and cofinitely on edges.

 \item[\sf 3)] (Osin's definition) The group $G$ is relatively finitely
  presented and admits a relative linear isoperimetric
  inequality relatively to
   a finite system of   subgroups of infinite index.

\item[\sf 4)] $G$ admits the star graph decomposition (1) where
the central vertex group $G_0$ is a finitely generated relatively
hyperbolic group with respect to those   edge groups $Q_i$ which
are infinite.
\end{itemize}
\bigskip


  The implication $1)\Rightarrow 4)$
follows from Theorem A and $3)\Rightarrow 4)$ from \cite[Theorem
2.44]{Os}. The implication $2)\Rightarrow 1)$ is proved in
\cite{Ge2}. The following proposition  yields the implication
$4)\Rightarrow 2).$

\medskip

\noindent  {\bf Proposition} (Proposition \ref{inverse}). {\it
Suppose that a group $G$ admits a graph of groups decomposition
(1) where the group $G_0$ is finitely generated and relatively
hyperbolic with respect to the subgroups $Q_i\ (i=1,,,n)$. Then
$G$ satisfies Bowditch's definition.}

\medskip

     The proof of the
implication $2)\Rightarrow 3)$ is an easy use of the common
methods for hyperbolic metric spaces \cite{Gr}. We include it for
the completeness avoiding the references to the sources in which
it is not clear that the assumption of finite generatedness is
inessential.

We resume all this discussion in the following diagram.

\vskip-7pt
\begin{picture}(300,110)(-120,-5)
\put(100,73){\makebox(0,0)[cb]{\cite{Ge2}}}
\put(0,70){\oval(80,30)}
\put(-35,67){\sf 2)}
\put(5,70){\makebox(0,0)[c]{\vbox{\hbox{\kern3pt
Bowditch's}\hbox{\kern5pt definition}}}}
\put(200,70){\oval(80,30)}
\put(167,66){\sf 1)}
\put(205,70){\makebox(0,0)[c]{$\rh$}}
\put(45,70){\vector(1,0){110}}
 \put(0,52){\vector(0,-1){35}}
\put(200,52){\vector(0,-1){35}}
\put(-2,35){\makebox(0,0)[cr]{Proposition {\ref{bowos}}}}
\put(203,35){\makebox(0,0)[cl]{Theorem A}}
\put(333,35){\makebox(0,0)[cl]{$(*)$}}
 \put(0,0){\oval(80,28)}
\put(-32,-4){\sf 3)}
\put(5,0){\makebox(0,0)[c]{\vbox{\hbox{\kern7pt
Osin's}\hbox{definition}}}}
\put(107,2){\makebox(0,0)[cb]{\cite{Os}}}
\put(45,0){\vector(1,0){120}}
 \put(200,0){\oval(60,28)}
\put(74,40){\makebox(0,0)[lb]{Proposition {\ref{inverse}}}}
\put(177,-3){\sf 4)}
\put(205,0){\makebox(0,0)[c]{\vbox{\hbox{\kern2pt
Star}\hbox{graph}}}}
 \put(168,8){\vector(-3,1){133}}
\end{picture}
\vskip20pt

We   note that the vertex groups corresponding to  the non-central
vertices of the star graph in the above Proposition \ref{inverse}
can be uncountable. This provides a construction of an uncountable
relatively hyperbolic group  too.  The proof of Theorem A
presented below does not depend on the cardinality of $G$ nor on
the metrisability of the space on which it acts, and is
self-contained. It is based on the theory of \it entourages \rm of
a compactum $T$ which are the neighborhoods of the diagonal of
$T^2$. In Section 3 using a $G$-orbit $A$ of entourages on $T$ we
construct a graph $\Ga$ on which $G$ acts and whose set of
vertices is $A$.
 The subgroup
$G_0$ will be chosen as the stabilizer of a connected  component
of a refined graph $\ti\Ga$ having the same set of vertices: $\ti
\Ga^0=\Ga^0=A.$ We will use a system of  \it tubes \rm and \it
horospheres \rm on $\Ga$ to establish the existence of the
requested splitting of $G$ as a star-graph of groups.

\vspace*{2mm}

\noindent {\bf Historical remarks and comments.} For the
completeness of the exposition we provide  a short survey of known
results related to Theorem \ref{equiv}. We start with a less
general (but more standard) case of a finitely generated group and
then describe briefly what is known when $G$ is a non-finitely
generated group.

\vspace*{2mm}

 \noindent {\it Case 1. $G$ is finitely generated.\rm}
The equivalence of the conditions 2) and 3) was proved by
F.~Dahmani \cite{Da} and D.~Osin \cite{Os}.

B.~Bowditch proved that the condition 2) implies that $G$ acts
properly discontinuously by isometries on a proper
 hyperbolic metric space $X$, and the action on the boundary  $\partial X$
 is {\it geometrically finite}
 meaning that every point of $\partial X$
  is either conical or
bounded parabolic. A strengthened converse statement was proved by
A.~Yaman \cite{Ya}. She showed that a group  that possesses a
 geometrically finite convergence action  on a non-empty metrisable
perfect compactum $T$  such that   the stabilizers of parabolic
points are all finitely generated satisfies the condition 2). Note
that the finite generatedness of the maximal parabolic subgroups
implies  by Corollary \ref{fingen} that the whole group $G$ is
finitely generated.

From the other hand  a minimal action on a metrisable compactum is
an $\rh$-action if and only if it is geometrically finite. Indeed
the sufficiency  follows from P.~Tukia's result \cite[Theorem
1.C]{Tu3}. The converse  statement is a partial case of \cite[Main
Theorem, b]{Ge1}.

So the conditions 1), 2) and 3) are equivalent  if $G$ is finitely
generated.  By \cite[Lemma 2.46]{Os}   the implication
$4)\Rightarrow 3)$ is  true  for any $G$ (we thank the referee for
this reference). On its turn 4) trivially holds for every finitely
generated relatively hyperbolic group.

Note that an alternative proof of Yaman's theorem in the finitely
generated case is given in \cite[Corollary of 7.1.1]{GePo3}.

\vspace*{1mm}

\noindent {\it Case 2. $G$ is countable.\rm} In \cite{Hr}
C.~Hruska\ pointed out   that the proofs of the equivalence
between
  the conditions 2) and 3) given in \cite{Da} and \cite{Os}
 remain true for   countable groups. However their relation
 with the   geometrical finitenness  is
 more delicate already in this case.  C. Hruska noticed that
    the
    proof of the above theorem of  Bowditch
    does not work if the parabolic subgroups are not finitely
    generated \cite[Remark after  5.6]{Hr}. He indicated
   how to  generalize  the methods of the paper of D.~Groves and J.F.~Manning
    \cite{GrMa} and to prove the implication $3)\Rightarrow  2)$
    (in fact the argument gives a stronger statement that
    3) implies Gromov's definition  of the relative hyperbolicity denoted by (RH-3) in
\cite{Hr}).

It is   claimed  without proof in \cite{Hr} that Yaman's theorem
remains valid in the countable case (note that in \cite{Hr} this
statement was misleadingly denoted    (RH-1) $\Rightarrow$ (RH-2)
but it  should   be (RH-1)$\Rightarrow$ (RH-4)). The main part of
Yaman's proof consists in generalizing the statements of the paper
\cite{Bo3} about hyperbolic groups to the case of  relatively
hyperbolic groups.  In particular Yaman uses some lemmas of
\cite{Bo3} when the group is a posteriori finitely generated. So
it seems to us that these arguments still require further
explanations in the case when the group is not finitely generated.

Note that if one admits that the proof of Yaman works in the
countable case then it would give a proof of  our Theorem A in
this case. Indeed $1)$ yields that the action is geometrically
finite by \cite{Ge1}. Then Yaman's theorem would imply  2).
Finally we obtain $2)\Rightarrow 3)$ as it is described in
\cite{Hr} and  $3)\Rightarrow 4)$ by \cite[Theorem 2.44]{Os}.

The argument of Hruska generalizing the theorem of \cite{GrMa}
 requires  the metrisability of the compactum $T$ which is
homeomorphic to the boundary of a hyperbolic space (given by
Gromov's definition (RH3)). Any group $G$ admitting a
3-discontinuous action on the metrisable compactum  $T$  is
countable \cite[ Corollary 2, section 5.3]{Ge1} (note that the
converse statement is  true
  if  one supposes that the action of $G\act T$ is
 $\rh$ then the countability of $G$ implies the metrisability of $T$
 \cite[Main Theorem, c]{Ge1}).  So the condition to be countable
for a relatively hyperbolic group seems to be unavoidable in this
approach.

\vspace*{1mm}

 \noindent {\it Case 3. $G$ is  an arbitrary group.\rm}

  As we have mentioned  the equivalence
$3)\Leftrightarrow 4)$ is true for any group \cite[2.44 $\&$
2.46]{Os}. The proof  of the equivalence of the conditions 3) and
4) to the condition 2) is not so difficult. Since the arguments
are spread in different papers and
 sometimes  require   modifications,
we included them in  Propositions \ref{inverse} and \ref{bowos}.

 The relation with the dynamical condition $\rh$
(or with the geometrical finitenness)  was not
 known  before. Thus the main
result of the Section is  Theorem A
 which establishes (with the statement $2)\Rightarrow
1)$ from \cite{Ge2})   the equivalence of the condition $1)$
 to all  other conditions.

 One of the
  difficulties of the situation is that the
condition $\rh$ still implies the geometrically finiteness by
\cite{Ge1} but the converse statement is not known in this case
(the argument of Tukia certainly needs the metrisability of the
compactum on which the group acts).

The implication $1)\Rightarrow 3)$ follows from the above
Corollary. It   generalizes  Yaman's theorem to the case of an
arbitrary group admitting an $\rh$-action (see also
\cite[Proposition 7.1.2]{GePo3}).  In particular if $G$ is
countable together with Tukia's theorem it yields a proof of
  Yaman's theorem in this case.

We note that despite that   Osin's theorem  and Theorem A have the
same conclusion (condition 4)) and their assumptions (conditions
3) and 1) respectively) are equivalent, none of them
 is a corollary of the other one as
the proof of this equivalence uses both  statements.

 Notice also
that the star-graph decompositions of relatively hyperbolic groups
have been  used in \cite{Os} to reduce the case of a non-finitely
generated relatively finitely presented group to the case of a
finitely generated one.

We finish this discussion by the following question asking whether
 Tukia's theorem  remains valid without assuming the
metrisability of the space:

\vspace{2mm}

\noindent {\bf Question.} Is   it true that a geometrically finite
non-elementary minimal action on a compactum  is an $\rh$-action
?\bx

\vspace{2mm}  A positive answer to this question would imply in
particular that Yaman's theorem is true in the non-metrisable case
too.

\bigskip

\noindent {\bf Part II of the Paper.}
It
deals with \it finitely generated \rm
 relatively hyperbolic groups. It is
based on the methods developed in the first part.
Starting with Section 4 we use the Floyd completion of locally
finite graphs.
Let $\G$ be a
locally finite, connected graph admitting a cocompact and
discontinuous action of a finitely generated group $G$ (e.g. a
Cayley graph of $G$ or the graph of entourages $\Ga$). According to
W.~Floyd by rescaling the graph distance $d$ of $\G$ by  a scalar
function $f:\N\to \R_{\geq0}$ one obtains the Cauchy completion
$\OFG$ of the metric space $(\G, \d_f)$ where $\d_f$ is the
rescaled metric. We call this space {\it Floyd completion} (see
Section 4). The action of $G$ extends continuously to
$\OFG$. By \cite{Ge2} there exists an equivariant continuous map
$F$ from the \it Floyd boundary \rm $\partial_f\G=\OFG\setminus
\G$ to the space $T$. The kernel of the map $F$ was described in
[GePo1, Theorem A].
 Namely if the preimage of a point $p$
is not a single point then $p$ is
parabolic and the preimage coincides with
 the topological boundary of the
stabilizer ${\rm Stab}_Gp$ of $p$. We
denote by $\partial_f{\rm Stab}_Gp$
  the Floyd boundary of ${\rm Stab}_Gp$ corresponding to a
function $f$.

A subset $X$ of $\G$ is called {\it Floyd (r-)quasiconvex} if
every Floyd geodesic (with respect to the metric $\d_f$) with the
endpoints in $X$ belongs to $r$-neighborhood $N_r(X)$ for the
graph metric $d$ and some $r>0.$ In particular if
$f$ is the identity then  the Floyd quasiconvexity
means the standard one.
 It is well-known  that the parabolic subgroups are quasiconvex with respect to
$d$ \cite{DS}   (for another proof see e.g. [GePo1, Corollary
3.9]). Our next Theorem establishes the Floyd quasiconvexity of
the parabolic subgroups.

\bigskip

 \noindent {\bf Theorem C.} {\it Let  $G$ be  a finitely
generated group acting $3$-discontinuously and $2$-cocompactly on
a compactum $T.$
 Let $\G$
be a locally finite, connected graph  admitting a cocompact
discontinuous action of $G.$ Then there exists a  Floyd scaling
function $f$, such that   every parabolic subgroup $H$ of $G$ is
Floyd quasiconvex for the Floyd metric $\d_f$.}\bx

\bigskip

\noindent As a consequence of Theorem C we obtain the following
Corollary which  answers  our question [GePo1,
  1.1]:
\bigskip

\noindent{\bf Corollary \ref{kermap}} {\it For a scaling function
$f$ satisfying conditions $(1-3)$ (see Section (7)) one has

$$F^{-1}(p)=\partial_f({\rm Stab}_Gp)$$

\noindent for every parabolic point  $p\in T$.} \bx

\bigskip

 Note that it was already known that the map $F$
 is $1$-to-$1$ at conical points \cite{Ge2}.
 Corollary \ref{kermap} gives a complete description of  the preimage of a
 parabolic point by $F$ as the Floyd boundary of its stabilizer. It
 gives rise to a complete generalization of
 the Floyd theorem \cite{F}
 to the case of relatively hyperbolic groups.

 The proof of Theorem C (and Corollary \ref{kermap}) in Section 7 and is based on a
description of a family of \it tight \rm curves  which are
quasigeodesics locally everywhere and geodesic outside the
horospheres (see Definition \ref{tcurve}). Their properties are
described in the following Theorem (see Section 6 for  more
details):

\bigskip

\noindent {\bf Theorem B.} {\it For every  tight curve $\ga$ in
the graph  of entourages $\Ga$ there exists   a quasigeodesic
$\al\subset A$ such that every non-horospherical vertex of $\ga$
belongs to a uniform neighborhood of $\al.$}

\bigskip

The main step in proving Theorem C is to show that every Floyd
quasigeodesic is tight. We notice that the graph of entourages
$\Ga$ plays here a special role and in the proofs of   Theorems B
and C we deal mainly with  it.

 \noindent This is our second paper in a series of papers about
relatively hyperbolic groups. Keeping  the same definition of the
relative hyperbolicity here  we apply however different  methods
 based on the theory of discrete
systems of entourages not used in  \cite{GePo1}.

\medskip

\medskip
{\bf Acknowledgements.} During the work on this paper both authors
were partially supported by the ANR grant ${\rm BLAN}~07-2183619.$
We are  grateful to the Max-Planck Institute f\"ur Mathematik in
Bonn, where a part of the work was done. We also thank the CNRS and the
Brasilian-French cooperation grant for support.

The authors are thankful to Wenyuan Yang for very useful remarks
and corrections. We also thank the referee for a careful reading
of our paper and providing us a list of corrections.

\section{Convergence Groups} By {\it compactum} we   mean a
compact Hausdorff space. Let $\bsn$ denote the quotient
 of the   product space
$\underbrace{T{\times}\dots{\times}T}_{\mbox{$n$ times}}$ by the action of the
permutation group on $n$ symbols. The elements of $\bsn$ are
generalized  unordered $n$-tuples (i.e. an element may belong to a
tuple with some multiplicity). Let $\Th^nT$ be the subset of
$\bsn$ whose elements are non-ordered $n$-tuples with all distinct
components. Put  $\bden=\bsn\setminus\Th^nT$, the set $\bde2$ is
just the diagonal of $T^2$.

\bigskip

\noindent {\bf Convention.} {\it If the opposite is not stated all
group actions on compacta are assumed   to   have the convergence
property. }

\bigskip

We refer to  \cite{Bo2}, \cite{GePo1}, \cite{GM}, \cite{Fr},
\cite{Tu2} where standard facts related to the convergence groups
are proved. We recall below some facts that are used in the paper.

 The limit set $\Lambda (G)$ is the set of accumulation (limit) points of
 the $G$-orbit for the action of $G$ on $T.$
It is known that either $\vert\La(G)\vert\in\{0,1,2 \}$ in which
case the action $G\act T$ is called {\it elementary } or it is a
perfect set and the action is not elementary \cite{Tu2}.

  An elementary action of a group  on $T$ is called  {\it parabolic}
if there is unique fixed point    called parabolic fixed point.

A limit point $x\in \La(G)$ is called {\it conical} if there
exists an infinite sequence $g_n\in G$ and distinct points $a,
b\in T$ such that
$$\forall y\in T\setminus\{x\}\ :\ g_n(y)\to a\ \wedge g_n(x)\to
b.$$

A parabolic fixed point $p\in \La(G)$ is called {\it bounded
parabolic} if the quotient space\hfil\penalty-10000 $(\La(G)\setminus \{p\})/{\rm
Stab}_G p$ is compact.

A set $M$  is called {\it $G$-finite } if $M/G$ is a finite set.

 An action  of a group $G$ on a compactum  $T$  is called {\it
geometrically finite} if every limit point of $T$ is either
conical or bounded parabolic.  As we have pointed out in the
Introduction if   $G\act T$ is a $3$-discontinuous and
$2$-cocompact action then it is also a geometrically finite one.
The opposite statement is also true if one assumes that $T$ is
metrizable.
\medskip

\noindent {\bf Notation.} From now on we  fix the notation $\P$
for the set of parabolic points for the $3$-discontinuous and
$2$-cocompact action $G\act T.$

\bigskip

\section{Exhaustion of non-finitely generated relatively hyperbolic groups
by finitely generated ones.}

\subsection{Entourages, shadows, betweenness relation.}

The following definition is motivated by \cite{Bourb} and
\cite{W}.

\begin{dfn}
\label{ent}   {\rm  Let $T$ be a compactum. Any (not necessarily
open) neighborhood of the diagonal $\bde2$ in $\bs2$ is called
{\it entourage} of $T.$ The set of all entourages of $T$ is
denoted by $\ent$.}
\end{dfn}

 \noindent {\bf Convention.}  By definition an entourage   consists of {\bf non-ordered}
 pairs.  However     sometimes we identify an  entourage
 $\be\in\ent$
 with the
 {\bf symmetric} neighborhood $\ti\be$ of the diagonal in $T\times T$.

\medskip

\noindent   We  denote the entourages   by   bold small
characters.

\medskip

An entourage $\bf e$ determines a graph whose vertex set is  $T,$
and two vertices $x, y$ are joined by an edge if and only if $\{x,
y\}\in \bf e$. Denote by $\de$ the corresponding graph distance
which is the maximal distance function with the property
$\{x,y\}{\in}\mathbf e\implies\boldsymbol\Delta_{\mathbf
e}(x,y){\leqslant}1.$  Note that $\de(x,y)=\infty$ if and only if
$x$ and $y$ belong to different connected components of the graph.
A set $U\subset T$ is called ${\bf e}$-small if its $\be$-diameter
is at most $1.$

 The set of all $\bf e$-small sets is denoted  by Small$(\bf e)$. For
 subsets $a, b\subset T$ we define $\de(a, b)={\rm inf}\{\de(x,y)\
\vert\  x\in a,\ y\in b\}$ and $\ti\de(a, b)={\rm sup}\{\de(x,y)\
\vert\  x\in a,\ y\in b\}$. From the triangle inequality we have
the inequality
 $\de(a, b)\geq \de(a, c) -\ti\de(c,b)$ frequently used further.

  For a subset $a\subset T$ define its $\bf
e$-neighborhood $a\bf e$ as $\{x\in T\ \vert\ \Delta_{\mathbf
e}(x,a){\leqslant}1 \}$.

For a  subset $o$ of $T$   its "{\it convex hull}" in
 $T{\sqcup}\mathsf{Ent}T$ is the set $$\widetilde o=o{\cup}\{\mathbf
e{\in}\mathsf{Ent}T:o'{\in}\mathsf{Small}(\mathbf
e)\},\hfill\eqno(\dagger)$$ \noindent where $o'$ denotes the
complement of $o$.

 We equip the space
$T{\sqcup}\mathsf{Ent}T$ with the topology generated by the "{\it
convex hulls}" of open subsets of $T$ and the single-point subsets
of $\ent.$ Namely a set $w$ in $T{\sqcup}\mathsf{Ent}T$ is
declared {\it open} if
 for every point $t\in w\cap T$
 there exists  and open subset $o$ of $T$ such that $t\in o$ and $\widetilde o\subset w.$  In particular $\ent$ is a
discrete open subset  and $T$ is a closed subspace of
$T{\sqcup}\mathsf{Ent}T$.

\medskip

\noindent {\bf Example 1.}   The definition of the topology on
$T{\sqcup}\ent$ can be illustrated in  terms of the open subsets
of the compactified real hyperbolic space
$\hn\cup\partial_{\infty}\hn$.
  Let
$B$ be a bounded subset of $\hn$. Define an entourage
$\be_B\in{\rm Ent }(\partial\hn)$ in the following way: $\{x,
y\}\in\be_B$ if and only if the geodesic $\ga(x,y)$ with the
endpoints $x$ and $y$ misses   $B.$  So
  a set $o\subset T=\partial_{\infty}\hn$ is $\be_B$-small if and
only if $B$ is contained in the convex hull of $o$ in $\hn$ (see
Figure 1).  Thus $B$
 is close to $a$ in the topology of $\hn\cup\partial_{\infty}\hn$
if and only if $\be_B$ is close to $a$ in the topology of
$T{\sqcup}\ent$. By the above definition   $\widetilde o$ is
obtained by adding to $o$ every entourage
 for which $o'$ is small. \bx

\centerline{
\ifnum\pdfoutput>0%
\pdfximage width 180pt height 200pt{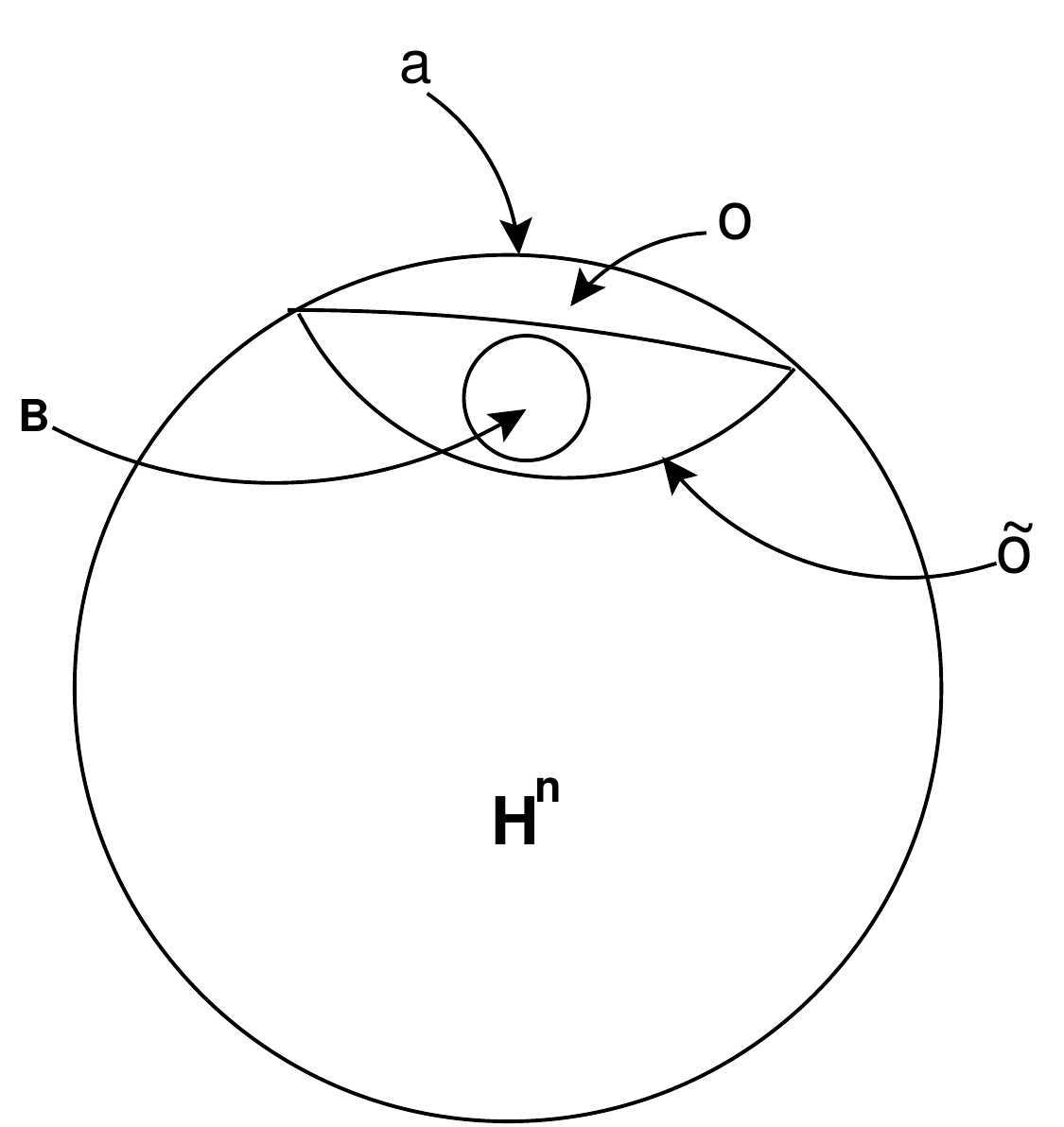}\pdfrefximage\pdflastximage\else\kern180pt\vbox to 200pt{\vss}\fi}

\begin{figure}[sh] \leavevmode

\caption{Bounded set  in $\hn$ and its visibility entourage.}
\end{figure}

\bigskip

\medskip

\begin{dfn}
\label{link} {\rm \cite{Ge1}
 Two entourages $\bf a$ and $\bf b$ are said to be {\it unlinked} if there exist
  $a\in {\rm Small}(\ba)$ and $b\in {\rm Small}(\bb)$
  such that $T=a\cup b$. We denote  this relation by $\bf a\bowtie\bf b$.
  In the opposite case
  we say that  $\bf a$ and $\bf b$ are {\it linked},
  and write $\bf a\# \bf b$.}
  \bx
\end{dfn}

\noindent Denote by $\lba$  the set $\{\bb\in\ent\ \vert\
\ba\#\bb\}.$ It is enough for our purposes to consider only
sufficiently small entourages implying the following.

\medskip

\noindent {\bf Convention.} All considered entourages are supposed
to be self-linked :

$$\ba\in\ent\ :\ \ba\#\ba. \eqno(1)$$

\bigskip

\begin{dfn}
\label{shadow}

{\rm
  \cite{Ge1}  Let $\bf a$ and $\bf b$ be two unlinked entourages. We define the following "shadow" sets :

$$\Sh_{\bf a} \bb=\{ a\in {\rm Small}(\ba)\ \vert\ a'\in {\rm Small}(\bb)\},$$

and

$$\sh_{\bf a}{\bf b}=\bigcap {\rm Sh}_{\bf a} {\bf b}= (\bigcup {\rm Sh}_{\bf b} {\bf
a})'.$$}

\end{dfn}

\medskip

 It is shown in  [Ge1, Lemma S0] that if
  $\ba\bowtie\bb$ and ${\rm diam}_\ba T > 2$ then
$\shab\not=\emptyset$; and if ${\rm diam}_\ba T > 4$ then $\shab$
has a nonempty interior.

\bigskip
 \noindent {\bf Convention.} We consider only the entourages
$\mathbf a$ with $\mathsf{diam}_{\mathbf a}T{>}4$. So every shadow
has non-empty interior.

\bigskip

\noindent {\bf Example 2.} \noindent Using the notations of
Example 1 let $\ba{=}\be_A$ and $\bb{=}\be_B$ for two disjoint
balls $A$ and $B$ in the hyperbolic space $\hn$. Then the shadow
$\sh_{\bf a}{\bf b}$ is given by the intersection  with
$\partial\hn$ of the boundaries of all hyperbolic half-spaces of
$\hn$ containing $B$ and not containing $A$ and similarly for
$\sh_{\bf b}{\bf a}$ (see Figure 2).

\centerline{
\ifnum\pdfoutput>0%
\pdfximage width 324pt height 166pt{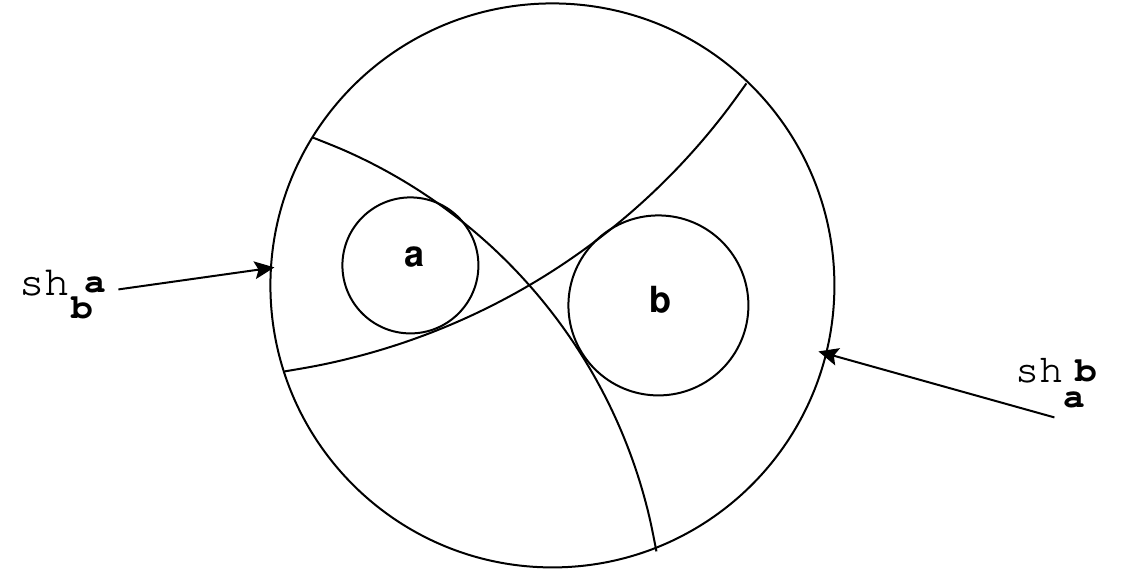}\pdfrefximage\pdflastximage\else\kern324pt\vbox to 166pt{\vss}\fi}

\begin{figure}[sh]
\leavevmode
\caption{Shadows $\sh_{\bb}\ba$ and $\sh_{\bb}\ba.$}
\end{figure}

\begin{dfn} (\it Betweenness relation).
\label{beet} Let $k$ be a  positive integer. \label{between} {\rm
\begin{itemize}

  \item[\sf 1)] Suppose $\ba, \bb, \bc\in\ent.$ We say that an
   entourage $\bf b$ lies
 {\it between} (or   {\it $k$-between})   ${\bf a}$ and  ${\bf
 c}$, and write $\ba - \bb - \bc\ (k)$   (or simply $\ba - \bb - \bc$),
if $\ba\bowtie\bb\bowtie\bc$ and $\db (\shba, \shbc) > k.$

\medskip

  \item[\sf 2)] Let $\mathbf{a,b}{\in}\mathsf{Ent}T$ and  let $p\in T$.
  We say that $\mathbf b$ lies {\it between} (or $k${\it-between}) $\mathbf a$ and $ p$
if $\mathbf{a{\bowtie}b}$ and  $\db(\mathsf{sh}_{\mathbf b}\mathbf
a,b){>}k$ for any $\mathbf b$-small neighborhood $b$ of $ p$

We write $\ba - \bb - p\ (k)$   (or simply $\ba - \bb - p$) in
this case.
\medskip

  \item[\sf 3)] Let $\mathbf{b}{\in}\mathsf{Ent}T$ and  let $p, q\in T$ be two
   distinct points.
  We say that $\mathbf b$ lies {\it between} (or $k${\it-between}) $p$ and
  $q$, and write $q - \bb - p\ (k)$   (or simply $q - \bb - p$),
if $\db(b_1, b_2)
> k$ for any $\mathbf
b$-small   neighborhoods $b_1$ and $b_2$ of the points $ p$ and
$q$ respectively.

\medskip

\end{itemize}

}
\end{dfn}

\medskip

\begin{rems}
\label{betw} {\rm a)  The betweenness relations 2) and 3)
represent an  extension "by continuity"   of the   relation  1)
between entourages to the points of $T$. Note that the middle
object in the relation $\ba-\bb-\bc$ is always an entourage.

Note also that if $\db(\mathsf{sh}_{\mathbf b}\mathbf a,b_0){>}k$
for some $\bb$-small neighborhood $b_0$ of $p$ then for any such
$b$ we have $\db(\mathsf{sh}_{\mathbf b}\mathbf a,b)\geq
 \db(\mathsf{sh}_{\mathbf b}\mathbf a,b_0)-\db(b,b_0)> k-2$ as
 $p\in b\cap b_0$ and  $\ti\Delta_{\bb}(b,b_0)\leq 2$. Therefore we will
 always assume further that $k>2.$

\medskip

b) Definition \ref{beet} in cases 2) and 3)    differs  from the
corresponding definition  in [Ge1] where the condition
$\Delta_\bb(\shab,p)> k$ is stated instead of 2). The above
betweenness definition is stronger than that of \cite{Ge1} and so
is easier to use. However both of them are quite close:  the
$k$-betweenness 2) implies $k$-betweenness of \cite{Ge1}. On the
other hand since the diameter of any small neighborhood is less
than $1$ the $k+1$-betweenness of \cite{Ge1} implies (by the
triangle inequality) the $k$-betweenness 2). We will use results
of \cite{Ge1} keeping in mind this relation.}
\end{rems}

\medskip

\begin{lem} \label{approximation
points}({\it Continuity property}). Suppose that ${\bf a} - {\bf
c} - p(k)\ (k\in\N)$ where $\ba\in T{\sqcup}\ent,\ \bc\in \ent,\
p\in T.$ Suppose that $p\in T$ is an accumulation point for an
infinite subset $B$ of  $\ent$. Then there exists $\bb\in B$ such
that ${\bf a} - {\bf c} - {{\bf b}} (k).$

\end{lem}

 \proof Let  first $\ba\in \ent$ be an entourage. Let  $c=U_p$ be an
open $\bc$-small set containing $p$ such that $\dc(c, \shca) >k$.
 By definition of the topology of
$T{\sqcup}\ent$ the complement $c'$ is $\bb$-small for some
$\bb\in B$. Then $c'\subset \bigcup\Sh_{\bb}\bc,$ and $c\supset
\sh_{\bc}\bb=(\bigcup\Sh_{\bb}\bc)'.$ Thus $\dc (\sh_{\bc}\ba,
\sh_{\bc}\bb)> \dc (\shca, c)   > k.$

If now $a\in T$ then for a $\bc$-small neighborhood $U$ containing
$a,$ we  obtain similarly $\dc(U, \sh_{\bc}\bb) > \dc (U,U_p)   >
k$. So we still have  $a-\bc-\bb\ (k)$ for $\bb\in B.$ \bx

\bigskip

 \noindent
\begin{dfn} \label{tubes}({\it Tubes}). \cite{Ge1} {\rm A
sequence $P$ of elements ${\ba_n}$  of $T\sqcup\ent$ is called
$k$-{\it tube} (or  {\it tube}) if

 $$\forall n:\ ({\ba}_n\bowtie  {\ba}_{n+1})\ {\wedge}\ ({\ba}_{n-1} - \ba_n - {\ba}_{n+1}(k))$$

\noindent  whenever ${\bf a}_{n\pm1}$ are defined.}
  \bx

\end{dfn}

\begin{lem}
\label{order}

\begin{itemize}

  \item[\sf 1)] (Ordering) For any three entourages at most one can be between
   the others.

\medskip

\item[\sf 2)] (Convexity) If $\ba-\bb-\bc (4)$ and $\ba,\bc\in
\ld$ then $\bb\in  \ld.$
\end{itemize}
\end{lem}

\proof 1) Indeed if not, we obtain $\ba-\bb-\bc$ and $\ba-\bc-\bb$
for some $\ba, \bb, \bc.$ The transitivity of the betweenness
relation \cite{Ge1} would imply  $\ba-\bb-\ba$ and so
$\ba\bowtie\ba$ which is impossible by our convention (1).

\medskip

\noindent 2) Otherwise $\bb\bowtie\bd$ and we have $T=b\cup d=
a\cup b_1=c\cup b_2$ where $b_i,b \in {\rm Small}(\bb)\ (i=1,2),\
d\in {\rm Small}(\bd), a\in {\rm Small}(\ba), c\in {\rm
Small}(\bc)$. It follows that $b\cap b_1=\emptyset$ or $b\cap
b_2=\emptyset$ since otherwise $\Delta_{\bb}(b_1, b_2)\leq 2$ and
we would have $\Delta_{\bb}(\sh_{\bb}\ba, \sh_{\bb}\bc)\leq
2+\ti\Delta_{\bb}(\shba, b_1)+\ti\De_{\bb}(\shbc, b_2)\leq 4$ (as
$\shba\subset b_1$ and $\shbc\subset b_2$) which is impossible.
If, for instance, $b\cap b_1=\emptyset$ then $b_1\subset d$ and
$\ba\bowtie\bd.$ A contradiction. \bx

\subsection{Discrete sets of entourages. Horospheres.}

Until the end of Section 3 we fix a 3-discontinuous 2-cocompact
action $\GAT$ of a   group $G$ on a compactum $T$.

\begin{dfn}
\label{discr}{\rm
   A  set $A$  of entourages on $T$ is called {\it discrete}
   if

   $$\forall {\bf w}\in\ent\ :\ \vert\{\ba\in A\ :\ \ba\#
   {\bf w}\}\vert<\infty.\hfill\eqno(1)$$}

   \end{dfn}

\noindent By    [Ge1, Proposition P] the set $\{g\in G\ :\ g\ba\#
\bw\}$ is finite for all $\bw,\ba\in \ent$. This  property is
called {\it Dynkin property} \cite{Fu}. Hence every $G$-finite set
is discrete.
\medskip

Let  $A\subset\ent$ be  a $G$-finite set of entourages. Denote by
$\ti T$ the subspace  $T{\sqcup}A$ of $T{\sqcup}\ent.$ Since  $A$
is discrete $\ti T$ is compact [Ge1, Proposition D].
\medskip

\begin{dfn}
\label{ggraph}  {\rm    Let $\Ga=\Ga_A$ be the graph whose vertex
set $\Ga^{0}$ is $A$ and the edge set $\Ga^{1}$ is the set of
pairs $\{\ba, \bb\}$ such that $\ba\#\bb.$ Denote by $d_A$ the
corresponding graph distance.}
\end{dfn}

Since $G$ acts on $T$ by homeomorphisms it acts isometrically on
$(\Ga, d_A).$

\medskip
\begin{lem}
\label{connect} The group $G$ is finitely generated if and only if
there exists a connected graph $\Ga_A.$
\end{lem}

\proof  Suppose first that $G$ admits a finite set of generators
$S\ ({\rm id}\in S$). Since $A$ is $G$-finite  we have
$\displaystyle A=\bigcup _{i=1}^l G(\ba_i).$ Any entourage $\ba_i$
contains a sub-entourage $\ba'_i$ such that $$\forall s\in S\ \ :\
\ba'_i\# s\ba'_j\ \ (i,j\in\{1,...,l\}).$$

 So up to choosing the
entourages $\ba_i\ (i=1,...,l)$ to be sufficiently small we can
assume that the above property is satisfied. Then all vertices in
the set  $\displaystyle \bigcup_i S\ba_i$ are pairwise connected
by edges. For any vertex $\bv\in \Ga_A$ there exists
$i\in\{1,...,l\}$ and $g\in G$ such that $\bv=g(\ba_i)$ and
$g=s_{i_1}s_{i_2}...s_{i_k}\ (s_{i_j}\in S)$. Then $\Ga_A$
contains  the edges $e=(s_{i_k}(\ba_i),\ba_i)$,
$e'=(s_{i_{k-1}}(\ba_i), \ba_i),$
 and so the path $s_{i_{k-1}}e\cup e'$
 between
$\ba_i$ and $s_{i_{k-1}}s_{i_k}(\ba_i)$. Continuing in this way we
obtain a path between $\bv$ and $\ba_i.$

Conversely suppose that $\Ga_A$ is connected.  Let $S$ be the set
$\{s\in G\ \vert\ s\ba_j\# \ba_i,\ 1\leq i\leq l\}$ where
$\displaystyle A=\bigcup_{i=1}^lG\ba_i$. By Dynkin property the
set $S$ is finite. For any $g\in G$ there is a path $l=\{\ba_i,
\bb_2, ...,\bb_{n-1}, \ba\} \subset \Ga_A$ between the vertices
$\ba=g(\ba_i)$ and $\ba_i.$ Then $\bb_2\# \ba_i$ so $\exists\
s_1\in S \ :\ \bb_2=s_1(\ba_j)\ (1\leq j\leq l)$. Thus
$s_1^{-1}\bb_3\#\ba_j$ and $\exists\ s_2\in S\ :\
\bb_3=s_1s_2\ba_k\ (1\leq k\leq l).$ Continuing in this way we
obtain $\ba=s_1s_2...s_n\ba_r\ (1\leq r\leq l).$ Then
$g^{-1}s_1s_2...s_n(\ba_r)=\ba_i$ and so $g^{-1}s_1s_2...s_n$
belongs to $S$ (by (1) of 3.1). The lemma is proved. \bx

\bigskip

It follows from Dynkin property and our convention (1) that the
stabilizer of each edge and each vertex of $\Ga$  is finite. The
action $\GAT$ is 2-cocompact so by    [Ge1, Proposition E] we can
suppose that the set $A$ is a single orbit $G(\ba_0)\ (\ba_0\in
\ent)$ having the following
 properties   :

\begin{itemize}

  \item[\sf i)] {\it $m$-separation property}:

$$\forall (p, q)\in\Th^2T\ \exists\ \ba\in A\ : p-\ba-q (m),\hfill\eqno(2)$$

\noindent for a fixed $m\in\N$.

\medskip

  \item[\sf ii)] {\it generating property:}

   $$ \ \displaystyle\forall~\bu\in\ent\ \exists~ \ba_i\in A\ (i=1,..,l):\
  \bu\supset\bigcap_{i=1}^l\ba_i.\hfill\eqno(3)$$

i.e.  $A$  generates  $\ent$ as a
  filter.
  \end{itemize}
\bigskip

\begin{conv}
\label{setA} {\rm From now on we fix an unlinked entourage
$\ba_0\in\ent$ (see (1) of 3.1) such that its orbit  $A=G(\ba_0)$
  satisfies $m$-separating and generating
properties. The value of $m$ can be easily restored in each
statement. Keeping in mind that this value might be needed to be
increased further we just suppose that $m$ is sufficiently large.

Furthermore if $G$ is finitely generated we will always assume (by
Lemma \ref{connect}) that the graph $\Ga$ is connected.}
\end{conv}

\medskip

\noindent {\bf Remarks.} The  graph  $\Ga$  plays the role
 of the Cayley
graph $\Ca(G)$ if $G$ is finitely generated, however by Dynkin
property it is always a locally finite graph.  The space $\ti
T=T{\sqcup}A$ is a compactification of $A=\Ga^0$ similar to the
Floyd completion (see Section 4). Every action $G\act T$ can be
naturally extended to the space $\ti T$.

\bigskip

\begin{lem}
\label{hausd} The space $\ti T=T{\sqcup}A$ is a compactum.
\end{lem}

\proof The space $T$ is Hausdorff. To prove that $\ti T$ is
Hausdorff we will consider three different cases. Let first $x, y$
be distinct points of $T$ then there exist disjoint closed
neighborhoods $U_x$ and $U_y$ in $T.$    Their convex hulls $\ti
U_x=U_x\cup \{\be\in A :\ U'_x\in \mathsf{Small}(\be)\}$ and $\ti
U_y=U_y\cup \{\bd\in A :\ U'_y\in \mathsf{Small}(\bd)\}$ are
neighborhoods of these points  in the  topology of $\ti T$ induced
from $T{\sqcup}\ent$ (see ($\dagger$) of 3.1). If $\ba\in A\cap\ti
U_x\cap \ti U_y$ then $U'_x$ and $\ti U'_y$ are both $\ba$-small.
Since $U_x$ and $U_y$ are disjoint   we have $U'_x\cup U'_y=T$ and
so $\ba\#\ba$ contradicting our Convention (1) of $\S 3.1$. Hence
$\ti U_x\cap \ti U_y=\emptyset$ in this case.

 If now    $\zx\in A$ and $y\in T$ then by the
same reason any $\zx$-small neighborhood of $y$ in $\ti T$ cannot
contain $\zx.$ Since every entourage is open in $\ti T$ we are
done in this case too. If finally both points are entourages they
coincide with their disjoint neighborhoods. So $\ti T$ is
Hausdorff.

The compactness of $\ti T$ follows from [Ge1, Proposition D]. \bx

\bigskip

\begin{prop}
\label{convT} If a group $G$ acts 3-discontinuously  on a
compactum $T$ then the induced action on $\ti T= T{\cup}A$ is also
3-discontinuous.
\end{prop}

\medskip

\noindent {\bf Remark.}  In [Ge2, Theorem 5.1] it is proved that
there is a unique topology on the compactified space $\ti T$ with
respect to which the action is 3-discontinuous. The argument below
provides a simple proof of this for the induced  topology on $\ti
T\subset T{\sqcup}\ent$  introduced above.

\medskip

\proof For a  subset $X\subset T$  denote by $\ti X= X\cup\{\ba\in
A\ \vert\ X'\in {\rm Small}(\ba)\}\subset\ti T$ its convex hull in
$\ti T$. In case if $X=\{\ba\}$ where $\ba\in A$ is an entourage
we put $\ti X=\ba.$ For every $g\in G$ denote by $\ti g$ its
natural extension to $\ti T.$

Every point   $x\in \Th^3\ti T$ admits a closed neighborhood
 which is a "cube"    $\ti K =
\ti  X\times \ti Y\times \ti Z$ where $X, Y$ and $Z$ are either
disjoint closed subsets of $T$ or some of $\ti X, \ti Y,\ti Z$ are
isolated entourages (in the latter case we call the corresponding
cube {\it degenerate}). Every compact subset of $\Th^3\ti T$ is a
finite union of such cubes. So it is enough to prove that for two
cubes $\ti K_i=\ti X_i\times \ti Y_i \times \ti Z_i\subset
\Theta^3\ti T\ (i=0,1)$ the following set is finite:

$$S=\{g\in G\ \vert \ :\ \ti g\ti X_0\cap \ti X_1\not=\emptyset,\ \ti g\ti Y_0\cap \ti
Y_1\not=\emptyset,\ \ti g\ti Z_0\cap \ti Z_1\not=\emptyset\}.$$

\noindent  Suppose to the contrary that $S$ is infinite. Since the
action $G\act T$ is 3-discontinuous,   every accumulation point
of $S$ with respect to Vietoris topology  is a {\it cross} $<p,
q>^\times=p\times T\sqcup T\times q$  [Ge1, Proposition P].
Consider now all possible  cases.

\medskip

 \noindent {\it Case 1.} Both cubes are not degenerate, i.e.
$X_i, Y_i, Z_i\ (i=0,1)$ are all closed disjoint subsets of $T$.

\medskip

Note that at least one of the "squares" $X_0\times X_1,\ Y_0\times
Y_1$ or $Z_0\times Z_1$ does not meet the cross. Indeed otherwise
two of them intersect both either $p\times T$ or $T\times q$ which
is impossible as $X_i, Y_i$ and $Z_i$ are pairwise disjoint    for
$i\in\{0,1\}$.

Let us assume that e.g. $Z_0\times Z_1\cap <p,
q>^\times=\emptyset$. Let $g\in S$ be  a homeomorphism whose graph
is contained in the neighborhood $T^2\setminus Z_0\times Z_1$ of
$<p,q>^\times$. Then $gZ_0\cap Z_1=\emptyset$. However $\ti g\ti
Z_0\cap \ti Z_1\not=\emptyset.$ So there exists $\ba\in \ti
Z_0\setminus Z_0$ such that $\ti g\ba\in  \ti Z_1$. By definition
of the convex hull $Z'_0$ and $(g^{-1}(Z_1))'$ are $\ba$-small.
Since $(g^{-1} Z_1)'\cup Z'_0=T$ we obtain that $T$ is the union
of two $\ba$-small sets, so $\ba\#\ba$ contradicting our
Convention \ref{setA}.

\medskip

\noindent {\it Case 2.} At least one of the cubes is degenerate.

\medskip

Then some of the sets $\ti X_i, \ti Y_i, \ti Z_i$ are entourages.
Note that since $g\ti X_0\cap \ti X_1\not=\emptyset$ for
infinitely many $g\in S$, by Dynkin property  $\ti X_0$ and $\ti
X_1$ cannot be entourages simultaneously. The same is true for
$\ti Y_i$ and $Z_i\ (i=0,1)$. So there could be at most 3
  entourages among these 6 sets. We   consider
  all the  possibilities below.

\medskip

\noindent {\it Subcase 2.1.} There is only one degenerate cube.

\medskip

We can assume that  $\ti X_0=\ba$ for some $\ba\in A.$ Then
$\forall g\in S$ we have $g\ba\in \ti X_1$. So $g^{-1}X'_1$ is
$\ba$-small. For a limit cross $<p, q>^\times$ for the set $S$ and
$\ba$-small neighborhoods $U_p$ and $U_q$ of the points $p$ and
$q$ respectively there exists   $g\in S$ such that $gU'_p\subset
U_q$ or $g^{-1}U'_q\subset U_p$. If now $U_q\cap X_1=\emptyset$
then $T$ would be the union of $\ba$-small sets $g^{-1}X'_1$ and
$g^{-1}U'_q$  contradicting the unlinkness condition $\ba\#\ba.$
So for every $\ba$-small neighborhood
  $U_q$ of $q$ we have $U_q\cap
X_1\not=\emptyset$. Since $X_1$ is closed it follows that $q\in
X_1$.

At most  one of the disjoint sets $Y_0$ or $Z_0$ can contain the
other
 point $p$ of the cross, let $p\not\in Z_0.$ Then for any neighborhood $U_q$
  and  for infinitely many  elements $g\in S$
we have $gZ_0\subset U_q$. If $gZ_0\cap Z_1\not=\emptyset$ for
infinitely many $g\in S$ then $q$ is an accumulation point for
$Z_1,$ and since $Z_1$ is closed we obtain that $q\in Z_1\cap X_1$
which is impossible. So for almost all $g\in S\ :\ gZ_0\cap
Z_1=\emptyset$ and this situation has been excluded in Case 1.

\medskip

\noindent {\it Subcase 2.2.} There are two degenerate cubes.

\medskip

 Note that they cannot belong to the same level, namely if $\ti X_0=\ba\in A$
and $\ti Y_0=\bb\in A$ then by the  argument of Subcase 2.1 we
must have $q\in Y_1\cap X_1$ which is impossible.

So let $\ti Y_1=\bb\in A$ and $\ti X_0=\ba\in A.$ By the argument
of Subcase 2.1 applied now to the inverse elements of $S$ we
obtain that $p\in Y_0.$ Hence for almost all elements $g\in S$ we
still have $gZ_0\cap Z_1=\emptyset$ which is impossible by Case 1.

\medskip

\noindent {\it Subcase 2.3.} There are three degenerate cubes.

\medskip

Then there are at least two of three entourages which are among of
the sets of the same level: $\ti X_i, \ti Y_i, \ti Z_i\ (i=0\ {\rm
or}\ i=1)$ which is  impossible. So neither case can happen.
 The  proposition is proved. \bx

 \medskip

\bigskip

\begin{lem}
\label{bbsequences} Let $B$   be an infinite subset of $A$ and $C=
N_d(B)$ where $N_d(B)$ is a $d$-neighborhood of $B$ in $\ti T$.
Then the  topological boundaries of $B$ and $C$ coincide.

In particular, if $(\bb_n)_n$ and $(\bc_n)_n$ are two sequences in
$A$ such that $d_A(\bb_n, \bc_n)$ is uniformly bounded, then
$(\bb_n)_n$ converges to a point $p\in T$ if and only if $\bc_n\to
p$.
\end{lem}

\proof The second claim directly follows from the first one. So to
prove the lemma we need only to show that every accumulation point
of $C$ is also an accumulation point of $B.$ Suppose not and there
exists a point $r\in\partial C\setminus
\partial B.$ Then for every neighborhood $U_r$ of $r$ in $\ti T$ there
exists an infinite subset $C_0\subset C$ such that $\forall \bc\in
C_0$ we have $\bc\in U_r$ implying that $U'_r\subset c$ for some
$c\in {\rm Small}(\bc).$

Arguing by induction on $d$ without loss of generality we may
assume that $d=1$. So $\forall \bc\in C\ \exists\ \bb\in B\ :\
\bc\# \bb.$ Then there exists a subset $B_0\subset B$ such that
$d_A(B_0, C_0) \leq 1.$ Since $C_0$ is infinite by discreteness of
$A$ the set $B_0$ is infinite too. Let $p\in T\setminus \{r\}$ be
an accumulation point of $B_0$. Then  for every neighborhood $U_p$
of $p$ there exists $\bb\in B_0$, corresponding to some $\bc\in
C_0$,  for which $U'_p\subset b$ where $b\in {\rm Small}( \bb).$
Choosing $U_p$ to be disjoint from  $U_r$   we obtain $b\cup c=T$
and so $\bb\bowtie\bc.$ A contradiction.\bx

\begin{dfn} \cite{Ge1} ({\it Horospheres, Conical and Parabolic Points}). {\rm Let
$k$  be a fixed positive integer, and let  $A$ be the above
discrete set of entourages.

 \label{hor}
\begin{itemize}

  \item[\sf 1)] We say that  a point $p\in T$  and an entourage   $\be$ are
  {\it neighbors} (with respect to  $A$) and
  write   ${\bf e} \underset{A, k}\#
  p$,
  if there is no
$\ba\in A$ such that ${\bf e} - {\bf a} - p(k)$.

\item[\sf 2)] The {\it horosphere} $T_{A, k} (p)$ (or $T_k(p)$ or
$T(p)$) at the point $p\in T$ is the set
$$T_{A, k}(p)=\{{\bf e}\in A ~\vert ~{\bf e}\underset{A, k}\# p\}.$$

\medskip

\item[\sf 3)] A point $x\in T$ is called {\it $(A, k)$-conical}
(or just {\it conical}) if $T_{A, k} (x)=\emptyset.$

\medskip

\item[\sf 4)]  A point $p\in T$ is called  $(A, k)$-{\it
parabolic} (or just {\it parabolic}) if   $T_{A, k} (p)$ is
infinite.
\end{itemize}
}

\end{dfn}

\noindent It is shown in \cite{Ge1} that the notions of $(A,
k)$-conical and $(A, k)$-parabolic points for $k\geq 3$ (see also
Remark \ref{betw}) are equivalent to the standard definitions (see
Section 2) of conical
 and bounded parabolic points respectively.

\medskip

 \begin{lem}
 \label{geomfin}\cite{Ge1} If the action $G\act T$ is
 $3$-discontinuous and $2$-cocompact then every  limit point of this action is either
conical or bounded parabolic. Furthermore the set of  non-conical
points is $G$-finite  and for every parabolic point $p\in T$ the
set $T(p)$ is ${\rm Stab}_Gp$-finite.
\end{lem}

\bigskip

The next lemma is proved in [Ge1, Lemma P2] for closed entourages.
We prove it below in a general form.

\begin{lem}
\label{unicityacp} For every $d>0$ the   parabolic point $p$ is
the unique accumulation point of the $d$-neighborhood $N_d(T_{A,
k}(p))$ of the horosphere $T_{A, k}(p)$.
\end{lem}

\proof  By Lemma \ref{bbsequences} it is enough to prove the
statement for the horosphere $T_{A, k}(p).$ Suppose it admits two
distinct accumulation points $p$ and $q.$ Since the set  $A$ is
$m$-separating there exists $\ba\in A$ such that $p-\ba-q (k)$ for
some $k \leq m.$ Then by Lemma \ref{approximation points} there
exists $\bb\in T_{A, k}(p)$ such that  $p-\ba-\bb(k)$ which is not
possible. \bx

\medskip

\bigskip

\noindent We have the following transitivity property:

\begin{lem}
\label{busorder} If $\ba, \bb, \bc\in \ent,\ p\in T$ and $k > 2$.
Then   $\ \ba-\bb-p(k)$ and $\bb-\bc-p(k)$ imply $\ba-\bc-p(k).$
\end{lem}

\proof If $a\in\Shab$ and $c\in \Sh_{\bc}\bb$, then the sets
$b=a',b_1=c'$ are $\bb$-small and $a\cup b=b_1\cup c=T.$ Then for
a $\bc$-small neighborhood  $c_0$ of $p$ we have
$\Delta_{\bc}(c,c_0)\geq \Delta_{\bc} (\shcb,
c_0)-\ti\Delta_{\bc}(\shcb, c) > k-1.$ So $\Delta_{\bc}(c,p) > k-1
> 0$ and $p\in b_1.$ Note that $b\cap b_1=\emptyset$ since
otherwise $\db(b_1,\shba)\leq \ti\db(b_1, b\cap b_1) + \db(b\cap
b_1, \shba)\leq 2$ which is impossible as $\ba-\bb-p(k)$ and
$k\geq 2.$ Thus $b_1\subset a$ and $a\cup c =T$. Since $c$ was an
arbitrary element of $\Sh_{\bc}\bb,$   it follows that
$\Shcb\subset \Shca$ and $\sh_{\bc}\ba\subset \sh_{\bc}\bb$. Thus
$\dc(\sh_{\bc}\ba, c_0)> k$.\bx

\bigskip

\noindent The above notions allow us to introduce the following
relation on the set $\ent$.

\begin{dfn}
\label{busord} ({\bf Busemann order}) {\rm For $\ba,\bb\in\ent,\
{\rm and}\ p\in T$ we say that $\ba$ and $\bb$ are {\it Busemann
ordered} with respect to $p$ if

$${\rm either}\ \ba=\bb,\ \ {\rm or}\ \ \ba-\bb-p(k).$$

\noindent We will denote this relation by $\ba\geq_{p,k}\bb.$}
\end{dfn}

 Lemma \ref{busorder} implies that this relation  is a partial
order on $\ent$. Using Busemann order we can reformulate the above
definitions of conical and parabolic points as follows.

\begin{lem}
\label{busem} A point $p\in T$ is $A$-conical if and only if its
Busemann order has no minimal elements. A point $p$ is
$A$-parabolic   if and only if its Busemann order has infinitely
many minimal elements.
\end{lem}

\subsection{Non-refinable tubes.}

\medskip

\begin{lem}
\label{finlab} The set $\Psi_k(\ba, \bb)=\{\bc\in A\ :\ \ba -\bc
-\bb(k)\}$ is finite for any $k\geq 1.$
\end{lem}

\proof Suppose that $\ba-\bc-\bb(k)$ and let us  prove that
$\bc\#(\ba\cap\bb).$ If it is not true, then we have
$\bc\bowtie(\ba\cap\bb)$, i.e. there exists $c\in\smc, \ w\in{\rm
Small}(\ba\cap\bb)$ such that $c\cup w =T.$ Thus $c\in
\Sh_{\bc}\ba\cap\Sh_{\bc}\bb$ and $\shca\subset c,\ \shcb\subset
c.$ Hence $\Delta_{\bc}(\sh_{\bc}\ba, \sh_{\bc}\bb)\leq 1$ which
is impossible. It follows that $\bc\#(\ba\cap\bb)$. The finiteness
of $\Psi_k(\ba, \bb)$ now follows from the discreteness of $A$ \bx

\medskip

\begin{dfn}
\label{nonsep} {\rm ({\bf Refinability}). A pair $\{\ba,
\bb\}\subset A$ is called ({$k$\it-)refinable} if $\Psi_k(\ba,
\bb)\not=\emptyset,$ and {\it ($k$-)non-refinable} otherwise.}

\end{dfn}

Proposition \ref{exent} below guarantees the existence of a finite
non-refinable tube between two given entourages in $A$.    To
prove it we need  the following:

\medskip

\begin{lem}
\label{nonref1} For every integer  $k\geq 2$, every pair $\{\ba,
\bb\}\subset A$ is either $k+1$-nonrefinable or there exists
$\bc\in \Psi_k(\ba, \bb)$ such that the pair $\{\ba, \bc\}$ is
$k+1$-nonrefinable.
\end{lem}
\medskip

\proof Suppose this is not true and let a pair $\{\ba, \bb\}$ be a
counter-example. By Lemma \ref{finlab} the set $ \Psi_k(\ba, \bb)
$ is finite so we can assume in addition that the number
$\vert\Psi_k(\ba, \bb)\vert$ is the minimal one among all such
counter-examples. So $\{\ba, \bb\}$ is $k+1$-refinable and there
exists $\bc\in\Psi_{k+1}(\ba, \bb)$ such that the pair $(\ba,
\bc)$ is $k+1$-refinable too. We now claim that

$$\Psi_{k+1}(\ba, \bc)\subset \Psi_{k+1}(\ba, \bb)\ (k>1).\hfill\eqno(1)$$

Let $\bd\in \Psi_{k+1}(\ba, \bc)$. By [Ge1, Lemma T2] we have
$\bd-\bc-\bb(k)$. Then    $\shdb\subset \shdc$ [Ge2, Lemma B1].
Therefore $\Delta_{\bd}(\shdb, \shda)\geq \Delta_\bd(\shdc,
\shda)$. So $\bd\in\Psi_{k+1}(\ba, \bb)$ and (1) follows.

As $\bc\in\Psi_k(\ba, \bb)\setminus \Psi_k(\ba, \bc)$ we obtain
that $\vert\Psi_k(\ba, \bc)\vert < \vert\Psi_k(\ba, \bb)\vert$.
Thus
 by the minimality of $(\ba,\bb)$ the pair $(\ba,\bc)$
cannot be a counter-example.   Then   $(\ba, \bd)$ is
$(k+1)$-nonrefinable. Since $\bd\in\Psi_{k+1}(\ba, \bb)\subset
\Psi_k(\ba, \bb)$ the pair $(\ba, \bb)$ cannot  be a
counter-example neither. A contradiction.\bx

\bigskip

\noindent For a tube  $P=\ba - \ba_1-...-\ba_n-\bb$   we   denote
by $\partial P$ its boundary $\{\ba, \bb\}.$

\bigskip

\begin{prop}
\label{exent}

For every pair $\{\ba, \bb\}\subset A$ and integer $k\geq 2$ there
exists a finite $k+2$-nonrefinable k-tube $P\subset A$ such that
$\partial P=\{\ba, \bb\}.$
\end{prop}

\proof Suppose this is not true. Let a pair  $\{\ba,\bb\}$ be a
 counter-example such that  it has the   minimal
cardinality $\vert\Psi_k(\ba, \bb)\vert$ among all such pairs.
Since  $\{\ba,\bb\}$ is $k+2$-refinable by the above lemma there
exists $\bc\in \Psi_{k+1}(\ba, \bb)$ such that  $\{\ba, \bc\}$ is
$k+2$-nonrefinable. Since the inclusion $\Psi_k(\bc, \bb)\subset
\Psi_k(\ba, \bb)$ is strict there   exists a $k+2$-nonrefinable
$k$-tube $Q$ with $\partial Q=\{\bc, \bb\}.$ By the transitivity
property  [Ge1, Lemma T2] the set $R=\{\ba\}\cup Q$ is a $k$-tube
with the boundary $\{\ba, \bb\}.$ It is $k+2$-nonrefinable by
construction. Thus the pair $\{\ba, \bb\}$ is not a
counterexample. We have a contradiction. \bx

\bigskip

\begin{dfn}\label{horproj}\cite{Ge1} (Horospherical projection). Let $p\in \P$ be a parabolic point and $T(p)$ be a
horosphere at $p$. Define a projection map $\Pi_{p} : A\to T(p)$
(or $\Pi_{p,k}$)  called {\it horospherical projection} as
follows. If $\ba\not\in T_k(p)$ then $\Pi_p(\ba)= \{\bp\in T_k(p)\
:\ \ba-\bp-p(k)\};$ and  if $\ba\in T_k(p)$ then $\Pi_p(\ba)=\ba
.$
\end{dfn}

\begin{prop}
\label{boundness} Let $\P$ denote the set of  parabolic points for
the action $G\act T.$ Then for any constants $k>3$ and $d>0$  the
following sets are
 $G$-finite:

 \begin{itemize}

 \item[\sf 1)]  $\forall \{\ba, \bb\}\subset A\ :\ \{ \{\bc, \bd\}\
 \vert\ \bc\in \Pi_p(g\ba),\ \bd\in \Pi_p(g\bb),\ p\in\P, g\in G\}$

\medskip

 \item[\sf 2)]  $ \A_1=\{(\ba,\bb)\ \vert\ \Psi_k(\ba,
  \bb)=\emptyset,\ \{\ba,\bb\}\not\subset T_{A,k}(p),\ p\in \P\}.$

\medskip

\item[\sf 3)] a) $\{\{p, q\}\subset\ \P \ \vert \
N_d(T_{A,k}(p))\cap N_d(T_{A,K}(q))\not=\emptyset\},$ and

 b) $\{N_d(T_{A,k}(p))\cap N_d(T_{A,k}(q))\ \vert\
\{p,q\}\subset\P\}.$

 \end{itemize}
 \end{prop}

\proof 1)  Suppose to the contrary that the set 1) is infinite.
Assume first that $\ba\not=\bb.$ Then there exist an infinite
sequence of elements $g_n\in G$, distinct entourages $\{\bc_n,
\bd_n\}\subset A$  such that

$$g_n\ba-\bc_n-p_n(k)\ {\rm and}\ g_n\bb-\bd_n-p_n(k),\ \bc_n\in
T_{A,k}(p_n),\ \bd_n\in T_{A,k}(p_n),\ p_n\in\P\hfill\eqno(2).$$

Since the set $\P$ is $G$-finite (Lemma \ref{geomfin})  we can
assume that $p_n=p.$ Since the stabilizer ${\rm Stab}_G p$ acts
cofinitely on $T_{A,k}(p)$ (Lemma \ref{geomfin}) we can also fix
$\bc_n=\bc\in T_{A,k}(p)$, and assume that $\bd_n=h_n(\bd),\
\bd\in T_{A,k}(p), h_n\in {\rm Stab}_G p.$ So (2) gives

$$g_n\ba-\bc-p(k),\ \ g_n\bb-\bd_n-p(k),\ \bc\in
T_{A,k}(p),\ \bd_n\in T_{A,k}(p),\ p\in\P.\hfill\eqno(2')$$

\noindent The following lemma implies that $p$ is a limit point of
$\{g_n\bb\}_n.$

\begin{lem}
\label{smconv} If $\ \bb_n-\bd_n-p(k)\ (k>1),\ \bd_n\in
T_{A,k}(p)$ and $\displaystyle \lim_{n\to\infty}\bd_n= p$ then
$\displaystyle \lim_{n\to\infty}\bb_n= p$.
\end{lem}

\proof We start with  the following.

 \bigskip

 \noindent {\bf Claim.} {\it For every $k>1$ there exists $\bd\in T_{A,k}(p)$ such that
$q-\bd-p(k).$}

\medskip

Indeed by $m$-separation property (2) there exists $\ba\in A$ such
that $q-\ba-p(k)\ 1 < k \leq  m).$ If $\ba\in T_{A,k}$ we are
done. If not let $\bp\in \Pi_p(\ba)$ so $\ba-\bp-p(k).$ Let $U_p$
be    a $\bp$-small neighborhood of $p$. Let also $b\in
\Sh_{\bp}\ba$, then  $a\cup b=T$  where $b$ is $\bp$-small and
$a=b'$ is $\ba$-small set respectively. We have
  $\Delta_{\bp}(b, U_p)
\geq \Delta_{\bp}(U_p, \sh_{\bp}\ba) - \ti\Delta(\sh_{\bp}\ba, b)
> k-1.$ Therefore $U_p\subset
a$ and so $U_p$ is $\ba$-small. Then for any $\ba$-small
neighborhood $U_q$ of $q$  we have $\Delta_{\ba}(U_q, U_p)
>k.$ Hence $\Delta_{\ba} (U_q, a) > k-1$. We have proved that
 $U_q\subset b$ for any $b\in\Sh_{\bp}\ba$. Thus $U_q$ is $\bp$-small and
$U_q\subset \sh_{\bp}\ba$. It implies that   $\Delta_{\bp}(U_q,
U_p)\geq\Delta_{\bp}( \sh_{\bp}\ba, U_p)
>k.$ The claim follows.

\bigskip

 {\it Proof of the lemma.} Suppose by contradiction that
there exists an accumulation point $q\in T$ of the set
$\{\bb_n\}_n$ distinct from $p.$ By the claim there exists $\bp\in
T_{A,k-1}(p)$ such that  $q-\bp-p(k-1)$. Since $\bd_n\to p$ and
$\bb_n\to q$  by Lemma \ref{approximation points} we obtain
$\bb_n-\bp-\bd_n(k-1)\ (n>n_0)$. So
$\sh_{\bd_n}\bp\supset\sh_{\bd_n}\bb_n.$ Since $\bb_n-\bd_n-p(k)$
we obtain $\bp-\bd_n-p(k-1).$ This is impossible as $\bp\in
T_{A,k-1}(p).$ The lemma  is proved. \bx

\bigskip

It follows from $(2')$ that for any $(g_n\ba)$-small set
 $a_n\in\Sh_{\bc}(g_n\ba)$ and a $\bc$-small neighborhood $U_{p}$
 of $p$ we have $\Delta_{\bd_n}(a_n, U_{p}) > k-1>0\ (n\in\N)$.
 Thus $U_p\subset a'_n$ and $U_p$ is
 $g_n\ba$-small for all $n\in \N.$

From the other hand by Proposition \ref{convT} we have that $G\act
\ti T$
 is a convergence action. Then by [GePo1, Lemma 5.1] for every pair of
 distinct non-conical points $\{x, y\}\subset\ti T$ the accumulation points of
 the orbit $G(x, y)$   belong to  the diagonal $\Delta^2\ti T$. By
 Lemma \ref{smconv}
 $\displaystyle\lim_{n\to\infty}g_n(\bb)=p$ so $\displaystyle\lim_{n\to\infty}
 g_n(\ba)=p$. Hence for the above  neighborhood $U_p$ we also have that $U'_p$
 is $(g_n\ba)$-small for some $n\in\N.$
  This
   is impossible by our Convention (1) of  3.1.
  Part 1) is proved.\bx

\bigskip

\noindent 2) Suppose  that   $\{(\ba_i, \bb_i)\in  A\times A\
\vert\ i\in I\}$ is an infinite set such that for every $i\in I$
there is no  $\bc_i\in A$ such that $\ba_i - \bc_i- \bb_i (k)$.
The set $ A$
 is $G$-finite so we can fix $\ba=\ba_i$ and assume
that $\bb_i=g_i(\bb)\ :\ g_i\in G$. Since the space $\ti T$ is
compact, the set $\{\bb_i\}_{i\in I}$ admits an accumulation point
$p$ which  is a limit point for the geometrically finite action
$G\act \ti T.$ By Lemma \ref{geomfin} $p$ is either $k$-conical or
$k$-parabolic point for some (any) $k>1$. Consider these two cases
separately.

Let first, $p$ be a  $k$-conical point. Then there exists ${\bf
c}\in A$ such that ${\bf a} - {\bf c} -p (k).$  By Lemma
\ref{approximation points} we have $\ba - \bc -\bb_i(k)\ (i\in I)$
contradicting the $k$-non-refinability of the pair $\{\ba,
\bb_i\}.$

Let us now suppose that $p$ is  $k$-parabolic. We will now show
that for almost all $i\in I$ the entourages $\ba$ and $\bb_i$
belong to the same horosphere $T_{A,k}(p).$ We claim first that
${\bf a}\in T_{A,k}(p).$ Indeed if not, then there exists $\bc\in
A$
 such that ${\bf a} - {\bf c} - p(k)$ contradicting  by the same argument the
 $k$-non-refinability of the pair $\{\ba, \bb_i\}\ (i\in I).$ So
  $\ba\in T_{A,k}(p).$

Suppose by contradiction that there exist   ${\bb_i}\not\in
T_{A,k}(p)$ for infinitely many  $i\in I$. Then there exist
 $\bc_i\in T_{A,k}(p)$ such that $${\bf b}_i-{\bf
  c}_i
 -p(k).\hfill\eqno(*)$$

We first note that in (*)  we cannot have the same entourage
 $\bc_{0}$ for infinitely many different $\bb_i$. Indeed if not,
 then from (*) we have $\Delta_{\bc_0}(\sh_{\bc_0}\bb_i,c_{0}) > k$
($i\in I)$ for a
 $\bc_{0}$-small set $c_{0}$ containing $p$.
 Since
  $p$ is an accumulation point
 for the set $\{\bb_i\}_{i\in I}$ then $c'_0$ is $\bb_i$-small for infinitely
 many $i\in I$.  Thus
  $c_0\supset \sh_{\bc_0}\bb_i,$ and
  $\Delta_{\bc_0}(c_0, \sh_{\bc_0}\bb_i) \leq 1$ which is
 impossible.

So  we can assume that $\bc_i$ are all distinct.
 By  Lemma
\ref{geomfin} the quotient $T_{A,k}(p)/{\rm Stab_Gp}$ is finite,
so there exists $h_i\in {\rm Stab_Gp}$ such that
$h_i(\bc_i)=\bc\in T_{A,k}(p).$ Hence $h_i(\bb_i)-\bc - p(k)$ for
every $i\in I_1$ where $I_1$ is an infinite subset of $I.$ Since
$\ba\in T_{A,k}(p)$ by Lemma \ref{unicityacp} $p$ is an
accumulation point for the set $\{h_i(\ba)\}_{i\in I_1}.$  Then by
Lemma \ref{approximation points} we obtain $h_i(\bb_i)-\bc -
h_i(\ba) (k)$ and so $\bb_i-h^{-1}_{i}\bc_i-\ba (k)$ which is
impossible.

So
  $\bb_i\in T_{A,k}(p)$ for almost all $i\in I.$ This shows that
the set $A_1$ is $G$-finite.   Part 2)  is proved.  \bx

\bigskip

\noindent 3) a) We omit the index $k$ below. Suppose that the
first set is infinite. Then there exists an infinite set of
$G$-non-equivalent pairs of parabolic points $(p_i, q_i)\in\P^2$
for which $N_d(T(p_i))\cap N_d(T(q_i))\not=\emptyset\ (i\in I)$.
Since the action of $G$ on $\Th^2T$ is cocompact there exist
$g_i\in G$ such that the pair $ (g_i(p_i), g_i(q_i))$ belong to a
compact subset of $\Th^2T$. So without lost of generality we may
assume that the sets $\{p_i\}_{i\in I}$ and $\{q_i\}_{i\in I}$
admits two distinct accumulation points $p$ and $q.$  It follows
from  [Ge1, Lemma P3] that there cannot  exist an entourage
belonging to the intersection of infinitely many distinct
horospheres (for a more general system of horospheres this is also
true, see [GePo3, Corollary of 4.4.2]). So there is an infinite
sequence of distinct entourages $\bb_i\in N_d(T(p_i))\cap
N_d(T(q_i))\ (i\in I)$. The set $\{\bb_i\}_{i\in I}$ admits an
accumulation point $ x\in T.$ Let $(\bc_i)_i\subset T(p_i)$ and
$(\bd_i)_i\subset T(q_i)$ be two subsets for which $d_A(\bb_i,
\bc_i)$ and $d_A(\bb_i, \bd_i)$ are bounded by the constant $d.$
Thus $\da(\bc_i, \bd_i) \leq 2d$ and by Lemma \ref{bbsequences} we
have $p=q=x$. A contradiction.

b) If now the second set is not $G$-finite then for a fixed
parabolic point $p\in \P$ by the part a) we obtain $q\in\P$ such
that the set $ N_d(T(p))\cap N_d(T(q))$ is infinite. Then by
\ref{bbsequences} we must have $p=q.$ The proposition is
proved.\bx

\bigskip

\begin{cor}
\label{constbound}

Suppose that $G$ is a finitely generated group acting
3-discontinuously and 2-cocompactly on a compactum $T.$ Then there
exists a constant $C>0$ such that the $d_A$-diameter of each of
the sets 1), 2) and 3b) of Proposition \ref{boundness} is bounded
by $C$.
\end{cor}

\proof Since $G$ is finitely generated by Lemma \ref{connect} the
graph $\Ga$ is connected. So $d$ is a real distance. The Corollary
follows from the above proposition.\bx

\bigskip

\noindent From   Proposition \ref{boundness},2)  we immediately
have:

\begin{cor}
\label{pairs} Let $\GAT$ be a 3-discontinuous and 2-cocompact
action satisfying the above conditions. Then if for a fixed
$\ba\in A$ and infinitely many $\bb_n\in A$   the pairs $(\ba,
\bb_n)$ are all non-refinable then for all but finitely many $n$
 one has  $(a,\bb_n)\subset T(p)$.\bx

\end{cor}

We  will now obtain few more finiteness properties characterizing
the horospherical projection $\Pi_p : A\to T_{A,k}(p)\ (p\in \P)$.
The following definition is motivated by Lemma \ref{approximation
points}.

\begin{dfn}
\label{visnbd} For a fixed $k > 3$ a {\it visibility neighborhood}
of the point $\bp\in \Pi_p(\ba)\subset T_{A,k}(p)$  from the point
$\ba\in A$ is the following set

$${\mathcal N}(\ba, \bp, p)\ =\ \{\zx\in T_{A,k}(p)\ \vert\ \ba-\bp-p(k)\
\wedge\ \lnot\
 \ba - \bp - \zx (k-1)\},$$
\end{dfn}

\noindent where $\lnot$ denotes the opposite logical statement.

The following proposition establishes the $G$-finiteness
properties of two more sets (by continuing the notations of
\ref{boundness}):

\begin{prop}
\label{finproj} For every $k>1$ the following sets are $G$-finite:

 \begin{itemize}

 \item[\sf 1)] $ {\mathcal A}_2=\{(\zx, \bp)\in T_k^2(p) \vert\ \zx\in {\mathcal N}(\ba,
\bp, p),\ \ba\in A,\ p\in{\mathcal P}\}.$

\medskip

 \item[\sf 2)] ${\mathcal A}_3=\{\Pi_p(T_k(q))\ \vert \ \{p, q\}\subset \P\}.$
\end{itemize}
\end{prop}

\bigskip

\proof 1) Suppose by contradiction that it is not true and  $\A_2$
is not $G$-finite for some $k>1.$ Since $A$ is one $G$-orbit up to
taking an infinite subset of $\A_2$ we  can fix the entourage
$\bp.$ By
  [Ge1, Lemma P3] $\bp$
 can belong to at most
  finitely many different horospheres. So up to a passing to a new infinite subset
   we can fix the parabolic point $p\in \P$.

  If first the set of entourages $\{\ba\ \vert\ (\zx, \Pi_p(\ba))\in\A_2\}$
  is finite, up to choosing a new infinite subset of $\A_2$
  we have $\ba - \bp -p (k)$ and $\lnot\ \ba-\bp-\zx (k-1)$ for a fixed $\ba$.
  Then the set of the first coordinates $\{\zx\ \vert\ (\zx,   \cdot)\in
  \A_2\}\subset T(p)$
  is infinite and by Lemma \ref{unicityacp} its accumulation point is $p.$
  Then by Lemma \ref{approximation
points} there exists $\zx$ in this set such that $\ba-\bp-\zx
(k)$. A contradiction.

If now the set $\{\ba\ \vert\ (\zx, \Pi_p(\ba))\in\A_2\}$ is
infinite let $q\in T$ be its accumulation point. Taking a
$\bp$-small neighborhood $U_q$ of $q$ we obtain that $U'_q$ is
$\ba$-small for every $\ba\in U_q.$ Thus $U_q\supset \sh_\bp\ba.$
Since $\ba-\bp-p (k),$ so $\Delta_\bp(U_q, U_p)
> k-1$ for a $\bp$-small neighborhood $U_p$ of $p.$ It yields
$q-\bp-p(k-1).$ There are infinitely many $\zx\in T(p)$
corresponding to the points $\ba\in U_q$. Since $p$ is the unique
accumulation point of $T(p)$  we must have $\zx\in U_p$ for most
such $\zx.$ Hence
$\Delta(\sh_\bp\ba,\sh_\bp\zx)\geq\Delta_\bp(U_p, U_q)
>k-1$. Therefore $\ba -\bp-\zx (k-1).$ Again a contradiction. \bx

\bigskip

2) Suppose not. Since the set of parabolic points $\P$ is
$G$-finite we can fix the point $p\in \P.$ Using the action of
${\rm Stab}_Gp$ on $T_k(p)$ we can  also assume that there is a
fixed  entourage $\bc\in T(p)$ such that for every $q\in\P\ :\
\bc\in \Pi_p(T(q)).$ So there exists an infinite set $\{\bd_i\in
\Pi_p(T(q_i))\ \vert\ i\in I, q_i\in \P\}$  such that for all
$i\in I$ we have

$$\bb_i-\bd_i-p (k),\ \ba_i - \bc - p (k), \{\ba_i, \bb_i\}\subset
T(q_i).$$

\noindent  Since $p$ is the unique accumulation point of $T(p)$,
up to passing to an infinite subsequence of $I$, we may assume
that $\displaystyle\lim_{i\to\infty}\bd_i= p$. Then by Lemma
\ref{smconv} we have $\displaystyle\lim_{i\to\infty}\bb_i= p$. Let
$q\in T$ be an accumulation point of the set $\{q_i\}_{i\in I}$.
We claim that $q=p$. Indeed if not then there exists an entourage
$\ba\in A$ such that $q-\ba-p (k)$. Hence for infinitely many
$i\in I$ we have $q-\ba-\bb_i (k)$ (Lemma \ref{approximation
points}). Then $\Delta_\ba(U_q, \sh_\ba\bb_i)
> k$ for every $\ba$-small   neighborhood $U_q$ of $q$.
  So for some $i\in I$ we have  $q_i\in U_q$ and hence
$q_i-\ba-\bb_i(k)$. The latter one  is impossible since $\bb_i\in
T_k(q_i).$ If the set $\{\ba_i\}$ has an accumulation point $r$
different from $p$ then   $\exists\ \ba\in A\ : \ r-\ba-p(k)$. So
as above we have $q_i-\ba-\ba_i(k)$ which is impossible by the
same reason. So $\displaystyle\lim_{i\to\infty}\ba_i=p.$ Then for
every $\bc$-small neighborhood $U_p$ of $p$ we have that $U_p$ is
also $\ba_i$-small and $U'_p$ is  $\ba_i$-small for infinitely
many $i.$ This is impossible. The proposition is proved. \bx

\bigskip

The following Corollary gives a uniform bound on the cardinality
of the intersection of the stabilizers of parabolic points for a
geometrically finite action.

\begin{cor} \label{parinter}
 Let
$G$ be a group admitting a 3-discontinuous and 2-cocompact action
on a compactum $T$. Then there is a constant $C$ such that for
every pair of distinct parabolic points $p_i$ and $p_j$  for the
action $G\act T$ one has
$$\vert {\rm Stab}_Gp_i\cap {\rm Stab}_Gp_j\vert\leq C,$$

\end{cor}

\proof Denote $H_i={\rm Stab}_Gp_i.$ Suppose the statement is not
true. By Lemma \ref{geomfin} the set of   parabolic points for the
action $G\act T$ is $G$-finite. So up to conjugation we can
suppose that there exists a sequence of the stabilizers of
parabolic points $H_0, H_n\ (n\in\N)$ such that $\vert H_0\cap
H_n\vert\to\infty.$ Let $T_n$ be a   horosphere at $p_n \ (n\in
\N\cup\{0\})$. Then the projection $\Pi_{p_0}(T_n)$ of $T_n$ on
$T_0$ is invariant under $H_0\cap H_n.$ Since the action $H_0\act
T_0$ is discontinuous we have $\vert \Pi_{p_0}(T_n)\vert\to
\infty$ which is impossible by Proposition \ref{finproj}.2. \bx

\medskip

\noindent {\bf Remark}. The above Corollary is also true if $G$ is
a countable group acting 3-discontinuously on a compactum $T$ such
that every point $T$ is either conical or bounded parabolic.
Indeed in this case by [Ge1, Main Theorem.c] the space $T$ is
metrisable. So by \cite{Tu3} the action $G\act T$ is 2-cocompact
and the above Corollary holds.

\subsection{Proof of Theorem A.}

The aim of this subsection is the following.

\medskip

 \noindent {\bf Theorem A.} {\it  Let $G$ be a relatively hyperbolic group  with respect to a
 collection of
parabolic subgroups $\{P_1,..., P_n\}.$ Then   $G$  is the
fundamental group of the following finite ``star graph''

\vskip-5pt
$$\vcenter{\pdFile{grph0}{200}{180}}
\eqno(1)$$ \vskip5pt
 \noindent  whose central
  vertex group  $G_0$ is   finitely generated relatively hyperbolic  with respect to
  those edge groups $Q_i=P_i\cap G_0$ which are infinite,  all other vertex groups
of the graph   are $P_i\ (i=1,...,n)$.

Moreover for every finite set $K\subset G$ the subgroup $G_0$ can
be chosen to contain $K$.

}

\medskip

\proof  Recall that $\displaystyle A=G(\ba_0)\ (\ba_0\in \ent)$ is
a discrete orbit of entourages  forming the vertex set of the
graph $\Ga$ satisfying our Convention \ref{setA}. Without lost of
generality we can assume that the group $G$ is not finitely
generated and $\ba_0\in K.$ So the graph $\Ga$ is not connected
(see Lemma \ref{connect}). The distance $\da(\zx, \by)$ is a
pseudo-distance being infinity if and only if $\zx$ and $\by$
belong to different connected components of $\Ga.$ By Lemmas
\ref{geomfin} and \ref{unicityacp} the set $\P$ of parabolic
points for the action $G\act T$ is $G$-finite; and for every
$p\in\P$ the stabilizer $H_p=\rm Stab\sl_Gp$ acts cofinitely on
its horosphere $T(p).$

Let  $\A_i\ (i=1,2,3)\subset A^2$ be the $G$-finite sets
introduced in Propositions \ref{boundness}.2 and \ref{finproj}.

We now construct a new graph $\ti\Ga$ whose   set of vertices is
$A$ and the set of edges is given by   the   pairs of entourages
belonging to the following sets:

\bigskip

a) the finite set   $K^2$ and the set of all its horospherical
projections $\{\Pi_p(K^2)\ \vert\ p\in\P\};$

\medskip

b)  the   set $\A_1$ and the set of all its horospherical
projections $\{\Pi_p(\A_1)\ \vert\ p\in\P\};$

\medskip

c)    the set $\A_2;$

\medskip

 d)   the  set
  $\A_3$.

\bigskip

All these sets are $G$-finite. Indeed the set $\A_1$ is $G$-finite
by Proposition \ref{boundness}.2. So  by Proposition
\ref{boundness}.1 the set $\{\Pi_p(\A_1)\ \vert\ p\in\P\}$
consisting of the projections
 of finitely many $G$-orbits of
 pairs is $G$-finite too.
 The sets $\A_2$  and $\A_3$ are $G$-finite by Proposition
 \ref{finproj}.

\begin{lem}
\label{exhaust} There exists a finitely generated  subgroup $G_0$
of $G$ containing any finite  subset $K\subset G$ and which is
relatively hyperbolic with respect to $Q_i=P_i\cap G_0\
(i=1,...,n).$
\end{lem}

 \proof Let $\Ga_0$ be
the connected component of $\ti \Ga$ containing $K.$ Set $ G_0=
{\rm Stab}_G\Ga_0$  and $ A_0=\Ga_0^0.$ By Lemma \ref{connect} the
group $G_0$ is finitely generated. We are left to prove that $G_0$
is relatively hyperbolic with respect to the subgroups
$\{Q_i\}_{i=1}^k.$

Let $T_0$ be a subset of $T$ which is the limit set of $G_0.$ We
will first show that the action $G_0\act T_0$ is 2-cocompact. By
[Ge1, Prop. E] the 2-cocompactness is equivalent to the
$k$-separation property:

$$\forall p, q\in T_0\ :\ p\not=q\ \exists\ \bb\in
A_0\ : p-\bb - q (k),\hfill\eqno(1)$$ \noindent for some $k>0.$
Since the action of $G$ on $T$ is 2-cocompact, the property (1) is
true for some $\bb\in A.$
 If $\bb\in A_0$ we are done, so suppose that $\bb\not\in
A_0.$ Let $U_p$ and $U_q$ be $\bb$-small neighborhoods of the
points $p$ and $q$  such that $\Delta_{\bb}(U_p, U_k) >k.$ Since
$p$ and $q$ are accumulation points of $A_0$ there exist
entourages $\ba, \bc\in A_0$ such that $U'_p$ is $\ba$-small and
$U'_q$ is $\bc$-small. So $U_p\supset \shba$ and $U_q\supset
\shbc.$ Hence
$$\ba-\bb-\bc (k).\hfill\eqno(2)$$

By Proposition \ref{exent} up to refining the pair $\{\ba, \bb\}$
  we can suppose that the pair $\{\ba, \bb\}$ is
$k+2$-nonrefinable. Since $\bb\not\in A_0,$ by operation b) above
the pair $\{\ba,\bb\}$ must belong to an horosphere $T_{k+2}(r)\
(r\in \P)$. As $\{\ba, \bc\}\subset A_0$ and $\Ga_0$ is connected
there exists a path $\ga=\ga(\ba, \bc)\subset \Ga_0.$ Let
$\be=\Pi_r(\bc).$ Note that  for every edge $l\in\Ga_0^1$ we have
$\Pi_r(l)\in\Ga_0^1$. Indeed  if $l$ joins two vertices of $A_0$
then  by the operations a), b) and d) all their horospherical
projections are joined by edges too. So $\Pi_r(\Ga_0)\subset
\Ga_0.$ Since $\{\ba, \be\}\subset T(r)\cap \Pi_r(\ga)$ we have
$\be\in A_0.$

Operation c) then implies that $\bb\not\in {\mathcal N}(\bc, \be, r).$
By Definition \ref{visnbd} we have
$$\bb-\be-\bc (k+1).\hfill\eqno(3)$$   So
$\shbc\subset\sh_{\bb}\be$ and (2) yields $\Delta_{\bb} (\shba,
\sh_{\bb}\be)> k-1$ and  $ \ba-\bb-\be (k-1).$ Thus
$\sh_{\be}\ba\subset\ \sh_{\be}\bb$ and by (3) we have $\ba-\be
-\bc (k-1)$ with $\be\in A_0.$ We have proved that the action
$G_0\act T_0$ is $(k-1)$-separating and so is 2-cocompact [Ge1,
Prop. E].

By [Ge1, Main Theorem] every  point of $T_0$ is either conical or
parabolic for the action of $G_0$ on $T_0.$ Let $p\in T_0$ be a
parabolic point for this action. We need the following.

\bigskip

\noindent {\bf Claim.} {\it The point $p$ is also parabolic for
the action of $G$ on $T.$}

\medskip

 {\it Proof of the claim.}  Suppose not. Let $T(p)\subset
A_0$ be a
 horosphere
 for the action $G_0\act \Ga_0$. Let us choose $\bb\in T_{k-3}(p)\subset
  T(p)\ (k>3)$ where $T_{s}(\cdot)$ denotes the
  "sub-horosphere" of $T(\cdot)$
of order $s$ (see \ref{hor}).

  Suppose
first that $\bb$ does not belong to any horosphere $\ti
T_{k}(q)\subset A$ for the action $G\act \ti\Ga$. Since $p$ is
conical for the action $G\act T$ there exists $\bc\in A$ such that
$\bb-\bc-p (k-1).$ Note that $\bc\not\in A_0$ as otherwise
$\bb\not\in T_{k-1}(p)$ which is impossible as $T_{k-3}(p)\subset
T_{k-1}(p).$ By the sublemma below we can also suppose up to
refining the couple $(\bb, \bc)$ that it is not $k$-refinable ($k>
3$). Since $\bb$ and $\bc$ do not belong to one horosphere in
$\ti\Ga$, by operation b) above, $\bc$ and $\bb$ are joined by an
edge in $\ti\Ga.$ So $\bc\in A_0$ and we have a contradiction in
this case.

We affirm now that   there exists  $h\in {\rm Stab}_{G_0}p$ such
that $h(\bb)$ does not belong to any horosphere $\ti T_{k}(q)$
where $q\in\P$. Suppose not, then $\bb\in T_{k-3}(p)\cap \ti
T_k(q)$ for some $q\in\P.$ Again since $p$ is conical for the
action on $T$ there exists $\bc\in A\setminus A_0$ such that
$\bb-\bc-p (k-1).$ By the argument above we can assume  that
$\bc\in \ti T_k(q)$ too. Up to choosing $h\in {\rm Stab}_{G_0}p$
so that $\bb_1=h(\bb)\in T_{k-3}(p)$ is sufficiently close to $p$
we can also assume that $\bb-\bc-\bb_1 (k-1)$ (Lemma
\ref{approximation points}). As the distance $d_{A_0}(\bb, \bb_1)$
is large, by Proposition \ref{boundness}.3b we have that
$\bb_1\not\in \ti T_k(q)$. Then there exists $q_1\in
\P\setminus\{p,q\}$ such that $\bb_1\in \ti T_k(q_1)$. By the
argument above giving the formula (3)  it follows that there
exists $\be\in \Pi_q(\bb_1)\cap A_0$ such that
$\bb_1-\be-\bc(k-1)$ and so $\bb-\be-\bb_1(k-2).$ Continuing in
this way we obtain an infinite sequence $\bb_n=h_n(\bb)\in
T_{k-3}(p)\cap \ti T_k(q_n)$ where $h_n\in {\rm Stab}_{G_0}p$ and
$q_n=h_n(q)$ are all different parabolic points. By Proposition
\ref{boundness}.1  it follows that the subset $\displaystyle
B=\bigcup_{n\in \N}\Pi_q(h_n(\bb))$ of $\ti T_k(q)\cap A_0$ is
finite. So up to choosing a new subsequence for a fixed $\be\in B$
we have $\bb-\be-\bb_n (k-2)\ (n\in\N).$ Since $p$ is the
accumulation point of $\{b_n\}_{n\in\N}$, for any $\be$-small
neighborhood $U_p$ of $p$ its complement $U'_p$ is $\bb_n$-small
for infinitely many $n.$ Thus $U_p\supset \sh_{\be}\bb_n$ and so
$\Delta_{\be}(\sh_{\be}\bb, U_p) > k-3$ implying $\bb-\be-p(k-3).$
Therefore $\bb\not\in T_{k-3}(p)$ which is a contradiction proving
the claim.\bx

\medskip

We have ${\rm Stab}_{G_0}p={\rm Stab}_Gp\cap G_0$. Lemma
\ref{exhaust} is proved modulo the following lemma.
\begin{sublem}
\label{refhor} If $\ \bb-\bc-p(k-1)$ and $\bb-\bc_1-\bc (k)$ then
$\bb-\bc_1-p(k-1)\ (k > 3).$

\end{sublem}

\proof Let us first show that $\bc_1 - \bc - p(k-2).$ Indeed the
second assumption implies that $\sh_{\bc}\bc_1\supset
\sh_{\bc}\bb.$ So for a $\bc$-small neighborhood  $U_p$ of $p$
using the first assumption for any $c\in \Sh_{\bc}\bc_1$ we have

$$\Delta_{\bc} (c, U_p) > \Delta_\bc (\shcb, U_p) -
\ti\Delta_\bc (\shcb, c) > k-2.$$

 So $U_p\subset c'\in\Sh_{\bc_1} \bc$ and $\ti\Delta_{\bc_1}(\sh_{\bc_1}\bc, U_p) \leq 1.$
 Hence
$\Delta_{\bc_1}(\sh_{\bc_1}\bb, U_p) >
\Delta_{\bc_1}(\sh_{\bc_1}\bb, \sh_{\bc_1}\bc) -
\Delta_{\bc_1}(\sh_{\bc_1}\bc, U_p) > k-1.$ The lemma and  the
proposition are proved.\bx

\bigskip

The  following lemma finishes the proof of the Theorem.

\begin{lem}
\label{graphreq} The action $G\act\ti{\mathcal G}$ induces an action
on a bipartite simplicial tree $\T$ such that the graph $X={\T}/G$
satisfies Theorem \ref{graph}.
\end{lem}

 \proof Using the graph $\ti\Ga$ we construct the tree $\mathcal T$ to have vertices belonging to two subsets ${\mathcal C}$ and $\mathcal H.$
 The elements of ${\mathcal C}$ are  components of $\ti\Ga$ and the elements of
 ${\mathcal H}$ are the horospheres of $A=\ti\Ga^0$.
 We call them {\it
non-horospherical}
 and {\it horospherical} respectively. Two vertices $C$ and $H$ of
 $\mathcal T$
 are joined by an edge if and only if $C\in {\mathcal C}$, $H\in{\mathcal
 H},$ and $C\cap H\not=\emptyset.$

Let us first show that $\T$ is connected. Indeed by construction
every horospherical vertex is joined with a non-horospherical one.
So it is enough to prove that every two non-horospherical vertices
can be joined by a path. Let $C_i\ (i=1,2)$ be the corresponding
connected  components of $\ti\Ga$ and let us fix two entourages
$\ba\in C^0_1$ and $\bb\in C^0_2$. By Proposition \ref{exent}
there exists a non-refinable tube between them: $P=\ba -
\bb_1-...-\bb_n-\bb\subset A.$ By operation b)  above every
non-refinable pair $(\bb_i, \bb_{i+1})$ either belongs to an
horosphere $T(p)$ or corresponds to an edge in the graph $\ti\Ga.$
In the latter case it stays in the same component of $\ti\Ga.$ In
the former case the horosphere $T(p)$ corresponds to a single
vertex of the graph $\T$. So the tube $P$ produces a path in  $\T$
between the corresponding vertices. Thus $\T$ is connected.

Let us now show that $\mathcal T$ is a tree. Suppose not and  it
contains a simple loop $\al$. Since the vertices of two types
alternate on $\al$ we can fix a horospherical vertex $H$
corresponding to the horosphere $T(p)$ and having two
non-horospherical neighboring vertices $C_1$ and $C_2.$ Let
$\al_1$ be a subpath of $\al$ containing the vertices $H$, $C_1$,
$C_2,$  and $\al_2$ be the closure of $\al\setminus \al_1$. The
path $\al_2$ corresponds to an alternating sequence of components
of $\ti\Ga$ and  horospheres. So we can choose a sequence of tubes
$P_i\subset C_i$ where each  $C_i\ (i\geq 3)$ is a component of
$\ti\Ga$ corresponding to a non-horospherical vertex of $\al_2$.
The tube $P_i$  connects two entourages from $C_i$ each belonging
to horospheres $T(q_i)$ and $T(q'_i)$ intersecting $C_i.$ Note
that these horospheres differ from the initial horosphere $T(p)$
as $\al$ is a simple loop. By operations b) and d) above it
follows that that there exists the path $\displaystyle\bigcup_i
\Pi_p(P_i\cup T(q_i)\cup T(q'_i))$   on $T(p)\cap \ti\Ga.$ It
implies that the vertices $C_1$ and $C_2$ correspond to the same
connected component of $\ti\Ga$ which is impossible. So $\T$ is a
tree.

  By
Lemma \ref{exhaust} we can assume that the stabilizer $G_0$ of a
component $\Ga_0\in {\mathcal C}$   is finitely generated and contains
the fixed finite set $K\subset G.$  The group $G$ acts
transitively on $A$ and so on $\mathcal C.$   Then every  element of
$\mathcal C$ is stabilized by a subgroup conjugate to $G_0$. So in the
graph $X={\mathcal T}/G$ there is only one non-horospherical vertex
$v_0={\mathcal C}/G$ whose vertex group is $G_0.$

The set of horospheres on $T$ is $G$-finite (Lemma \ref{geomfin})
so $X$ contains $n$ vertices of non horospherical type each
representing the $G$-orbit of an horosphere $T(p)\ (p\in \P).$
Every one of them is connected with $v_0$ by a unique edge. So
every vertex group of horospherical type is $P_i$ and the edge
groups are $Q_i=P_i\cap G_0\ (i=1,...,n)$.
 The Theorem
is proved.\bx

\subsection{Corollaries of Theorem A}

 Theorem A admits several immediate corollaries describing
different type of finiteness properties of relatively hyperbolic
groups.

\medskip

\begin{cor}
\label{fingenexhaust}
 Let $G$ be a relatively hyperbolic group with respect to the system $P_j\ (j=1,...,n)$.
 Then there exists an exhaustion
$\displaystyle G=\bigcup_{i\in I} G_i$ where $G_i$ is a finitely
generated group which is relatively hyperbolic   with respect to
the system $P_j\cap G_i\ (j=1,...,n).$\bx
\end{cor}

\medskip

\begin{dfn}
\label{relfingen} A group $G$ is called relatively finitely
generated with respect to a system $\frak P$ of subgroups if it is
generated by the system $\frak P$ and a finite set $S$ generators.

Furthermore $G$ is relatively finitely presented with respect to
$\frak P$ if there are at most  finitely many relations between
the elements of $S.$
\end{dfn}

\medskip

\begin{cor}
\label{fingen}
 Let $G$ be a group acting 3-discontinuously and
2-cocompactly on a compactum T. Then $G$ is relatively finitely
generated with respect to the stabilizers of the parabolic points.
\end{cor}

\proof Indeed by Theorem A the group $G$ is generated by a
finitely generated subgroup $G_0$ and by the parabolic subgroups
$H_i\ (i=1,...,n)$. The Corollary follows.\bx

\medskip

\begin{cor}
\label{noparab}

A group $G$ acting 3-discontinuously and 2-cocompactly on a
compactum $T$ without parabolic points   is finitely generated.\bx
\end{cor}

\noindent {\bf Remark.} If in particular   $G$ acts
3-discontinuously and 3-cocompactly
 on $T$ without isolated points then  every
 point of $T$ is conical [GePo1, Appendix]. So by Corollary \ref{noparab} $G$ is
 finitely generated in this case. By a direct argument
 one can now deduce that $G$ is word-hyperbolic  [GePo1, Appendix].
This  provides a new proof of a theorem due to B.~Bowditch
 \cite{Bo3}.\bx

\medskip

 Before we state the next corollary let us recall two more definitions of relative
 hyperbolicity valid for   infinitely   generated groups. The
 first one is due to B.~Bowditch:

\begin{defi}\cite{Bo1} \label{bodef} A graph $\Gamma$ is called fine
if there are at most finitely many simple arcs of a bounded length
with fixed endpoints.

An action of a group $G$ on a graph $\Gamma$ is proper on the set
 of edges $\Gamma^1$ if the stabilizers of edges are   finite, the
action is called cofinite if $\vert\Gamma^1/G\vert < \infty$.

A  group $G$ is called relatively hyperbolic with respect to a
system  of subgroups $\frak P$ if $G$ acts non-parabolically   on
a connected fine hyperbolic graph $\Gamma$ cofinitely and properly
on edges such that   $\frak P$ is a maximal set of non-conjugated
infinite
  stabilizers of vertices.
\end{defi}

\medskip

\noindent The second definition is due to D.~Osin:

\begin{defi}\cite[Definions 2.3, 2.30, 2.35]{Os}\label{osdef}
A group $G$ is relatively hyperbolic with respect to a finite
system $\frak P$ of subgroups of infinite index if it is finitely
presented with respect to $\frak P$ and the corresponding relative
Cayley graph $\Gamma={\rm Cay}(G, \frak P\cup S)$ admits a linear
relative Dehn function i.e. the relative area of a cycle in
$\Gamma$ of length $\leq n$ is bounded by a linear function of
$n$.
\end{defi}

\medskip

 The main corollary of Theorem A is the following result
establishing that   all known definitions   of the relative
hyperbolicity valid for a group  of any cardinality   are
equivalent to the existence of the star-graph of groups
decomposition (1) of Theorem A.

\begin{theor}\label{equiv} The following conditions are equivalent
for a group $G.$

 \begin{itemize}

   \item[\sf 1)] Definition \ref{ourdfn}.

 \item[\sf 2)] Bowditch's definition \ref{bodef}.

 \item[\sf 3)] Osin's definition \ref{osdef}.

\item[\sf 4)] $G$ admits the star-graph of groups decomposition
(1) of Theorem A  where the central vertex group $G_0$ is a
finitely generated relatively hyperbolic group with respect to
those edge groups $Q_i$ which are infinite.\bx
\end{itemize}
\end{theor}

\medskip

\noindent {\bf Remark.} In fact in \cite{Os} the finiteness of the
system $\frak P$ is not required. We need it to have the
equivalence of the Osin's  definition to all others definitions
which all imply this assumption.

 \proof As it was mentioned in the
Introduction modulo Theorem A and known facts the proof of the
theorem goes according to  the diagram (*) from the Introduction.
It remains only to show Propositions \ref{inverse} and \ref{bowos}
whose proofs are given below.

 \medskip

\begin{prop}\label{inverse}   Let a group $G$ admit  the graph $\mathfrak G$ of groups
decomposition (1).  Suppose that the subgroup $G_0$ acts on a fine
$\delta$-hyperbolic graph $\Gamma$ properly and cofinitely on the
set of edges    and the groups $Q_i$ are the stabilizers of
$G_0$-non-equivalent vertices  of $\Gamma$.   Then there exists an
action of $G$ on a fine $\delta$-hyperbolic graph $\Delta$
properly and cofinitely on $\Delta^1$ such that $P_1,P_2,\dots$
are the stabilizers of $G$-non-equivalent vertices.

In particular, if $G_0$ is relatively hyperbolic in the sense of Bowditch
with respect to the subgroups $Q_i$ then $G$ is
relatively hyperbolic in the same sense with respect to the
subgroups $P_i\ (i=1,...,n)$.
\end{prop}

\medskip

Notice that the statement above is more general than the
implication $4)\Longrightarrow 2)$ as we do not need to assume
that $G_0$ is finitely generated and that all the subgroups $P_i$
are infinite. The group $G$ is relatively hyperbolic with respect
to those which are infinite.

\bigskip

 {\it Proof of the proposition.} We will construct the
graph $\Delta$ as the quotient of another graph
$\widetilde\Delta$.

Let $\mathcal T$ be the universal covering tree of the graph
$\mathfrak G$ and let $\tau:\mathcal T\to \mathfrak G$ be the
covering map. For a vertex $v\in\tau^{-1}\{G_0\}{\subset}\mathcal
T^0$ denote by $G_{0,v}$ its stabilizer $\mathsf{St}_Gv$ in $G$
and by $\widetilde\Gamma_v$ a copy of $\Gamma$ on which the group
$G_{0,v}$ acts. We can assume that all these copies are disjoint.

There is a bijection between the $\mathcal T$-edges incident to
$v$ and the vertices of $\widetilde\Gamma_v$. Using this bijection
we can replace a small neighborhood of $v$ in $\mathcal T$ by
$\widetilde\Gamma_v$ joining the $v$-endpoint of an edge with the
corresponding vertex of $\widetilde\Gamma_v$. This implantation
can be made $G$-equivariantly, so the obtained graph
$\widetilde\Delta$ is connected and acted upon by $G$.

There are two kinds of the edges of $\widetilde\Delta$: those of
the graphs $\widetilde\Gamma_v$ and those of $\mathcal T$. Now we
collapse all the edges of the second kind. Let $\Delta$ denote the
resulting graph. Since the collapsing is $G$-equivariant there is
an induces action $G{\curvearrowright}\Delta$ and the projection
map $\pi:\widetilde\Delta\to\Delta$ is a $G$-equivariant morphism
of graphs. Hence the graph $\Delta$ is connected too.

\bigskip

We have the diagram
\begin{picture}(70,27)(-30,-20) \put(0,0){$\widetilde\Delta$}
\put(-2,-2){\vector(-1,-1){20}} \put(8,-2){\vector(1,-1){20}}
\put(-18,-11){$\pi$} \put(19,-11){$\sigma$}
\put(-30,-30){$\Delta$} \put(30,-30){$\mathcal T$}
\end{picture} of $G$-equivariant graph morphisms where $\sigma$ collapses
$\widetilde\Gamma_v$
to $v$.
\vskip3pt
Denote by $\Gamma_v$ the $\pi$-image of $\widetilde\Gamma_v$ is
isomorphic to $\widetilde\Gamma_v$. Unlike the
$\widetilde\Gamma_v$'s the subgraphs $\Gamma_v$ of $\Delta$ are
not disjoint.

\medskip

\noindent{\it Canonical lifting of the paths.} For a vertex
$w{\in}\Delta^0$ the subgraph $\pi^{-1}w$ of $\widetilde\Delta$ is
either
the star of a vertex $\widetilde w{\in}\mathcal T^0$ such that
$\tau\widetilde w$ is one of the vertices $P_i$ or a vertex in some $\Gamma_v$. We say that
$\widetilde w$ is the \it central representative \rm of $w$ in the first case.

\medskip

\sc Lemma\sl. For every locally injective path $\gamma:I\to\Delta$
between two vertices there exists a unique locally injective path
$\widetilde\gamma:J\to\widetilde\Delta$ between the central
representatives of the endpoints of $\gamma$ and a monotone map
$\iota:J\to I$ such that
$\gamma{\circ}\iota{=}\widetilde\gamma{\circ}\pi$. If $\gamma$ is
geodesic then $\widetilde\gamma$ also is\rm.

\proof Both existence and uniqueness follow from the fact that the
$\pi$-preimages of vertices are connected subtrees and the maximal
subpaths in $\Gamma_v$'s lift uniquely to paths in
$\widetilde\Gamma_v$. These lifted subpaths can be joined uniquely
in the corresponding subtrees providing the lift of the whole path.

 The statement about geodesic paths
follows from the fact that every locally injective path in a tree
is geodesic. The lemma is proved.

\medskip

\noindent{\it Verification of the properties of the action
$G{\curvearrowright}\Delta$.} By construction the action
$G{\curvearrowright}\Delta^1$ is proper and cofinite.

Recall that a \it circuit \rm in a graph $\Gamma$ is a subgraph
homeomorphic to the circle  $\Bbb S^1$. Every circuit of $\Delta$
is contained in a subgraph $\Gamma_v$ since otherwise its lift is
a circuit in $\widetilde\Delta$ containing $\mathcal T$-edges so
its $\sigma$-image is a non-trivial circuit in the tree $\mathcal
T$ which is impossible. This implies that $\Delta$ is fine.

It remains to verify that $\Delta$ is $\delta$-hyperbolic provided
that $\Gamma$ is. Let $\tau$ be a geodesic triangle in $\Delta$
and let $\widetilde\tau$ denote the triangle in $\widetilde\Delta$
obtained from $\tau$ by the canonical lifting of the sides. By
Lemma $\widetilde\tau$ is geodesic in $\widetilde\Delta$ and
$\sigma(\tau)$ is geodesic in $\mathcal T$. Thus $\widetilde\tau$
actually consists in pieces that either $\mathcal T$-edges or
geodesic bigons in subgraphs $\Gamma_v$ or (at most one) geodesic
triangle in some $\Gamma_v$. Each this piece is $\delta$-thin by
the hypothesis. The Proposition is proved.\bx

\medskip
 \noindent {\bf Remark.}
Proposition \ref{inverse} gives a generic construction
of non-finitely generated relatively hyperbolic groups.
By Theorem A any relatively hyperbolic
group can be constructed in this way.

\medskip

 Furthermore  the
implication $2)\Longrightarrow 1)$ of \ref{inverse} and
\cite[remark 9.1]{Ge2} imply that $G$ acts $3$-discontinuously and
$2$-compactly on a compactum such that the $\{P_i:i=1,...,n\}$ is
a complete list of representatives of the stabilizers of parabolic
points.

\medskip

 The rest of the subsection is devoted to the proof of the
implication $2)\Longrightarrow 3)$.
  The argument is rather
standard. However we did
not find an adequate reference for non-finitely generated groups,
so for the reader's convenience we provide it here.
The argument below is motivated by
\cite{Bo4}.

\medskip
In this section we
 consider graphs
as 1-dimensional CW-complexes and cycles as non-oriented cycles.

Let $\Xi$ denote the graph with two vertices $P,Q$
and three edges $a,b,c$ that join $P$ with $Q$.

Every continuous map $\varphi:\Xi\to\Gamma$ determines three
 cycles in $\Gamma$; denote them by
$\varphi_{ab}$, $\varphi_{bc}$, $\varphi_{ca}$.

\noindent\kern390pt\begin{picture}(100,94)(0,-10)
\put(0,0){\pdFile{50EC2DBA0}{100}{80}}
\put(46,82){$P$}
\put(3,37){$a$}
\put(43,37){$b$}
\put(92,37){$c$}
\put(46,-10){$Q$}
\end{picture}

\vskip-85pt
\skip255\rightskip
\rightskip120pt

A non-negative function $\boldsymbol\alpha$ on the set
\{cycles in $\Gamma$\} is called\hfil\penalty-10000 a
\it pseudo-area function \rm on $\Gamma$ if
$\boldsymbol\alpha(\varphi_{ac})\leqslant
\boldsymbol\alpha(\varphi_{ab}){+}\boldsymbol\alpha(\varphi_{bc})$
for every continuous map $\varphi:\Xi\to\Gamma$.

Let $\mathcal C$ be a set of circuits in a graph $\Gamma$. Denote
by $\Gamma{+}\mathcal C$ the CW-complex obtained from
$\Gamma$ by attaching a 2-cell to each circuit
in $\mathcal C$.

For a locally injective map $\gamma:\Bbb S^1\to\Gamma$ denote by
$\ell(\gamma)$ the natural \it length \rm of $\gamma$: this is the
number of edges that $\gamma$ consecutively passes.

\rightskip\skip255

\begin{lem}\label{finitecomp}  Let a group $G$ act on a fine
hyperbolic graph $\Gamma$ properly and cofinitely on the set  of
edges $\Gamma^1$ of $\Gamma$. Then there exists a $G$-finite set
$\mathcal C$ of circuits in $\Gamma$ such that the complex
$\Gamma{+}\mathcal C$ is simply-connected and for every
$G$-invariant pseudo-area function $\boldsymbol\alpha$ on $\Gamma$
there exists a number $M$ such that
$\boldsymbol\alpha(\gamma)\leqslant M\ell(\gamma)$ for every
locally injective map $\gamma:\Bbb S^1\to\Gamma$.
\end{lem}

\vskip-8pt
\noindent\kern315pt\hbox{\pdFile{50EC2DB90}{170}{136}
\begin{picture}(0,0)(0,30)
\put(-6,132){$A$}
\put(-181,53){$B$}
\put(-180,102){$C$}
\put(-155,23){$D$}
\put(-96,96){$\gamma_1$}
\put(-169,66){$\gamma_2$}
\end{picture}}
\vskip-140pt
\skip255\rightskip
\rightskip200pt
{\it Proof. }\rm
Consider a locally injective map
$\gamma:\Bbb S^1\to\Gamma$. Choose points $A,B$ in the image of
$\gamma$ such that
$|AB|{=}\mathsf{diam\,Im}(\gamma)$.
We will
show that $B$ belongs to a non-geodesic piece of
$\gamma$ of
a bounded length.

Suppose that the $d$-neighborhood in $\gamma$ of the point $B$ is
a geodesic segment. Then we have $|BC|{=}|BD|{=}d$ and $|CD|{=}2d$.
This implies:\hfil\penalty-10000
 $|AB|{+}|CD|\geqslant
d{+}\mathsf{max}\{|AC|{+}|BD|,|AD|{+}|BC|\}$.\hfil\penalty-10000
 Since our graph is
hyperbolic, the value of $d$ is bounded by some constant $\delta$
\cite[1.1.A]{Gr}.
 \vskip5pt\rightskip\skip255 Hence for
$d=\delta+1$ the arc $CBD$ is not geodesic. We join $C$ with $D$
by a geodesic segment $\lambda$.
 Let \vbox to17pt{\vfil}$\gamma_1=\lambda\cup{\rm arc}(CAD)$
and $\gamma_2=\lambda\cup{\rm arc}(CBD)$. The length of $\gamma_2$
is bounded by $4\delta+3$.

Thus the cycle $\gamma_2$ possesses a tesselation by at most
$4\delta{+}3$ simple cycles of length at most $4\delta{+}3$.

 Since
$\Gamma$ is fine and $\Gamma^1/G$ is finite, the set $\mathcal
C{=}\{$circuits of length ${\leqslant}4\delta{+}3\}$ is
$G$-finite. As $\boldsymbol\alpha$ is $G$-invariant there exists a
constant $M$ such that $\boldsymbol\alpha(\gamma_2)\leqslant M$.
By the definition of pseudo-area
$\boldsymbol\alpha(\gamma)\leqslant \boldsymbol\alpha(\gamma_1)
{+}\boldsymbol\alpha(\gamma_2)$. Since
$\ell(\gamma_1)\leqslant\ell(\gamma){-}1$ the desired inequality
follows by induction on $l(\ga)$. The lemma is proved.

\medskip

 Let $S$ denote a (finite) relative system of generators and
$P_i\ (i\in I)$ be the system of all maximal parabolic subgroups
of $G$. Recall  that  the set of vertices of the coned-off Cayley
graph ${\rm Cof}(G, S\cup\P)$ is $G\cup\P$ where $\P$ is the set
of the parabolic points $p_i$ whose  stabilizer is $P_i$. The set
of edges is $E_1\cup E_2$ where $E_1$ is the set of edges of the
absolute Cayley graph ${\rm Cay}(G, S)$ and $E_2$ is the set of
edges joining every parabolic element in $P_i$ with the point
$p_i$ fixed by $(i\in I)$ \cite{Fa}.

\begin{lem}
\label{chgraph} Let a group $G$ admit a cofinite and proper on
edges action on a fine hyperbolic graph $\Gamma$.
Then there exists a finite system $S$ of relative generators of
$G$ such that the coned-off Cayley graph ${\rm Cof}(G, S\cup\P)$
is also fine and hyperbolic and $G$ acts on it   cofinitely and
properly on edges, where   $\P$ is a set of parabolic vertices of
$\Gamma.$
\end{lem}

\proof We present a finite algorithm for
passing from $\Gamma$ to ${\rm Cof}(G, S\cup\P)$ by keeping all the above
properties valid.

Consider first an intermediate graph $\Delta$ obtained from
$\Gamma$ as follows. Set $\Delta^0=\Gamma^0\sqcup G.$ We call the
vertices \rm and the edges of $\Gamma$ \it blue \rm
and the elements of $G$ \it red\rm.
 Denote by $F\subset\Gamma^0$ a fundamental
set for the action $G\act\Gamma^0$  containing one representative
in each $G$-orbit of blue vertices.
 Join the vertex represented by the element $1$ of $G$
with each vertex in $F$ by a \it red \rm edge. Let
$\eta$ denote this set of red edges. By applying to
this new edges the elements of $G$ we
obtain a $G$-invariant set $G\eta$ of red edges. Put
$\Delta^1=\Gamma^1 \cup G\eta$.

To construct  a new  graph $\widetilde \Delta$ we have to
add a finitely many orbits of new
edges and eventually remove all the
blue edges preserving the connectedness
of the graph. We proceed as follows.

Let $e=(x,y)\in\Gamma^1\subset \Delta^1$ be a blue edge. If
$\Delta\setminus Ge$ is connected then put $\widetilde
\Delta=\Delta\setminus Ge$. Suppose that
$\Delta\setminus Ge$ is not connected and the
endpoints $x$ and $y$ of a blue edge
 belong to different
connected components of $\Delta\setminus Ge$.
We choose
  red vertices $x'$ and $y'$ adjacent to $x$ and $y$ respectively,
  join them with
  a new \it yellow \rm edge $e'$ and put $\widetilde \Delta=(\Delta\setminus Ge)\cup
  Ge'.$

We need to  show that adding or deleting
  the orbit of one edge keeps the properties of   $\Gamma$ valid.
We proceed by induction keeping the notation $\Delta$ for the
graph
  of the previous step for which
  all the requested
   properties were true (at the beginning  $\Delta=\Gamma$); and denote by $\widetilde\Delta$ the
   graph $\Delta\pm Ge$ where $e$ is a red, a yellow or a blue edge   (we add
   the red and the yellow edges and delete the blue ones).
   By construction at least one of the
   vertices  of each new edge $e$   has finite stabilizer  for  $G\act\Delta$,
   so the action on the new graph is still proper and cofinite.
   Since we first add a yellow edge   and then delete
the corresponding
 blue one   the graph $\widetilde\Delta$ remains connected.

To prove  the finess of the new graph we will use a result of
\cite{Bo1} which we briefly state now for the  completeness. A
collection $\mathcal L$ of subgraphs of a graph is
called {\it edge-finite} if for every edge $e$ the set $\{L\in\mathcal L\ \vert\ e\in
\L^1\}$ is finite. We need the following

\begin{lem}\cite[Lemma 2.3]{Bo1}\label{bowlem}. Suppose that $K$ is a fine
graph, and $\mathcal A$ is a collection of arcs of bounded length
in $K$. Then the graph $K[\mathcal A]$ obtained by adding the
edges joining the endpoints of the arcs in $\mathcal A$ is also
fine.
\end{lem}

By this lemma after the adding of the
orbit of an edge to $\Delta$ the obtained
 graph   remains fine. Indeed in the  orbit $\mathcal A=\{gL:g\in G\}$ of a finite
 arc $L\subset \Delta$ the set $\{g\in G\ \vert\  e\in gL^1\}$ is finite.
 Otherwise, since the set $\Delta^1/G$ is finite, there would exist an edge
 of $L$ with infinite stabilizer.

It is obvious that the operation of
deleting of an orbit of blue edges
preserves the finess.

To check the hyperbolicity  let us check that
  there is a map between the set of vertices $\Delta^0$ and $\widetilde\Delta^0$
  of the graphs $\Delta$ and $\widetilde\Delta$ respectively  which is at most
$K$-bilipschitz. Here $K=\max\{r+1, 3\}$ and $r$ is the
diameter of  $F$.
 Indeed
the map  is at most $1$-bilipschitz at the beginning
when we add the first orbit of red edges. Here the direct map is
the identity on $\Gamma^0$ and the inverse map is the projection
of $G$ to $\Gamma^0$ which does not increase the distances ($G$ is
not contained in $\Gamma^0$). On the next steps the direct map is
still an isometry on $\Delta^0$. For the inverse map  to return
from $\widetilde\Delta$ to $\Delta$ we delete the orbit $Ge$ of the edge
$\{x,y\}$ where $y $ is red and $x\in F$ is a blue vertex. Note
that $F$ contains a vertex $x_0$ already joined with $1$, so the
distance between $1$ and $x$ in $\Delta$ is at most $1+r$
asserting that the map is $r+1$-bilipschitz.

  In its turn   adding   of an orbit of yellow or blue edges is
at most $3$-bilipschitz as we replace a path of length at most $3$
by an edge. So on each step there is a quasi-isometry between the
graphs $\Delta$ and $\widetilde\Delta$. Since the number of orbits is
finite the process completes after a finitely many steps
and the final graph is still
 hyperbolic \cite{Gr}.

   To obtain a coned-off Cayley graph from the
final graph $\widetilde\Delta$ it remains  to remove each blue vertex $z$
having finite stabilizer (the blue vertices with infinite
stabilizers will be the parabolic vertices of the coned-off Cayley
graph). We also remove
 all the red edges incident
 to $z$ and join every two vertices  adjacent to $z$
 in $\widetilde\Delta$
 by a yellow edge.
By  the same  argument as above
the obtained graph  is connected, fine and hyperbolic, and
 the $G$-action
  on it is cofinite
and proper on edges.
Let $S$ denote the set of the elements of $G$ joined with the element 1
by yellow edges. It follows that the obtained graph
  is the coned-off Cayley graph
${\rm Cof}(G, S\cup\P)$ where $\mathcal P$ is the set of
parabolic vertices. The lemma is proved. \bx

\medskip

The following proposition finishes the proof of Theorem
\ref{equiv}

\begin{prop}\label{bowos} Let a group $G$ act on a fine
hyperbolic graph $\Gamma$ properly and cofinitely on
$\Gamma^1$. Then $G$ is relatively finitely presented with respect
to any maximal set $\mathfrak P$ of non-conjugate infinite
stabilizers of vertices for the action. Furthermore
it admits a  linear relative Dehn function.
\end{prop}

\proof By Lemma \ref{chgraph} there exists a finite set $S\subset
G$ such that the coned-off Cayley graph ${\rm Cof}(G,
S\cup\mathcal P)$ is fine and hyperbolic. The group $G$ acts on it
cofinitely and properly on edges.

By Lemma \ref{finitecomp}
there exists a simply connected complex
${\mathcal W}={\rm Cof}(G, S\cup\mathcal P)\cup\mathcal C$
where $\mathcal C$ is
 a $G$-finite set of circuits
that bound 2-cells. For every such 2-cell $D$  whose vertices  do
not all belong  to the star of one parabolic vertex we do the
following surgery. Once the boundary $\partial D$  contains  two
consecutive red edges passing through a parabolic vertex $p$
we replace this pair of edges by one yellow edge
and consider the component of $D$ not
containing $p$. We cut in this way  all parabolic vertices on
$\partial D$ and obtain a 2-disk $D'$ whose boundary lies in the
relative Cayley graph. Proceeding similarly with all  2-cells of
$\mathcal W$   we obtain a 2-complex $\mathcal S$ containing a
$G$-finite set of 2-cells attached to circuits in the  graph ${\rm
Cay}(G, S\cup({\cup}\mathfrak P))$.
 It follows that every singular disk in $\mathcal S$ whose boundary is not contained
in the star of a parabolic vertex can be tesselated by a finite
number of the  2-cells $D$ obtained above. So these cells give
rise to a finite relative presentation for $G$ relatively to the
parabolic subgroups $P_i\ (i\in I)$.

To estimate the Dehn function consider  a circuit $\gamma$  in the
complex $\mathcal S$.
We can assume that no three consecutive vertices of $\gamma$
belong to a star of a parabolic vertex.
 We
now replace every yellow boundary edge of $\gamma$
whose endpoints belong to the same coset of a parabolic subgroup
$P_i$ by a pair of red edges in
 ${\rm Cof}(G, S\cup\mathcal P)$ passing through its fixed point $p_i$ ($i\in I$).
 Let $\gamma'$ be the obtained
  circuit  in $\mathcal W$. We have
$\boldsymbol\alpha(\gamma)=\boldsymbol\alpha'(\gamma')$ and
  $\ell'(\gamma')\leq 2\ell(\gamma)$
where $\boldsymbol\alpha$ and $l$ (respectively
$\boldsymbol\alpha'$ and $l'$) denote the relative area  and the
length of a circuit  in the relative Cayley graph (respectively
coned-off graph). By Lemmas \ref{finitecomp} and \ref{chgraph}
$\boldsymbol\alpha'(\gamma')\leq M \ell'(\gamma')$ and so
$\boldsymbol\alpha(\gamma)\leq 2M \ell(\gamma).$

 The proposition
and Theorem \ref{equiv} are proved.\bx

\section{Floyd  metrics and shortcut metrics.}

 \noindent From now on we will assume that $G$ is a finitely generated group acting
 3-discontinuously and 2-cocompactly on a compactum $T.$
 Let us first recall  few   standard definitions
concerning Floyd compactification   (see \cite{F}, \cite{Ka},
\cite{Tu1}, \cite{Ge2}, \cite{GePo1} for more details).

We  will deal with abstract graphs even without assuming any group
action (in particular it can be the Cayley graph or the entourage
graph $\Ga$ considered in Section 3).

Let   $\G$ be a locally finite connected graph.  For a finite path
$\al :I\to \G\ (I\subset \Z)$ we define its \it length \rm to be
$\vert I\vert -1.$ We denote by $d(,)$ the canonical shortest path
distance function on $\G$, and by $B(v, R)$   the ball   at a
vertex $v\in \G^0$
 of radius $R$.

Let $f:\N\to \R_{>0}$ be a function
  satisfying the following conditions :

\begin{itemize}

\item[] $\exists\ {\lambda > 0}\ \forall n\in\N\ :\ 1 < {f(n)\over
f(n+1)} <\lambda\hfill(1)$

\medskip

\item[] $\displaystyle \sum_{n\in \N} f(n) < +\infty.\hfill(2)$

\end{itemize}

\bigskip

Define the Floyd length $\lfv(\al)$ of a path $\al=\al(a,
b)\subset\G $ with respect to a vertex $v$ as follows:

$$\displaystyle \lfv(\al)= \sum_{i} f(d(v, \{x_i, x_{i+1}\})).\hfill\eqno(*) $$

 \noindent where $\al^{0}=\{x_i\}_i$ is the set of vertices of $\al$
 (we assume $f(0):=f(1)$ to make it well-defined).

The  Floyd metric $\d_{f,v}$ is   defined to be the corresponding
shortest path metric:

 $$\displaystyle \d_{f,v}(a,b)=\inf_{\alpha} \lfv(\al),\hfill\eqno(**) $$

\noindent  where the infimum is taken over all   paths $\al$
between the vertices $a$ and $b$ in $\G.$ We denote by $\overline
{\G}_f$   the Cauchy completion of the metric space $(\G,
\d_{f,v})$ and call it {\it Floyd completion}.  Let
$$\DFGa=\overline {\G}_f\setminus \G$$ be its boundary, called {\it Floyd
boundary.}

If $\G$ is a Cayley graph $\Ca(G,S)$ of a group $G$ with respect
to a finite generating system $S$ we denote by ${\overline G}_f$
and by $\partial_f G$ the Floyd completion and the Floyd boundary
respectively. Then the condition (1) above implies that the
$G$-action extends to its Floyd completion $\overline {G}_f$ by
homeomorphisms \cite{Ka}. Therefore in this case for any $g\in G$
the Floyd metric $\d_g$ is the $g$-shift of $\d_1$:

$$\d_g(x, y) =\d_1(g^{-1}x, g^{-1}y),\ \ x,y\in\overline{G}_f, g\in G,$$

\noindent where $1$ is the neutral element of $G$. Every two
metrics $\d_{g_1}$ and $\d_{g_2}$ are bilipshitz equivalent with a
Lipshitz constant depending on $d(g_1, g_2).$ The same properties
are valid for any locally finite, connected and $G$-finite graph
$\G$ ($\vert\G^0/G\vert <\infty$).

\bigskip

Recall   that a {\it quasi-isometric map} (or {\it
$c$-quasi-isometric map}) $\varphi : X\to Y$ between two metric
spaces $X$ and $Y$ is a correspondence such that  :

$${1\over c}d_X(x, y) - c < d_Y(\varphi(x),\varphi (y))\leq cd_X(x, y)+ c,$$

\medskip
\noindent where $c$ is a uniform constant and $d_X,\ d_Y$ denote
the metrics of $X$ and $Y$ respectively.

If in addition $d_X({\rm id}, \psi\circ\varphi) \le {\rm const}$
for a   ($c$-)quasi-isometric map $\psi : Y \to X$  we say that
$\varphi$ is a ($c$-)quasi-isometry between $X$ and $Y.$

 A $c$-quasi-isometric
map $\varphi : I\to X$   is called {\it $c$-quasigeodesic} if $I$
is a  convex subset of $\Z$ or $\R$. A quasigeodesic path $\ga :
I\to \G$ defined on a half-infinite subset $I$ of $\Bbb Z$ is
called \it (quasi-)geodesic ray\rm; a (quasi-)geodesic path
defined on the whole $\Bbb Z$ is called \it (quasi-)geodesic
line\rm.

\medskip

 The following lemma will be often used.

\begin{lem}
\label{agraph} {\rm (Karlsson Lemma)}. Let $\G$ be a locally
finite connected graph. Then for every $\varepsilon
> 0$ and  every $c >0,$  there exists  a finite set $D $ such that $\d_v$-length
of every  $c$-quasigeodesic  $\ga\subset \G$   that does not meet
$D$  is less than $\ve$.\bx

\end{lem}

\medskip

\noindent {\bf Remark.}   {\rm    A.~Karlsson \cite{Ka} proved it
for geodesics in the Cayley graphs of finitely generated groups.
The proof of \cite{Ka} does not use the group action and is also
valid for quasigeodesics.}

\medskip

\medskip

 Consider now a set  $S$ of paths of the form $\ga:[0,n]\to \G$
 of unbounded length  starting at the point $a=\al(0)\in\G$.  Every
 $\ga\in S$ can be considered as an element of the product $\prod_{i\in
I} B(a, i)$. Since $\G$ is a locally finite graph the latter space
is compact in the Tikhonov topology. So every infinite sequence
$(\al_n)_n\subset S$ possesses a ``limit path''
$\delta:[0,+\infty)\to \G$ whose initial segments are initial
segments of $\al_n$.

The following lemma illustrates the properties of limits of
infinite  quasigeodesics of $\G$.

\begin{lem}
\label{geodline}\cite{GePo1} Let $\G$ be a locally finite
connected graph. Then the following
 statements are true:

 \begin{itemize}

  \item[\sf 1)] Every infinite ray $r:[0, +\infty[\to \G$
  converges to a point at the boundary: $\displaystyle \lim_{n\to\infty}r(n)=p\in\partial_f\G.$

\medskip

  \item[\sf 2)] For every point $p{\in}\partial_f\G$ and every $a\in \G$
there exists a geodesic ray joining $a$ and $p$.

\medskip

  \item[\sf 3)] Every two distinct points in $\partial_f\G$ can be joined by a
geodesic line.

\end{itemize}\bx
\end{lem}

Let $\G$ be a locally finite, connected graph on which a finitely
generated group $G$ acts cocompactly. Besides the Floyd metrics
the Floyd completion $\OG$ possesses a set of \it shortcut \rm
pseudometrics which can be introduced as follows (see also
\cite{Ge2}, \cite{GePo1}). Let $\omega$ be a closed $G$-invariant
equivalence relation on $\OG.$ Then there is an induced $G$-action
on the  quotient space $\OG/\omega$. A shortcut pseudometric
$\overline{\d}_g$ is the maximal element in the set of symmetric
functions $\varrho: {\overline \G}_f{\times}{\overline
\G}_f\to\Bbb R_{\geqslant0}$ that vanish on $\omega$ and satisfy
the triangle inequality, and the inequality
$\varrho{\leqslant}\d_g$.

For $p,q{\in}{\overline \G}_f$ the value $\odg(p,q)$ is the
infimum of the finite sums $\displaystyle\sum_{i=1}^n\d_g(p_i,
q_i)$ such that $p{=} p_1$, $q{=}q_n$ and $\tupl{ q_i,
p_{i+1}}{\in}\omega\ (i{=}1,\dots,n{-}1)$ [BBI, p. 77]. Obviously,
the shortcut
 pseudometric $\odg$  is
the $g$-shift of ${\overline\d}_1$. The metrics
$\overline{\d}_{g_1},\ \overline{\d}_{g_2}$ are bilipschitz
equivalent for the same constant as for $\d_{g_1}$, $\d_{g_2}$.

The pseudometric $\odg$ is constant on $\omega-$equivalent pairs
of points of $\partial_f\G,$ so it  induces a pseudometric on the
quotient space $\OG/\omega.$ We denote this induced pseudometric
by the same symbol  $\odg$.

Let $\G$ be a connected, locally finite and $G$-finite  graph. The
graph $\Ga$   given by the discrete system $A=G(\ba_0)\
(\ba_0\in\ent)$ of entourages (see Definition \ref{graph} and
Convention \ref{setA}) is also locally finite, $G$-finite and
connected (Lemma \ref{connect}). So there exists a
$c$-quasi-isometry $\varphi : \G\to\Ga.$ Let $f$ and $g$ be
scaling functions satisfying (1-2) and the condition:
$${g(n)\over f(c n)} < D\ (n\in \N),\hfill\eqno(3)$$

\noindent where $c$ is the above quasi-isometry constant. By
[GePo1, Lemma 2.5] the map $\varphi$ extends to a $G$-equivariant
Lipshitz map between the Floyd completions ${\overline \G}_{f}$
and ${\overline \Ga}_{g}$ of these graphs. We denote this map by
the same letter $\varphi.$ The following lemma is a direct
consequence of the main result of \cite{Ge2}:

\begin{lem}
\label{flmap} (\rm Floyd map). {\it Let  G be a finitely generated
group acting 3-discontinuously and 2-cocompactly on a compactum T.
Then
  there exist  $\mu\in ]0, 1[$  and  a continuous $G$-equivariant
map $F: {\overline \G}_f\to \ti T=A\sqcup T$ for the scaling
function $f(n)=\mu^n.$

Furthermore for every vertex $v\in \G^0$ the quantity
$\od_{\bv}(F(x), F(y))$  is a metric on $\ti T$ where
$x,y\in{\overline \G}_f$ and  $\bv=\varphi(v)=F(v).$}
\end{lem}

\proof It follows from  \cite{Ge2} that  there exists
$\nu\in]0,1[$ and a continuous $G$-equivariant map $\F:
{\overline\Ga}_g\to \ti T$ where $g(n)=\nu^n.$

Let  $\varphi : {\overline \G}_{f}\to {\overline \Ga}_g$  be the
$G$-equivariant Lipshitz map described above where $f(n)=\mu^n$
and $\mu=\nu^{1/c}$. Set    $F=\F\circ\varphi.$ The map   $F$
transfers the pseudometric $\od_v$ on ${\overline \G}_f$ to $\ti
T$ as follows:

$$\od_{\bv}(F(x), F(y))=\od_v(x,y),\ {\rm where}\ \bv=F(v),\ v\in\Ca(G, S).$$

\noindent By \cite{Ge2} each $\od_{\bv}$ is a metric on $\ti T.$
The kernel of $F$ is the closed $G$-invariant equivalence relation
on ${\overline \G}_f$ such that $\od_{\bv}(F(x), F(y))=0$. Indeed
since $\od_{\bv}$  is a metric on $\ti T$ the latter one yields
$F(x)=F(y)\ (x,y\in{\overline G}_f).$

  \bx

\bigskip

 \begin{rems}\label{shortcutvers}  1) We will call the obtained metric  $\overline\d_\bv$
 ($\bv=F(v)\in A$)  on $\ti T$    {\it shortcut (Floyd) metric.}

2)  Lemma \ref{flmap} is in particularly true for any polynomial
scalar function $f$. Moreover one can put  $f=g$ as
$f(cn)/f(n)={\rm const}$ in this case.

3) Since $\odg\leq\d_g$ the Karlsson Lemma \ref{agraph} is also
true when one replaces the Floyd $\d_v$-length by the shortcut
$\odg$-length.
\end{rems}

\section{Horospheres and tubes.}

  Let a finitely generated group $G$ act
$3$-discontinuously and $2$-cocompactly on a compactum $T$. Then
the graph of entourages  $\Ga$ is connected (Lemma \ref{connect}).
We will use the graph distance $d_A$ on $\Ga$ as well as the set
of shortcut metrics $\od_\bv\ (\bv\in \Ga)$ on the compactified
space $\ti T= T\cup A$ coming from Lemma \ref{flmap}  where
$A=\Ga^0.$

We obtain in this Section several properties of tubes and
horospheres which will be used later on.
\begin{lem}
\label{floydsep}  For any integer $k>1$ there exists a constant
$\nu>0$ such that

$$\forall\ \ba,  \bc\in \ti T= T{\sqcup} A,\ \forall \bb\in A\ :\
\ba-\bb-\bc(k)\ {\rm then}\ \od_{\bb}(\ba,\bc)\geq \nu.$$
\end{lem}

\proof For a fixed entourage $\bb\in A$ let $C_{\bb,k}$ denote the
closure of  the set $\{\{\ba, \bc\}\in \widetilde T\times
\widetilde T \ :  \ \ba -\bb -\bc(k)\}$  in $\ti T.$ We first
claim that the set $C_{\bb,k}$ does not intersect the diagonal of
$\ti T\times \ti T.$ Suppose not and $(p, p)\in C_{\bb,k}\cap
\Delta^2\ti T$.   Then there  exist two infinite sequences
$(\ba_n)_n$ and $ (\bc_n)_n$ in $C_{\bb,k}$ converging to $p$. By
discreteness of $A$ we may suppose that $p\in T$. By Lemma
\ref{approximation points} we have $\ba_n-\bb-\bc_n (k)$. Let $U$
be a $\bb$-small neighborhood of $p$. Then $U'$ is $\ba_n$-small
and $\bc_n$-small simultaneously for $n>n_0.$ Hence
$\sh_{\bb}\ba_n\cup\sh_{\bb}\bc_n\subset U$, and so
$\Delta_{\bb}(\sh_{\bb}\ba_n, \sh_{\bb}\bc_n) \leq 1$ which is
impossible. It follows  that $C_{\bb,k}\cap\Delta^2 \ti
T=\emptyset.$

Since $C_{\bb,k}$ is a closed subset of $\ti T\times \ti T,$  and
$\od_{\bb}$ is a metric on $\widetilde T$, there exists a constant
$\nu(\bb)
>0$ such that $\od_{\bb}(\ba,\bc)\geq \nu(\bb)$ on  $C_{\bb,
k}$. Thus our statement holds for the set $C_{\bb,k}$ of
entourages separated by   the fixed entourage $\bb.$

We have  $\displaystyle A= G(\ba_0)$.   If now $\ba-\bb-\bc (k)$
then $\exists~g\in G\ : \bb=g(\ba_0)$,  so
$g^{-1}\ba-\ba_0-g^{-1}\bc(k).$ Thus
$\od_{\bb}(\ba,\bc)=\od_{\ba_0}(g^{-1}(\ba),g^{-1}(\bc))\geq \nu,$
where   $\nu=\nu(\ba_0)$ is the above constant for $\ba_0.$ The
lemma   is proved.\bx

\medskip

The following lemmas  give a local description of
$C$-quasigeodesics around  tubes and horospheres.

\begin{lem}
\label{prqg} There exists a constant $D >0$ such that for every
$C$-quasigeodesic $\ga=\ga(\ba, \bc)$ in $\Ga$ with the endpoints
$\ba, \bc$ we have :

$$\forall\bb\in \Psi_k(\ba,\bc) \ :\ d_{A}(\bb, \gamma) \leq
D,\hfill\eqno(1)$$

\noindent where $\Psi_k(\ba, \bc)=\{\bb\in A\ :\ \ba -\bb
-\bc(k)\}.$
\end{lem}

\proof By Lemma \ref{floydsep} we have $\od_{\bb}(\ba,\bc)\geq
\nu,$ and so the Floyd length $L_{f, \bb}(\ga)$ of $\ga$ is at
least $\nu.$
 By Karlsson Lemma  \ref{agraph} (see also \ref{shortcutvers}.3) there exists a constant $D>0$
 such that $\gamma\cap B(\bb, D)\not=\emptyset$ for the  $d_A$-ball
  $B(\bb, D)$ in $\Ga$
centered at $\bb$ with the radius $D.$ The lemma is proved.\bx

\medskip

\begin{lem}
\label{prhor} The following statements are true :

\begin{itemize}
\item[\sf 1)] For any   $C>0$ and  $E\geq 0$ there exists  $L>0$
such that for any parabolic point $p\in T$ and any
$C$-quasigeodesic $\ga :[0,1]\to\Ga$  one has

 $$d_A(\ga(1), T(p)) \leq E \ \implies d_A(\ga,
\Pi_p(\ga(0))) \leq L\hfill\eqno(2).$$

\item[\sf 2)] There exists a constant $D>0$ such that for any
parabolic point $p\in T$ and any $C$-quasigeodesic $\ga
:[0,\infty[\to\Ga$  one has

 $$\lim_{n\to\infty}\ga(n)=p\ \implies d_A(\ga,
\Pi_p(\ga(0))) \leq D.\hfill\eqno(3)$$

\end{itemize}

\end{lem}


\centerline{
\ifnum\pdfoutput>0%
\pdfximage width 180pt height 150pt{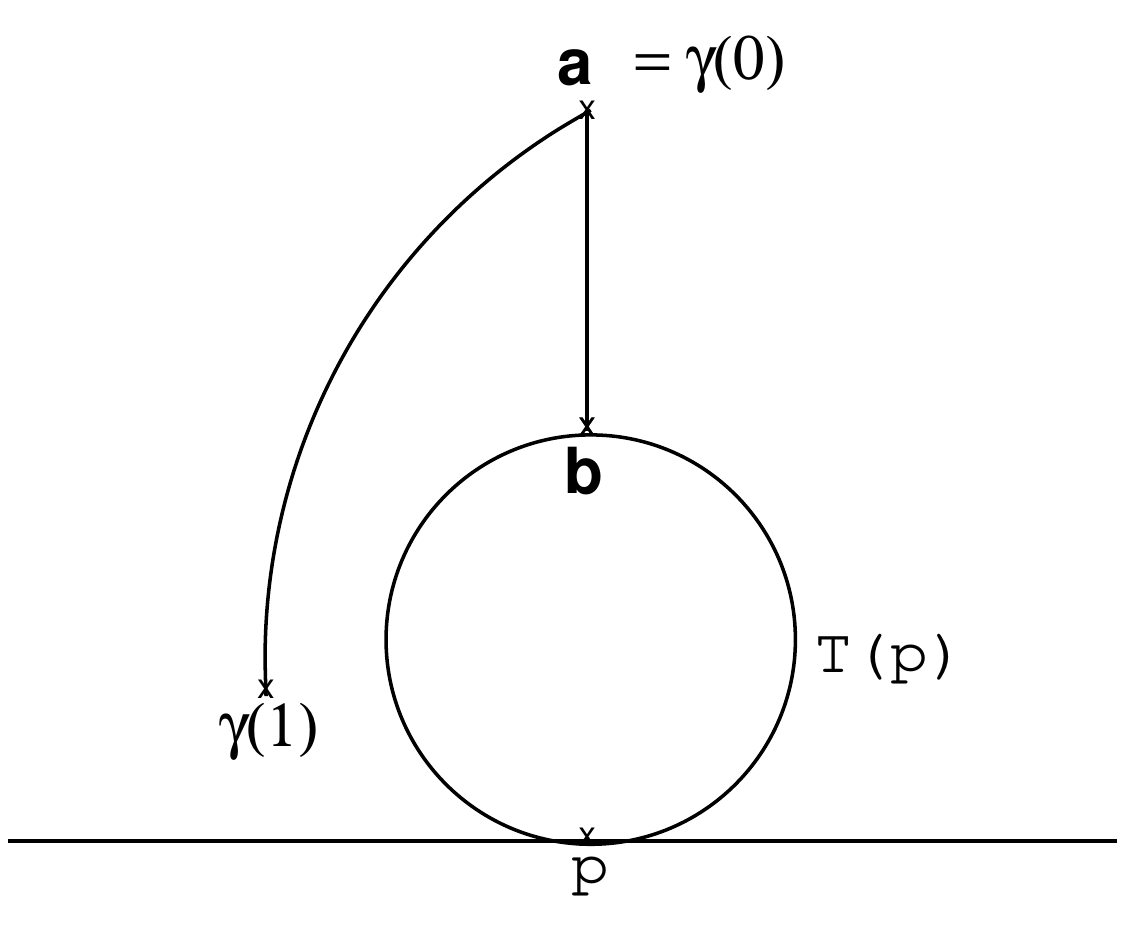}\pdfrefximage\pdflastximage\else\kern324pt\vbox to 166pt{\vss}\fi}
\nopagebreak
\begin{figure}[sh]
\caption{Quasigeodesics around horospheres.}
\end{figure}


\proof 1) Suppose not, then there exist  constants
$C$ and $E$ such that for any  $n$ there exist a parabolic point
$p_n$ and a $C$-quasigeodesic $\ga_n:[0, 1]\to \Ga$ such that
$d_A(\ga_n(1), T(p_n)) \leq E$ and $d_A(\ga_n,
\Pi_{p_n}(\ga_n(0))) > n$ for all $n\in\N.$ By Lemma \ref{geomfin}
 there are at most finitely many
$G$-non-equivalent parabolic points.  So  we may assume that
$p=p_n$ and let $\bb_n\in\Pi_p(\ga_n(0)).$  By the same lemma
  the group ${\rm Stab}_{G} p$ acts cofinitely on
$T(p)$  so we may also suppose that $\bb=\bb_n$.

 Since $d_A(\ga_n(1), \bb)$ is unbounded the set $\{\ga_n(1)\}_n$ is infinite.
 As
$d_A(\ga_n(1), T(p)) \leq E$  by Lemma \ref{unicityacp} up
  to passing to a subsequence we have $\ga_n(1)\to p$. Denote $\ba_n=\ga_n(0)$
  and $\bc_n=\ga_n(1).$
We have $\forall n\ \ba_n\not=\bb$ so $\ba_n\not\in T(p)$ and
$\ba_n-\bb-p.$ By Lemma \ref{approximation points} we obtain
$\ba_n-\bb-\bc_n (n>n_0)$. Thus
 Lemma \ref{prqg} implies that $d_A(\bb, \ga_n) \leq  D$   which is a
contradiction. The statement 1) is proved.

\bigskip

2) We have $\displaystyle\lim_{n\to \infty}(\ga(n)=\bc_n)= p$ and
without lost of generality we can suppose that $\ba=\ga(0)\not\in
T(p).$ Then arguing similarly we obtain $\ba-\bb-\bc_n\ (n>n_0)$
where $\bb=\Pi_p(\ba).$ From Lemma \ref{prqg} we have $d_A(\bb,
\ga) \leq D$.\bx

\bigskip

\noindent The following lemma is a generalization of the previous
one to the geodesics with  variable endpoints:

\bigskip

\begin{lem}
\label{nhor} The following statements are true :

\begin{itemize}
\item[\sf 1)] For any   $C>0$ and  $E\geq 0$ there exists  $M>0$
such that  for any parabolic point $p\in T$ and any
$C$-quasigeodesic $\ga :[-1,1]\to\Ga $  one has

 $$d_A(\{\ga(-1), \ga(1)\}, T(p)) \leq E \ \implies d_A(\ga(0),
\Pi_p(\ga(0))) \leq M\hfill\eqno(2').$$

\item[\sf 2)] There exists a constant $D>0$ such that for any
parabolic point $p\in T$ and any $C$-quasigeodesic $\ga
:[-\infty,+\infty[\to\Ga$  one has

 $$\lim_{n\to\pm\infty}\ga(n)=p\ \implies d_A(\ga(0),
\Pi_p(\ga(0))) \leq D\hfill\eqno(3').$$
\end{itemize}
\end{lem}

\bigskip

\proof 1) As before using the finiteness of $G$-non-equivalent
parabolic points, we fix a parabolic point $p.$ Let
$\gamma_{-}=\gamma([-1,0]),\ {\rm and}\ \gamma_{+}=\ga([0,1]).$ If
$\ba=\gamma(0)\not\in T(p)$ and $\bb=\Pi_p(\ba)$ then by the
statement 1) of Lemma \ref{prhor}  we have $d_A(\ga_{\pm}, \bb)
\leq L$. Let $\bz\in \g_+$ and $\by\in\g_-$ be the points
realizing these distances. Since there is a path from $\bz$ to
$\by$ through $\bb$ of length $2L$,
  the length  $l(\g(\bz, \by)))$ of the C-quasigeodesic $\g(\bz, \by)$ between
  $\bz$ and $\by$ is at most $2L(C+1)$. So at least for one of these
entourages, e.g. $\bz$,   we have  $l(\g(\ba, \bz))\leq L(C+1).$
By  the triangle inequality we obtain $d_A(\ba, \bb)\leq
M=L(C+2).$

\medskip

The same argument and \ref{prhor}.2 imply the statement 2). \bx

\bigskip

\noindent The following Corollary establishes the uniform
quasiconvexity of all horospheres and the quasiconvexity (simple)
of every parabolic subgroup (see also \cite{Ge1} and
\cite{GePo1}).

\begin{cor}
\label{quasiconv} Suppose $G$   acts 3-discontinuously and
2-cocompactly  $T.$ Then there exists $M>0$ such that for every
$p\in\P$ the horosphere $T(p)$ is a $M$-quasiconvex subset of $A$.

Furthermore for every $p\in\P$ there exists a constant $D_p$ such
that the parabolic subgroup $H_p={\rm Stab}_Gp$ is
$D_p$-quasiconvex.
\end{cor}

\proof Suppose first that $\gamma\subset A$ is a $C$-quasigeodesic
with $\partial\gamma\subset T(p)$
 for some $p\in\P$.   By Lemma
\ref{nhor}.1 for   $E=0$ there exists a uniform constant $M>0$
such that   $\gamma\subset N_M(T(p)),$ where $N_M(\cdot)$ denotes
the $M$-neighborhood with respect to the distance $d_A$.

To prove the second part note that since $G$ is finitely generated
it is enough to prove it for the graph $\Ga$ quasi-isometric to
the Cayley graph. By Lemma \ref{geomfin} for every $p\in\P$ the
set $T(p)/H_p$ is finite where $H_p={\rm Stab}_Gp$. So there
exists a constant $E=E(p)$ such that $H\subset N_E(T_p)$ and
$T(p)\subset N_E(H).$ So if $\gamma\subset A$ is a
$C$-quasigeodesic with $\partial\gamma\subset H$ then
$d_A(\partial\gamma, T(p))\leq E$. Then   again by \ref{nhor}.1
there exists a constant $M=M(p)$ such that $\gamma\subset
N_M(T(p))$. So $\gamma\subset N_{D_p}(H)$ where $D_p=M+E.$ \bx
\bigskip

 \noindent {\bf Remark.}  {\rm The above Lemmas \ref{prhor}
and \ref{nhor} are close to some lemmas contained in our work
\cite{GePo1} where the horospheres were defined without using the
entourages. We need the above results in terms of entourages to
apply them in the further argument where the language of
entourages is crucial.}

\medskip

By Proposition \ref{boundness}.3b we have that for every $d>0$
there exists $e=e(d)>0$ such that

$$\forall\ p,q\in\P\ {\rm diam} (N_d(T(p)\cap N_d(T(q)))\leq e.\hfill\eqno(4)$$

\begin{dfn}
\label{horosp}{\rm Let $\ga\subset \ti T$ be a $C$-quasigeodesic.
We call an  entourage $\bv\in \ga$  $d$-{\it horospherical} if
there exist parts $[\bv, \bc]$ and $[\ba, \bv]$ of $\ga$ of length
greater than the constant $e$ and  which are contained in a
$d$-neighborhood  $N_d(T(p))$ of  a horosphere $T(p)$.

The entourage $\bv\in \ga$ is called {\it non-horospherical} in
the opposite case.}
\end{dfn}

\medskip

\noindent {\bf Remark.} By (4) we can  suppose that the parabolic
point $p$ with respect to which the (non)-horosphericity is
considered is unique.

\medskip

\begin{lem}
\label{voistube} Let   $\gamma=\gamma(\ba, \bc)$  be a
$c$-quasigeodesic. Suppose that   $P=P(\ba, \bc)$ is a
non-refinable tube   having the same ending vertices $ \ba$ and
$\bc$ as $\ga$. For every sufficiently large $d>0$  there exists a
constant $E>0$  such that $d_A(\bg, P) \leq E$ for every
$d$-non-horospherical point $\bg\in \ga.$

\end{lem}

 \proof  Note that the non-refinable tube $P(\ba, \bc)$
 exists by Proposition
\ref{exent}. By Lemma \ref{prqg} there exists $D>0$ such that
  for every $\bp_i\in P$ we have $d_A(\bp_i, \ga) \leq D\ (i=1,...,m)$. So let us fix a
 non-horospherical entourage $\bg\in
 \ga$,  and let $\bg_i\in\ga$ be such that
 $d_A(\bp_i,\ga)= d_A(\bp_i, \bg_i)\ (i=0,...,m)$. Let us also assume that
 $\bg\in \ga(\bg_i, \bg_{i+1})$ where $\ga(\bg_i, \bg_{i+1})$ denotes the
 part of $\ga$ between $\bg_i$ and $\bg_{i+1}$.

By Corollary \ref{constbound}
 there exists a constant $C>0$ such that if $d_A(\bp_i,
 \bp_{i+1})>C $ then the pair $\{\bp_i, \bp_{i+1}\}$ is contained in a
 horosphere $T(p)$. In this case $\{\bg_i, \bg_{i+1}\}\subset N_D(T(p))$ and  by Lemma \ref{nhor}
we have that $\ga(\bg_i,\bg_{i+1})\subset N_L(T(p))$ for some
  $L=L(D)>0.$ Let $d$ be any number bigger than $L.$
If $\bg$ is $d$-non-horospherical then by \ref{horosp}
  $d_A(\bg,\bg_i)$ or $d_A(\bg,\bg_{i+1})$ is less
than $e$. Thus $d_A(\bg, P)\leq e+d.$

If now $\da(\bp_i, \bp_{i+1}) \leq C$ then
$d_A(\bg_i,\bg_{i+1})\leq c(C+2D)+c.$ So $\da(\bg, P)\leq
\da(\bg_i, \bg)+D\leq c(C+2D)+c+D.$

Put $E=\max\{e+d, c(C+2D)+c+D\}$. The lemma  is proved. \bx

\bigskip

\noindent {\bf Remark.} The constants $d$ and $e$ depend on the
constants $D$, $C$ and $L=L(D)$ given respectively by the
statements \ref{prqg}, \ref{constbound} and \ref{nhor}.

\section{Tight curves in $\Ga.$}

  Let a finitely generated group $G$ act
$3$-discontinuously and $2$-cocompactly on a compactum $T$.
 For a parabolic point
$p$   we denote by $N(T(p))$ a  neighborhood of the horosphere
$T(p)$ in the graph $\Ga$ (see Section 3.2). The notation ${\rm
diam}(\cdot )$ is used for the diameter of a set with respect to
the distance $d_A$ and $\vert\cdot\vert$ stands for the length of
a curve. We denote by $c^{-1}(n)$ the linear function
$\displaystyle {n\over c} - c$
  for some constant $c > 0.$

\begin{dfn} \label{tcurve}
{\rm For   positive integers  $l$ and $c$, a curve $\ga : I\to
\Ga$ is called {\it $(l,c)$-tight} (or just {\it  tight } when the
values of $l$ and $c$ are fixed) if for every $J\subset I$ the
following conditions hold:

\begin{itemize}

 \item[\sf 1.] $ \vert J\vert \leq l\implies
 \gamma\vert_{J}$ is a $c$-quasigeodesic.

 \medskip

 \item[\sf 2.] If  $\vert\ga(J)\cap N(T(p))\vert >  l$ for some
 $p\in\P$ then ${\rm diam} (\ga(\partial J)) > c^{-1}(l).$\bx

 \end{itemize}

 }
 \end{dfn}

\noindent The rest of the Section is devoted to the proof of the
following Theorem describing the non-horospherical points (see
Definition \ref{horosp}) of tight curves.

\bigskip

\noindent {\bf Theorem B.} {\it For every   $c>0$ and $d>0$ there
exist positive constants $l_0, w_0,  c_0$ such that for all $l\ge
l_0$ and  every $(l, c)$-tight curve $\ga\subset \Ga$
 there exists   a $c_0$-quasigeodesic $\al\subset A$ such that  every
$d$-non-horospherical vertex of $\ga$ belongs to the $w_0
$-neighborhood $N_{w_0}(\al)$ of $\al.$}

\bigskip

The following three lemmas are close  to the results of the
previous Section. We use below the notation ${\rm diam}_{\od_\bv}$
 for the diameter of a set with respect to the shortcut
metric $\od_\bv\ (\bv\in A)$ on $\ti T$  (see Lemma \ref{flmap})

\begin{lem}
\label{geodsep} There exist  positive  constants    $\rho$ and  $d
$ such that for every $c$-quasigeodesic  $\ga :I\to \Ga$ of
non-zero length and a $d$-non-horospherical point $\ga(0)\in \Ga$
one has:

$${\rm diam}_{\od_{\ga(0)}}(\ga(\partial
I)) >\rho.$$

\end{lem}

\proof Let us first prove that there exists a constant $r>0$ such
that for some $\rho=\rho(r)$ we have  $$\da(\ga(0), \ga(\partial
I))> r\ \Longrightarrow\ \od_{\ga(0)}(\ga(\partial I))
>\rho\hfill\eqno(*).$$

\noindent Suppose not. Then for every $d >0$ there exists a
sequence of quasigeodesics $\ga_n$ such that $\da(\ga_n(0),
\ga_n(\partial I))\to+\infty$ and $\od_{\ga_n(0)}(\ga(\partial
I))\to 0$ where $\ga_n(0)$ is a $d$-non-horospherical point of
$\ga_n.$

Up to choosing a subsequence we may suppose that the sequence
$(\ga_n)_n$ converges in the Tikhonov topology to a
$c$-quasigeodesic $\ga:\Z\to \Ga$ such that
$\displaystyle\lim_{n\to\pm\infty}\ga(n)=p\in T.$ Then $\ga$ is a
horocycle at $p$ and by [GePo1, Lemma 3.6] the point $p$ is
parabolic. By Lemma
 \ref{nhor}.2 for every $i\in\Z$ the distance $\da(\ga_n(i), T(p))$ is uniformly
 bounded by a constant $D > 0$. So the points $\ga_n(0)$ are $D$-horospherical for
 sufficiently large $n.$ The obtained  contradiction proves (*).

We are left now with the case when $\da(\ga(0), \ga(\partial
 I))\leq
 r$ where the constant $r$ satisfies (*).
  Suppose first that the distance between  $\ga(0)$ and both   endpoints of
  $\ga(\partial I)$ is
less than $r$. By translating $\ga(0)$ to a fixed basepoint
$\bv\in A$  we obtain that $\ga$ is contained in a finite ball
$B(\bv, r+c(r))$. Then the $\od_\bv$-length of $\ga$ is uniformly
bounded from below. If the distance
 between $\ga(0)$ and only one of its endpoints is bigger than $r$
 then   the $\od$-length
 of $\gamma$ is still bounded from below.

 Denoting by $\rho$ the minimum among all of
 these constants we obtain the lemma.  \bx

\bigskip

\noindent{\bf Remark.} Above we have used  Lemma 3.6 from
\cite{GePo1} stated there for the Cayley graphs. Since our graph
$\Ga$ is quasi-isometric to the Cayley graph this result   can be
applied.

\bigskip

 Recall that
$\displaystyle A=  G(\ba_0)$ is the vertex set of the graph $\Ga.$
Using a "refining" procedure we will  now introduce a new graph
$\Ga^*$ whose vertex set $A^*$ satisfies some additional
conditions.

From now on we fix the constant  $d$ and  $\rho=\rho(d)$ coming
from Lemma \ref{geodsep} and an integer  $k > 3$ which will be
used in the betweenness relation below. Let $\delta$ be a number
such that
$$\displaystyle 0 < \delta < {\rho\over k+2}.\hfill\eqno(**)$$

\bigskip

\noindent \underline{\it Definition of the set $A^*$} :  For every
$\bv\in A$ denote by  $\bv^* $ the entourage  $\{\{x, y\}\in \bs2\
:\ \od_\bv (x, y)<\delta\}.$

\bigskip

 It follows from the
following lemma that the compactifying topology on $T$ coming from
the  graphs $A^*$ and $A$  is the same.

\bigskip

\begin{lem}
\label{convA^*}{\rm

  $\forall p\in T\ \ba_n\to  p $ if and only if $\ba^*_n\to p.$}
\end{lem}

\proof Suppose first that $\ba_n\to p$ and $\ba^*_n\not\to p$.
Then there exists a neighborhood $U_p$ of the point $p$ such that
$U'_p$ is not $\ba^*_n$-small for $n>n_0.$ So $\exists\ \zx_n,\
\by_n\in U'_p\ :\ \od_{\ba_n}(\zx_n, \by_n) >\delta.$ It follows
that up to subsequences we have $\zx_n\to x\in T,\ \by_n\to y\in
T\ (n\to\infty)$ and $x\ne y\ne p\ne x.$ Let $U_x$ and $U_y$ be
closed neighborhoods of $x$ and $y$ such that $U_p\cap U_x\cap
U_y=\emptyset.$

Let $H(U_{x,y})\subset \Ga$ denote the set of geodesics whose
endpoints are situated in $U_{x,y}=U_x\sqcup U_y.$ By [GePo1, Main
Lemma]   $\overline {H(U_{x,y})}\cap T= \overline {U_{x,y}}\cap T$
where $\overline U_{x,y}$ means the closure of $U_{x,y}$ in $\ti
T=A\sqcup\ent.$ It follows that  the geodesics $\ga_n(\zx_n,
\by_n)\subset \Ga$ between $\zx_n$ and $\by_n$ do not intersect a
neighborhood $V_p\subset U_p$ of $p\ (n>n_0).$  Since $\ba_n\to p$
we have $\da(\ba_n, \ga_n)\to\infty.$ By Karlsson Lemma
\ref{agraph}  (see also Remark \ref{shortcutvers}.3) we obtain
that $\od_{\ba_n}(\zx_n, \by_n) < \delta\ (n>n_0)$ which is a
contradiction.

Suppose now $\ba^*_n\to p$ and $\ba_n\not\to p.$ Then up to a
subsequence we have $\ba_n\to q\not=p.$ Let $U_p$ be a
neighborhood of $p$ such that $U'_p$ is $\ba^*_n$-small ($n>n_0$).
We have $d_A(\ba_n, U_p)\to +\infty$. Then by Karlsson Lemma
$\forall\ x, y\in U_p\ \od_{\ba_n}(x,y) <\delta\ (n>n_0).$ So
$U_p$ and $U'_p$ are both $\ba^*_n$-small ($n>n_0$) which is
impossible (see (1) of 3.1). \bx

\bigskip

\noindent The  need of the graph $A^*$ is explained by  the
following:

\begin{lem}
\label{geodtube}

There exists  constant   $w>0$  such that for  every quasigeodesic
$\ga:I\to \Ga$ containing three vertices $\ba, \bb, \bc\in A$ the
following is true:

$$\bb\ {\rm is\ } {d-{\rm non}-}{\rm horospherical}\ \wedge\
 d_A(\bb, \{\ba, \bc\})
> w\ \Longrightarrow\ \ba^*-\bb^*-\bc^* (k).$$
\end{lem}

\proof Suppose not and there are sequences $\ba_n, \bc_n$ and
$\bb_n$ such that $\bb_n$ is $d$-non-horospherical, $\da(\bb_n,
\{\ba_n, \bc_n\})\to\infty$ and $\ba_n^*-\bb_n^*-\bc_n^* (k)$ is
not true. Since $A$ is $G$-finite we can suppose that $\bb_n=\bb.$
Up to a subsequence we have $\ba_n\to p,\ \bc_n\to q$. Let
$\ga_n=\ga_n(\ba_n, \bc_n)\subset \Ga$ be a geodesic between
$\ba_n$ and $\bc_n$. Since $\bb$ is non-horospherical we have by
Lemma \ref{geodsep} that $\od_{\bb}(p, q) >\rho,$ hence $p\ne q$.

Let $U_p$ and $U_q$ be disjoint $\bb^*$-small neighborhoods of $p$
and $q$ respectively. So  $\od_{\bb}(U,V)
> \rho - 2\delta,$ and (**) yields $\od_{\bb}(U_p, U_q)
\geq \rho- 2\delta
>   k\cdot \delta $. We obtain  $\Delta_{\bb^*}(U_p, U_q)>k.$ By Lemma \ref{convA^*} we  also have
$\ba^*_n\to p$ and $\bc^*_n\to q.$ So $U'_p$ and $U'_q$ are
$\ba^*_n$-small and $\bc^*_n$-small respectively ($n>n_0).$ Hence
$U_p\supset \sh_{\bb^*}\ba^*_n$ and $U_q\supset
\sh_{\bb^*}\bc^*_n.$ It follows that
$\Delta_{\bb^*}(\sh_{\bb^*}\ba^*_n, \sh_{\bb^*}\bc^*_n) \geq
(\Delta_{\bb^*}(U_p, U_q) > k$. Therefore $\ba^*_n-\bb^*-\bc^*_n
(k)$ which is a contradiction. \bx

\bigskip

 \begin{lem}
\label{convhor}

For every $d>0$ there exists a constant  $l_0$    such that for
every parabolic point $p$,
 and all entourages $\bb, \bc, \bd\in N_d(T(p))$, and $\ba\in A$ one has

$$\forall\ l
> l_0\ \ :\ \da(\bb, \bc) > l\ \wedge \ \da(\bb, \bd) > l\ \wedge\ \ba^*-\bb^*-\bc^*\  (k)\ \Longrightarrow \
\ba^*-\bb^* - \bd^*\ (k-1)\hfill\eqno (1)$$

\end{lem}

\proof Since by Lemma \ref{geomfin} the set of parabolic points is
$G$-finite it is enough to prove the statement for a fixed
parabolic point $p\in T.$ By Lemma \ref{unicityacp} the
 parabolic point $p$
is the unique limit point of $N_{d}(T(p)).$
 By definition of the topology
 on $T{\sqcup}\mathsf{Ent}T$ for sufficiently large
 $l_0$ our assumption implies that
 the entourages $\bc$ and $\bd$ are sufficiently close to $p$. By
 Lemma \ref{convA^*} the entourages
 $\bc^*$ and $\bd^*$ are also  close to $p.$
 So for every
   $\bb^*$-small neighborhood $U_p$ of $p$ its complement  $U'_p$ is
 $\bc^*$-small and $\bd^*$-small for $l>l_0.$ Then $\sh_{\bb^*}\bc^*\subset U_p$ and
 $\sh_{\bb^*}\bd^*\subset U_p$. Therefore
 $\ti\Delta_{\bb^*}(\sh_{\bb^*}\bc^*, \sh_{\bb^*}\bd^*)\leq 1.$
  We obtain $\Delta_{\bb^*}(\sh_{\bb^*}\bd^*, \sh_{\bb^*}\ba^*) \geq
  \Delta_{\bb^*}(\sh_{\bb^*}\ba^*, \sh_{\bb^*}\bc^*)-
  \ti\Delta_{\bb^*}(\sh_{\bb^*}\bc^*, \sh_{\bb^*}\bd^*)  >
  k-1$.
 \bx

\bigskip
\begin{rem} (about the constants).
\label{allconst} {\rm Since now on we   assume that the tightness
constant $l$ is much bigger than the horosphericity
 constants $d$, $e=e(d)$ (see Definition
 \ref{horosp} and the Remark after it) and $w$ (see \ref{geodtube}). We will also suppose that the chosen constants satisfy
  the following
 relations: $$l_0 > 4w,\ w > e.$$
}
  \end{rem}

\bigskip

 {\it Proof of Theorem B.} Recall that for a fixed
constant $d>0$ by Lemma \ref{geodsep} we have found $\rho=\rho(d)$
 and have defined the set $A^*$ of vertices of a new graph of entourages.
Since now on the term "(non)-horosphericity" will  mean
"$d$-(non)-horosphericity".

Before going into the details we   outline the proof of the
theorem. We start by choosing  non-horospherical points $\bv_n$ of
the curve $\gamma$ which give by Lemma \ref{geodtube} an auxiliary
tube $P^*=\{\bv_n^*\}$ in the graph $A^*$. There is a
quasi-geodesic $\al^*\subset A^*$ whose non-horospherical points
are in a bounded distance from $P^*$ (Lemma \ref{voistube}). Since
the graphs $\Ga$ and $\Ga^*$ are $G$-finite the map $\varphi :
\bv\to \bv^*$ is a quasi-isometry between them.   This will give
us a quasi-geodesic $\al\subset A$ satisfying the statement of the
Theorem. All the remaining constants will be found in the course
of the proof.

To construct the tube $P^*$  we  proceed inductively by choosing
vertices of $\ga$ as follows. Let $\ga(0)$ be the first
non-horospherical point on $\ga$, then we put $\bv^*_0=\ga^*(0)$.
Suppose that a  point $\bv^*_n=\ga^*(n)$ is already chosen. Then
for the constant $w$ fixed above we choose $i_{n+1} \geq i_n+w$
such that $\ga(i_{n+1})$ is the first non-horospherical point on
$\ga$ after $\ga(i_n+w).$ We set $\bv^*_{n+1}=\ga^*(i_{n+1}).$ The
following proposition shows that for every $n$ each three chosen
neighboring vertices form a tube $\bv^*_{n-1}-\bv^*_n-\bv^*_{n+1}\
(k-2)$ for the integer $k$ fixed above. Then all the constructed
vertices will give a tube $P^* = \bv^*_0-\bv^*_1-...-\bv^*_m\
(k-2)$.

\begin{prop}
\label{3tube}
 For every $n\in \N$ one has $\bv^*_{n-1} - \bv^*_n - \bv^*_{n+1}\  (k-2).$
 \end{prop}

 {\it Proof of the proposition.} There are four different
cases depending on the lengths $\vert\ga\vert_{[i_n,
i_{n+1}]}\vert= i_{n+1}-i_n$ of the parts of $\ga\ (n\in\N).$

\bigskip

 \noindent {\it Case 1.} $\displaystyle i_n-i_{n-1} \leq l/2
 \wedge\   i_{n+1}-i_n \leq l/2,$

By definition of a tight curve the points  $\ga(i_{n-1}),
\ga(i_n), \ga(i_{n+1})$ belong to a $c$-quasigeodesic part of
$\ga$ so the result follows from Lemma \ref{geodtube}.

\bigskip

\noindent {\it Case 2.} $\displaystyle i_{n+1}-i_{n-1} > l.$

\medskip

There are three subcases.

\medskip

\noindent {\it Subcase 2.1.} $\displaystyle
 i_n-i_{n-1}\leq l/2 \
 \wedge\ i_{n+1}-i_n  > l/2,$

Since $\ga(i_{n+1})$ is the first non-horospherical point   on
$\ga$ after $\ga(i_n+w)$  and $w < l/2$ the point $\ga(i_n+w)$ is
horospherical.  Since $w>e$ by   the   Remark after \ref{horosp}
there exists a unique horosphere $T(p)$ such that $\da(\ga(i_n+w),
T(p))\leq d$. As $\ga\vert_{[i_n, i_n+w]}$ is a $c$-quasigeodesic
we have

$$\da(\ga(i_n), T(p)) < cw+c+d.\hfill\eqno(***)$$
 Furthermore
Lemma \ref{geodtube} yields:

$$\ga^*(i_{n-1})-\ga^*(i_n)-\ga^*(i_{n-1}+l)\ (k).\hfill\eqno(2)$$

\noindent   Since $ i_{n+1}- i_{n-1}  > l$ the point $\ga(i_{n-1}
+ l)$ is  also horospherical and
  $\displaystyle \ga(i_{n-1}+l)\in \ga\vert_{]i_n+w, i_{n+1}]}.$ The curve $\ga\vert_{[i_n, i_n+l]}$
  is still $c$-quasigeodesic so we have
  $$\da(\ga(i_n), \ga(i_{n-1}+l))  >
  {i_{n-1}+l- i_n\over c}-c
  \geq  {l\over 2c} - c > {l\over 4c},\hfill\eqno(3)$$

\noindent where we assume that  $l>   l_0 > 4c^2$ for the constant
$l_0$ from Lemma \ref{convhor}.

\vskip-10pt
\noindent\centerline{\kern180pt
\ifnum\pdfoutput>0%
\pdfximage width 486pt height
134pt{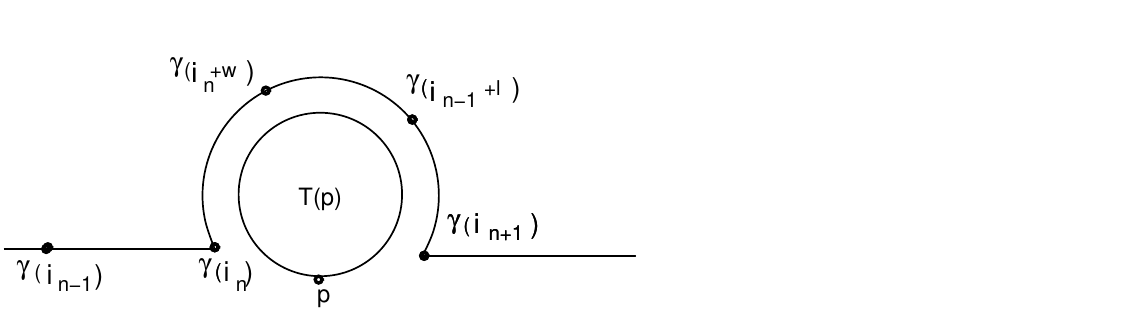}\pdfrefximage\pdflastximage\else\kern270pt\vbox to 134pt{\vss}\fi}

\centerline{Figure 3: Tight curves around horospheres.}

\bigskip

By construction we can also suppose that $\ga(i_{n+1})\in
N_d(T(p))$ for the $d$-neighborhood $N_d(T(p))$ of $p$. Indeed
otherwise there would exist another non-horospherical point on
$\ga$ after $\ga(i_n+w)$ and preceding $\ga(i_{n+1})$. So by (***)
$\{\ga(i_n), \ga(i_{n+1})\}\subset N_{d_0}(T(p)),$ where
$d_0=cw+c+d.$

  If, first,
  $ i_{n+1}-i_n\leq l$ then $\ga\vert_{[i_n,  i_{n+1}]}$   is a
$c$-quasigeodesic, and $\displaystyle \da(\ga(i_n), \ga(i_{n+1}))
> l/2c -c > l/ 4c.$ Hence by the choice of $l_0$ (see Remark \ref{allconst})   and all
$l>l_0$ we have  from (2), (3) and Lemma \ref{convhor}

   $$\ga^*(i_{n-1})-\ga^*(i_n)-\ga^*(i_{n+1})\ (k-1).\hfill\eqno(4)$$

  If now $ i_{n+1}-i_n
> l$ then applying   \ref{tcurve}.2   to $ N_{d_0}(p)$ we obtain $\da(\ga(i_n),
\ga(i_{n+1})) > c^{-1}(l) $ and again  (4) follows from (2), (3)
and  Lemma \ref{convhor}. \bx

\medskip

\noindent {\it Subcase 2.2.} $\displaystyle i_n-i_{n-1}\geq l/2 \
 \wedge\  i_{n+1} - i_n \leq l/2,$

\medskip

The argument is similar to that of Subcase 2.1 but it works in the
opposite direction. We have the tube $\ga^*(i_{n+1}) - \ga^*(i_n)
- \ga^*(i_{n+1}-l)\ (k).$ As above  if $i_n-i_{n-1}\leq l$ then
the curve $\ga\vert_{[i_{n-1},  i_{n}]}$ is $c$-quasigeodesic and
so its diameter is greater than $l/4c$. If not then using the
tightness property of it, we obtain that $\da(\ga(i_{n-1}),
\ga(i_{n}))
> c^{-1}(l) $ and (4) follows by the same argument as in Subcase 2.1.

\medskip

\noindent {\it Subcase 2.3.} $\displaystyle i_n-i_{n-1} \geq l/2 \
 \wedge\ i_{n+1} - i_n\geq l/2,$

\medskip

In this case we have that   the points $\ga(i_n - l/ 4)$ and
$\ga(i_n + l/ 4)$ preceding respectively  $\ga(i_n)$ and
$\ga(i_{n+1})$ are both horospherical. Indeed  $w < l/4$ and
$\ga(i_n)$ and $\ga(i_{n+1})$ are the first non-horospherical
points after
 $\ga(i_{n-1})$ and $\ga(i_n)$ respectively. So we can suppose that
 $\ga(i_n)\in N_d(T(p))$ and
$\ga(i_{n+1})\in N_d(T(q))$ where $p$ and $q$ are distinct
parabolic points. Since $\ga\vert_{[i_n-l/4, i_n+l/4]}$ is a
quasigeodesic
  by Lemma \ref{geodtube} we have
$$\ga^* (i_n-l/4) - \ga^*(i_n) - \ga^*(i_n+l/ 4 )\
(k).\hfill\eqno(5)$$

\noindent We also have $\da(\ga(i_{n-1}), \ga(i_n))$ and
$\da(\ga(i_n), \ga(i_{n+1}))$ are both greater than $l/4c$. Indeed
if $i_n-i_{n-1} > l$ then by $(l,c)$-tightness we have
$\da(\ga(i_{n-1}), \ga(i_n)) > c^{-l}(l)>l/4c.$ If  $i_n-i_{n-1}
\leq l$ then $\ga\vert_{[i_{n-1},  i_{n}]}$ is $c$-quasigeodesic,
and    as above  $\da(\ga(i_{n-1}), \ga(i_n)) > l/4c$. In the same
way we obtain  $\da(\ga(i_{n}), \ga(i_{n+1}))
> l/4c.$

Applying now Lemma \ref{convhor} to (5) two times for $l
> 4c l_0$ we  obtain

$$\ga^*(i_{n-1}) - \ga^*({i_n}) - \ga^*(i_{n+1})\ (k-2).$$

\noindent The proposition is proved.\bx

\bigskip

We  continue the proof of Theorem B.   By Proposition \ref{3tube}
the curve $\ga$ admits a set of non-horospherical points
$\bv_n=\ga(i_n)$ such that $\bv^*_n=\varphi(\ga(i_n))$ is a vertex
of the tube $P^*.$ Let $\bu=\ga(i)$ be a non-horospherical point
of $\ga$ which does not belong to the set $\{\bv_n\}_n.$  Then by
construction $i_n \leq i < i_n+w$ for some $i_n\in \{0, ...,m\}.$
Since $w < l$ the curve $\ga\vert_{[i_n, i_n+w]}$ is a
$c$-quasigeodesic so $\da(\bv_n, \bu) \leq cw+c.$ The map $\varphi
: \bu\in A\to \bu^*\in A^*$ is a quasi-isometry so $d_{A^*}(\bu^*,
\bv_n^*) \leq w_1$ for some uniform constant $w_1
>0.$ Let $\al^*$ be a geodesic in the graph $\Ga^*$ with the same
endpoints as $P^*$. Then by Lemma \ref{prqg} (applied to the graph
$\Ga^*$) there is a constant $D^*
>0$ such that $\forall\ \bv^*\in P^*\ :\ d_{A^*}(\al^*, \bv^*) \leq D^*.$
So for every non-horospherical point $\bu\in \ga$ we have
$d_{A^*}(\bu^*, \al^*) \leq d_{A^*}(\bu^*, \bv^*)+d_{A^*}(\bv^*,
\al^*)\leq w_1 + D^*$ where $\bv^*\in P^*.$ The map $\varphi^{-1}
: \bu^* \to \bu$ is a quasi-isometry too. Hence $\al
=\varphi^{-1}(\al^*)$ is a $c_0$-quasi-geodesic in $\Ga$ such that
for every non-horospherical point $\bu\in \ga$ we have $\da(\bu,
\al) \leq w_0$ for some positive constants $c_0$ and $w_0.$
Theorem B is proved. \bx.

\section{Floyd quasiconvexity of parabolic subgroups.}

Let $G\act T$ be a 3-discontinuous and 2-cocompact action of a
finitely generated group $G$  on a compactum $T.$ Let $\G$ be a
locally finite, connected graph on which   $G$ acts
discontinuously and cofinitely (e.g. its Cayley graph or the graph
of entourages). We denote by $d(,)$ the graph distance of $\G$.
Let $f:\N\to \R_{>0}$ be a scaling function esatisfying the
following conditions (1-2) (see Section 4):
\begin{itemize}

\item[] $\exists\ {\lambda > 0}\ \forall n\in\N\ :\ 1 < {f(n)\over
f(n+1)} <\lambda\hfill(1)$

\medskip

\item[] $\displaystyle \sum_{n\in \N} f(n) < +\infty.\hfill(2)$

\end{itemize}

To precise that $f$ satisfies (1) with respect to some $\lambda
\in ]1, \infty[$ we will say that the function $f$ is {\it
$\la$-slow}. Denote by $\d_f$ the corresponding Floyd metric on
$\G$ with respect to a fixed vertex $v\in\G^0$.

 By a standard
argument based on Arzela-Ascoli theorem it follows that the Floyd
completion $\overline {\G}_f$ of the graph $\G$ is a geodesic
(strictly intrinsic) space (see e.g. [BBI, Theorem 2.5.14]. We
call Floyd geodesic (or $\d_f$-geodesic) a geodesic in the space $
\overline {\G}_f$ with respect to the Floyd $\d_f$-metric. The
geodesics in $\G$ with respect to the graph distance $d$ we call
below $(d-)$geodesics.

 The set $\G^0/G=K$ is finite so we can
identify in $\G$ a subgroup $H$ of $G$ with the orbit
$HK=\bigcup_{h\in H}hK\subset \G^0$. Let $N_R(H)$ denote the
$R$-neighborhood of $HK$ in $\G$ for the graph metric.

\bigskip

\begin{dfn}\label{fqconv}
{\rm Let $\G$ be a locally finite, connected graph possessing a
$G$-finite action. A subgroup $H$ of $G$ is called {\it Floyd
quasiconvex} in $\G$ if there exists a constant $R=R(H) >0$ such
that every Floyd geodesic $\ga=\ga(h_1, h_2)\subset \G$ for the
metric $\d_f$ having the endpoints $h_i$ in $H$ belongs to
$N_R(H)$: $\forall x\in\ga\ :\ d(x, H) < R.$}
\end{dfn}

\bigskip

 By Corollary \ref{quasiconv}   every parabolic subgroup of $G$ is quasiconvex
 with respect to the word metric (see also \cite{Ge1}). The aim of
this Section is to prove the following Theorem stating the Floyd
quasiconvexity of parabolic subgroups.
\bigskip

\bigskip

\noindent {\bf Theorem C.} Let  $G$ be  a finitely generated group
acting $3$-discontinuously and $2$-cocompactly on a compactum $T.$
 Let $\G$
be a locally finite, connected graph  admitting a cocompact
discontinuous action of $G.$ Then there exists a constant
$\la_0\in ]1, \infty[$ such that  for every $\lambda\in ]1,
\la_0[$ and every  $\lambda$-slow Floyd scaling function $f$
satisfying (1-2), each parabolic subgroup $H$ of $G$ is Floyd
quasiconvex for the Floyd metric $\d_f$.\bx

\bigskip

\noindent  We start with two lemmas.

\begin{lem}
\label{fwgeod} For every $r>0$ there exists $\lambda_0 >1$ such
that $\forall\lambda\in ]1, \lambda _0[$ and every $\lambda$-slow
function $f$ the condition $d(x,y) \leq r\ (x, y\in \G^0)$ implies
that every Floyd $\d_f$-geodesic $\ga=\ga(x,y)\subset\G$ whose
endpoints are $x$ and $y$ is a geodesic in $\G$.
\end{lem}

\bigskip

\noindent {\bf Remark.} A similar statement for
$\delta$-hyperbolic spaces   is proved in [Gr, Lemma 7.2.1]

\bigskip

  \proof Let $v\in \G^0$ be a basepoint. Denote by $\om=\om(x,y)$ a $\G$-geodesic between $x$
  and $y$ for which $\vert \om\vert=r.$ Let   $m\in\omega$ such that $d(v,m)=d(v,\omega).$
  Then for at least one of the points $x$ or $y$, say $x$, we have $d(v,x)\geq d(m,x)$. Indeed otherwise
  $d(x,y)\leq d(v,x)+d(y,v) < r$ which is impossible.
  Put 
  $R=d(v,x).$ We have  $\displaystyle L_f(\om)=\sum
_{i=1}^r f(d(v, \{x_i, x_{i+1}\}))\leq r f(d(v,m)).$

Suppose by contradiction that $\ga$ is not $d$-geodesic and so
$\vert \ga\vert \geq r+1$. Let $\gamma'$ be the part of $\gamma$
in the ball $B(v, R+r+1)$ of radius $R+r+1$  centered at $v$. By
the triangle inequality
  we also have $\vert\gamma'\vert \geq r+1.$
So $L_f(\ga)\geq L_f(\ga')\geq (r+1) f(R+r+1).$ We obtain

$$\displaystyle {f(R+r+1)\over f(R-d(x ,m))}\leq{f(R+r+1)\over f(d(v,m))} \leq {
L_f(\ga)\over  (r+1)f(d(v,m))}\leq {L_f(\om)\over
(r+1)f(d(v,m))}\leq {r\over r+1}.$$

\noindent Since $f$ is $\la$-slow we have $\displaystyle
{f(r+R+1)\over f(R-d(x,m))} > {1\over \la^{r+d(m,x)+1}}>{1\over \la^{2r+1}}.$ Thus
$$\displaystyle {1\over \la^{2r+1}} < {r\over
r+1}.\hfill\eqno(*)$$

Then there exists $\lambda_0 > 1$ such that for $\la\in ]1,
\lambda_0[$ the inequality (*) is not true for a fixed $r>0.$ So
for such $\lambda_0$ we have a contradiction. The lemma is
proved.\bx

\bigskip

\noindent {\bf Remark.} Obviously if $r$ is not fixed and tends to
infinity  the above constant  $\la_0$ does not exist.

\bigskip

The group $G$ acts discontinuously and cofinitely on the graph
$\G$ and on the graph $\Ga$ of entourages (see Section 3). Since
the set $\G^0/G=K$ is finite and $\Ga^0/G=\{\ba_0\}\ (\ba_0\in A)$
the correspondence  $K\to \ba_0$ extends $G$-equivariantly to the
quasi-isometry $\psi : gK\to g\ba_0\ (g\in G).$ In the same way we
define the inverse quasi-isometric map $\psi^{-1} : \Ga\to \G$ for
which $\psi^{-1}(\ba_0)\in K.$

For a parabolic point $p\in\P$ let $H$ denote the stabilizer of
$p$ in $G$.

\begin{lem}
\label{parhor} The map $\psi$ extends continuously by the identity
map to the map $\G\sqcup \P\to \Ga\sqcup\P$. Furthermore for any
$d>0$ there exists $d'=d'(d, p)$ such that $\psi(N_d(H))$ belongs
to a $d'$-neighborhood $N_{d'}(T(p))$ of the horosphere
$T(p)\subset \Ga$; and vice versa $\psi^{-1}(N_{d}(T(p))\subset
N_{d'}(H).$

\end{lem}

\proof  It follows from [GePo1, Lemma 3.8] that the unique limit
point of $N_d(H)$ on $T$ is $p$. The set $\psi(N_d(H))$ is an
$H$-finite subset of $\Ga$ and so belong to $N_{d'}(T_p)$ for some
$d'=d'(d,p)$ (see also the proof of Corollary \ref{quasiconv}).
Since the unique limit point of the set $N_{d'}(T_p)$ is also
  $p$ the map $\psi$ extends identically to the set $\P$. The
  second statement is similar.\bx

\bigskip

\begin{lem}
\label{ltfloyd}

For every $l >0$ and $\ve >0$ there exists $\lambda _0 >1$ such
that for any $\la\in]1, \la_0[$ and   $\la$-slow function $f$
satisfying (1-2) one has: if  $\ga\subset\G$ is  $\d_f$-geodesic
then the curve $\psi(\ga)\subset \Ga$ is $(l,c)$-tight  where $c$
is the quasi-isometry constant of $\psi.$

\end{lem}

\proof For a fixed $l>0$ by   Lemma \ref{fwgeod} (applied to
$r=l$) there exists $\la_0
> 1$ such that for any $\la\in ]1, \lambda_0[$ and any $\la$-slow
function $f$, every part of $\ga$ of length less than $l$ is
geodesic in $\G.$ Then $\beta=\psi(\ga)$ is $c$-quasigeodesic on
every interval of length at most $l.$ So the first condition of
Definition \ref{tcurve} is satisfied for $\beta\subset\Ga.$

To prove   \ref{tcurve}.2  assume that $$ \vert\beta(J) \vert >
l,\hfill\eqno(**)$$

\noindent If first ${\rm diam} (\partial\ga(J)) \leq l$ then again
by Lemma \ref{fwgeod} $\ga\vert_J$ is geodesic in $\G$. So
$\beta\vert_J$ is $c$-quasigeodesic in $\Ga$. It follows from (**)
that ${\rm diam} (\partial(\beta(J)))
> c^{-1}(l)=l/c - c$

If now ${\rm diam} (\partial\ga(J)) > l$ then   we have
$\vert\partial \beta(J))\vert
> c^{-1}(l)$ since $\psi$ is a $c$-quasi-isometry.
The lemma is proved. \bx

\bigskip
\noindent Note that the proof of Lemma \ref{ltfloyd} does
 not use the horospheres to prove the tightness condition
 \ref{tcurve}.2. The needed property holds for any part of $\beta$ of length bigger
 than $l.$ The following Corollary shows that it remains valid for a
 curve in $\G$ close in the Floyd metric to a Floyd geodesic
if the latter one does not belong to the graph.

 \begin{cor}
 \label{flap} For every  $l>0$
 there exists $\lambda_0>1$  such that for every $\lambda\in ]1, \lambda_0[$ and
   $\lambda$-slow function
 $f$  if  the Floyd geodesic
 $\ga[x,y]\subset \G_f$     joining two distinct points $x$ and $y$ does not
 belong to $\G,$ then
  there exists a
 curve $\ti\gamma[x,y]\subset\G$ between $x$ and $y$ such that
 $\vert L_f(\ti\ga)-L_f(\ga)\vert\leq\ve$ and every part of $\ti\ga$ of length $l$ is $d$-geodesic.

 Furthermore   the curve
 $\psi(\ti\ga)\subset \Ga$ is
 $(l,c)$-tight for the quasi-isometry constant $c$.
 \end{cor}

 \proof For a fixed $l$ we choose $\lambda$-slow function $f$ such that $\displaystyle
 \lambda \geq {l\over l+1}$.
Suppose that a Floyd geodesic $\ga[x, y]$ intersects the Floyd
boundary $\partial_f\G$. Then for any $\ve>0$ there exists a curve
$\hat\ga:I\to \G$ such
 that $\hat\ga(\partial I)=\{x, y\}$ and
 $\vert L_f(\hat\ga)-L_f(\ga)\vert < \ve.$ Let $x'$ and $y'$ be two points
 on $\hat\ga$ such that $d(x',y')=l$. If the part $\hat\ga[x',y']$
 of $\hat\ga$
 between them is not $d$-geodesic we replace it by a $d$-geodesic
 $\omega=\omega[x',y']$ between $x'$ and $y'.$ Then the $d$-length of the
 obtained curve $\ti\ga$   is strictly less than that of $\hat\ga$.
  Furthermore  by Lemma
\ref{fwgeod} (applied to $r=l$) the curve  $\omega$ is
 also a Floyd geodesic. So we have $$L_f(\ga)\leq L_f(\ti\ga)\leq
 L_f(\hat\ga)\leq L_f(\ga)+\ve.$$    Repeating this procedure with every
 pair of points of $\ti\ga$ situated at the distance $l$ we strictly
 decrease  its $d$-length.
 Since     $d(x,y)\in \Z_{>0}$    after
finitely many steps we obtain a curve (still denoted by $\ti\ga$)
satisfying the first statement.

 Since $\psi:\G\to\Ga$ is a
$c$-quasi-isometry the last part  follows from the argument of
Lemma \ref{ltfloyd}.\bx

 \bigskip

 {\it Proof of  Theorem C.} The group $G$ acts
3-discontinuously and 2-cocompactly on a compactum $T$. Let $\G$
be a locally finite, connected graph  admitting cocompact
discontinuous action of $G.$

  Let $l_0$ and $\la_0$ be the constants given by
Theorem B and    Lemma \ref{ltfloyd} (or Corollary \ref{flap}).
Let $f$ be a $\la$-slow function for  $\la\in ]1, \la_0[$. Suppose
that
 $\ga=\ga(h_1,h_2)\subset \G$ is a $\d_f$-geodesic between two elements $h_1$
 and $h_2$
 in the parabolic subgroup $H.$ Then by Lemma \ref{ltfloyd} the curve
  $\beta=\psi(\ga)$ is $(l,c)$-tight in $\Ga.$

 A segment of a curve $\beta\subset \Ga$ having the extremities at points $\bh_i\in \Ga\ (i=1,2)$
we denote  by $\beta[\bh_1, \bh_2].$ By Lemma \ref{parhor} for
every $d>0$ and $p\in \P$ there exists $d'=d'(d, p)$ such that the
set
 $\psi^{-1}(N_d(T_p))$ belongs to $N_{d'}(H).$ So Theorem C follows from the
following.

\bigskip

\begin{prop}
\label{qconvflent}

For every $c>0$ there exist positive  constants $s$, $d$ and $l_0$
such that for all $l
> l_0$ every $(l,c)$-tight curve $\beta[\bh_1, \bh_2]\subset \Ga$
with $\bh_i\in N_d(T(p))\ (i=1,2)$  is situated in $N_s(T(p))$ for
some $p\in\P$.
\end{prop}

\bigskip

 {\it Proof of the proposition.} Since $\P$  is
$G$-finite it is enough to prove the statement for a fixed $p\in
\P.$ Suppose that $\beta$ is a $(l,c)$-tight curve where $l
>l_0$ and the constants $l_0$ and $c$ are given by  Theorem B.
So there exists a $c'$-quasigeodesic $\al\subset \Ga$ such that
every non-horospherical point $\bv$ of $\beta$ belongs to the
$w_0$-neighborhood $N_{w_0}(\al)$ with respect to the distance
$\da.$ By Lemma \ref{nhor}.1 we have $\forall i\in I\ :\
\da(\al(i), T(p))\leq {\rm const}.$
 Thus there exists a constant  $R=R(d)>0$
 such that for any non-horospherical point $\bv\in \beta$ we have
 $\da(\bv, T(p))\leq R$.

Let now $\beta[\zx, \by]$ be a $d$-horospherical part of $\beta$
lying in $N_d(T(q))$ of another parabolic point $q$. Up to
increasing the above part of $\beta$ we can suppose that both
extremal points $\zx$ and $\by$ are non-horospherical. So we have
$\da(\zx, T(p))\leq R$ and $\da(\by, T(p))\leq R$. Let $\zx_1$ and
$\by_1$ be points on $T(p)$ realizing these distances
respectively. Denote by $\al_1=[\zx, \zx_1]$ and $\al_2=[\by,
\by_1]$ the corresponding geodesics (see Figure below).


\centerline{
\ifnum\pdfoutput>0%
\pdfximage width 270pt height
220pt{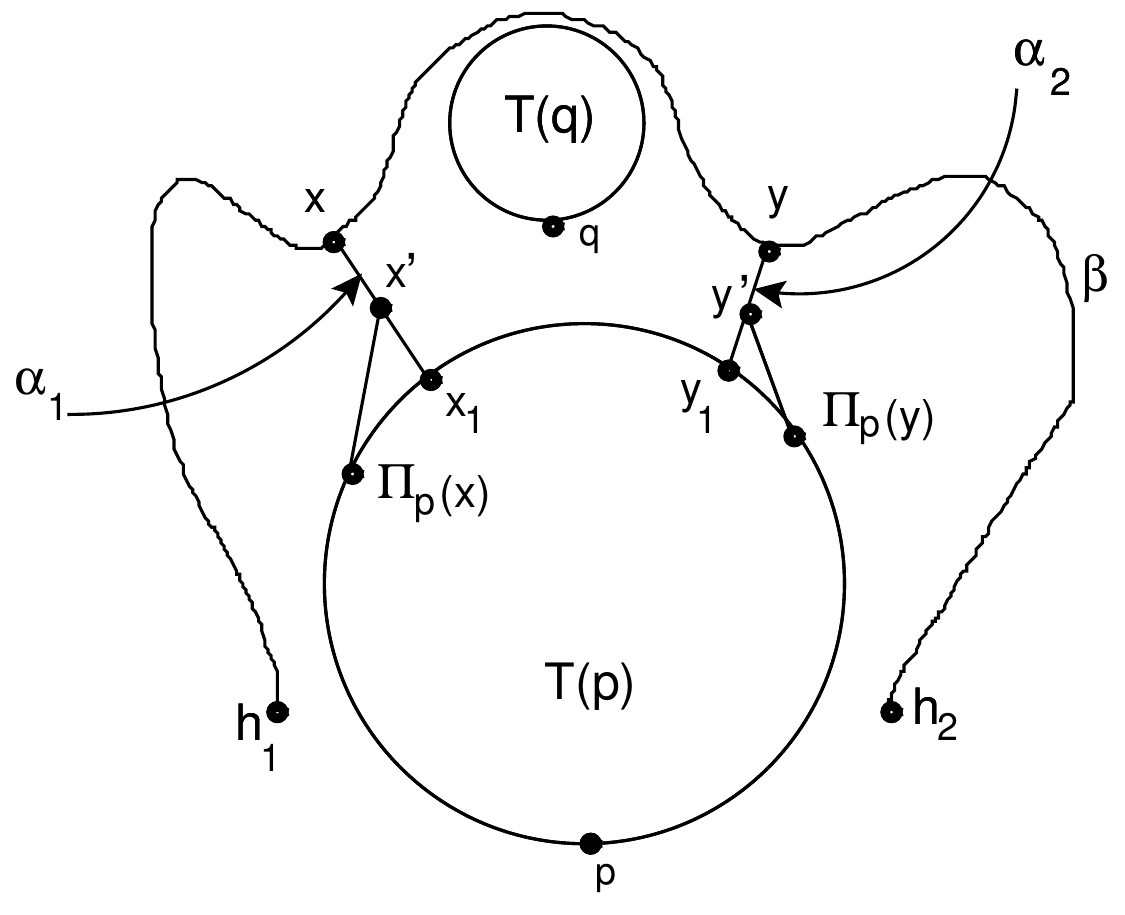}\pdfrefximage\pdflastximage\else\kern270pt\vbox to 220pt{\vss}\fi}

Let $\Pi_p(\zx)$ and $\Pi_p(\by)$ be the projections of $\zx$ and
$\by$ on $T(p)$. By Lemma \ref{prhor}.1 we have $\da(\al_1,
\Pi_p(\zx))=\da(\zx', \Pi_p(\zx))\leq L$ for some constant $L$
depending only on $R$, where $\zx'\in \al_1.$ Hence $\da(\zx,
\Pi_p(\zx))\leq R+L$ and similarly $\da(\by, \Pi_p(\by))\leq R+L.$
By Proposition \ref{finproj}.2 the set $\Pi_p(T(q))$ is finite and
so is $\Pi_p(N_d(T(q)))$. So there exists a constant $C>0$ such
that $\da(\Pi_p\zx, \Pi_p\by) \leq C$. Therefore $\da(\zx, \by)
\leq C+2R+2L.$ The above constants $C$, $R$ and $L$ depend only on
$p$ so we can choose   the parameter $l$ from Theorem B satisfying
$l
> \max (l_0, C+2R+2L).$ Then the segment $\beta[\zx, \by]$ is
$c$-quasigeodesic whose length is bounded by  $c(C+2R+2L)+c.$
Hence $\beta[\zx, \by]\subset N_s(T(p))$ where $s=R+c(C+2R+2L)+c.$
Theorem C is proved.\bx

\bigskip

 Since every parabolic subgroup $H$ is quasiconvex in $G$ there
 exists
   a quasi-isometric map $\varphi$ of the group  $H$ into the
  graph $\G$. We have the following.

\bigskip

\begin{cor}
\label{flbound}  For the  constant $\la_0$ from Theorem C and
every $\la\in ]1, \la_0[$ let $f$ be a $\lambda$-slow Floyd
function satisfying  in addition the following assumption:

$$\displaystyle {f(n)\over f(2n)}\leq \kappa \ (n\in \N)\hfill\eqno(3)$$

\noindent for some constant $\kappa>0.$ Let $p$ be a parabolic
point for the action of $G$ on $T$ and $H={\rm Stab}_Gp$ be its
stabilizer.
 Then $\varphi$
extends injectively to the Floyd boundaries:

$$\varphi\ :\ {\overline H}_f\to {\overline \G}_f.\hfill\eqno(4)$$

\noindent
\end{cor}\bx

\bigskip

\noindent {\bf Remark.} Note that every polynomial type function
$f(n)=(n+1)^{-k}\ (k>1)$ satisfies the conditions (1-3) for any
fixed $\la
>1$ and $\kappa>0\ (n>n_0).$

\bigskip

 {\it Proof of Corollary \ref{flbound}}.  We suppose that
$H\subset \G^0$ and $\varphi : H\hookrightarrow \G^0$ is the
identity map inducing the quasi-isometry between the word metrics.
Let
  $d'(,)$ and $d(,)$ be the graph distances of $H$ and  $\G$
  respectively. We also denote by $\d_{f, H}$ and $\d_{f, G}$ the
  corresponding Floyd distances with respect to a fixed basepoint
  $v\in H.$ Since $f$ satisfies (3)  by [GePo1, Lemma 2.5]   the  map
$\varphi$ extends to a Lipschitz map (denoted by the same letter)
$\varphi : \OFH\to \OFG$ between the Floyd completions of $H$ and
$\G$.

 Let $x, y\in H\subset \G$ be two distinct points. If the Floyd
 geodesic between $x$ and $y$ belongs to  $\G$ we denote it by
 $\gamma$; if not for any $\ve\in ]0,1[$  let  $\gamma$ be the
 $(l,c)$-tight curve ($l>l_0$) given by  Corollary \ref{flap}
 whose Floyd length is $\ve$-close
 to that of the Floyd geodesic.
 In the first case by Theorem C there exists a constant $R=R(H)$ such that
  that $\gamma\subset N_R(H),$ and in the second case the same
  conclusion for the curve $\gamma$ follows from Proposition \ref{qconvflent}.

We have $L_f(\gamma)=\sum_{i=1}^l f(d(v, \{x_i, x_{i+1}\})).$
Denote by $x'_i\in H$ one of the closest vertices  to $x_i$ in $H\
(i=1,...,l).$ By Theorem C there exists a constant $R
>0$ such that $d(x_i, x'_i) \leq R.$ Thus $d(x'_i, x'_{i+1}) \leq
2R+1.$ So for any vertex $x'_{ij}$ on a   geodesic in $H$ between
$x'_i$ and $x'_{i+1}$ we obtain $d(v, \{x_i, x_{i+1}\}) \leq
(3R+1)+d(v, x'_{ij}).$ Since $\varphi$ is quasi-isometric we have
$1/\alpha\cdot d'(v, x'_{ij})-\beta\leq d(v, x'_{ij})\leq \al
d'(v, x'_{ij}) + \beta$ for some constants $\al$ and $\beta$. Let
$\gamma'=\gamma'(x,y)\subset H$ be the curve between $x$ and $y$
obtained by connecting the vertices $x'_i$ and $x'_{i+1}$ by
geodesics segments in $H$ passing through $x'_{ij}$. We have
$\displaystyle \alpha\cdot (d(x'_i, x'_{i+1})+\beta)\cdot f(d(v,
\{x_i, x_{i+1}\}))\geq d'(x'_i, x'_{i+1})\cdot f(d(v, \{x_i,
x_{i+1}\}))= \sum_jf(d(v, \{x_i, x_{i+1}\}).$ Thus $\displaystyle
f(d(v, \{x_i, x_{i+1}\}))\geq {1\over 2\alpha R+\beta+\alpha}
\sum_jf(\al d'(v, x'_{ij})+m_1),$ where $m_1=\beta+3R+1$. The
conditions (1) and (3) yield

 $$L_{f,
G}(\gamma)\geq {L_{f,H}(\gamma')\over (2\alpha
R+\beta+\alpha)\lambda^{m_1}\kappa^{k_1}}\geq {\d_{f, H}(x,y)\over
(2\alpha
R+\beta+\alpha)\lambda^{m_1}\kappa^{k_1}},\hfill\eqno(5)$$

\noindent  where $\displaystyle k_1=\min\{k\ :\ 2^{k} > \al\}.$
Since for every $\ve\in ]0,1[$  there exists a curve $\gamma$
satisfying (5) and for which $L_{f,G}(\gamma)\leq \d_{f,
G}(x,y)+\ve$ we have
$$\forall x,y\in H\ \ \d_{f, G}(x,y)\geq {1\over (2\alpha R+\beta+\alpha)\lambda^{m_1}\kappa^{k_1}}\cdot \d_{f,
H}(x,y).\hfill\eqno(6)$$

\noindent By continuity the   inequality (6) remains valid for
every pair of distinct points $x, y\in\OFH.$ So the map
$\varphi:\OFH\to \OFG$ is injective. The Corollary is proved.\bx

\bigskip

 If $G$ acts on $T$ is 3-discontinuously and $2$-cocompactly then
   the kernel of the
equivariant Floyd map $F$ from the Floyd boundary $\BF$ of the
Cayley graph of $G$ to $T$ is described in [GePo1, Theorem A].
Namely if it is not a single point then it is equal to the
topological boundary $\partial({\rm Stab}_Gp)$ of the stabilizer
${\rm Stab}_Gp$ of a parabolic point $p\in T$. We denote by
$\partial_f{\rm Stab}_Gp$ the Floyd boundary of ${\rm Stab}_Gp$
corresponding to the function $f.$ By Corollary \ref{flbound} we
have that $\varphi\vert_{\partial_fH}$ is a homeomorphism. So the
following is immediate.

\bigskip

\begin{cor}
\label{kermap}

For every $\la\in ]1, \la_0[$ and each $\la$-slow function $f$
satisfying $(1-3)$   one has

$$F^{-1}(p)=\partial_f({\rm Stab}_Gp),\hfill\eqno(7)$$

\noindent for every parabolic point  $p\in T$.
\end{cor}\bx

Corollary \ref{kermap}  answers positively our question [GePo1,
  1.1] and provides complete generalization of the theorem
of Floyd \cite{F} for the relatively hyperbolic groups.

%
%
%
%
%
%
%
%
%


\bigskip



\bigskip


\begin{thebibliography}{Main25}
\addcontentsline{toc}{chapter}{Bibli ography}

\bibitem[BM]{BM}

 A. Beardon, B. Maskit, {\sl Limit sets of Kleinian groups and
finite sided fundamental polyhedra}, Acta Math. \bf 132 \rm(1974)
1--12.




\bibitem[Bo1]{Bo1} B. H. Bowditch, {\sl Relatively hyperbolic
groups,} Internat. J. Algebra Comput. 22 (2012), no. 3,
1250--1316.


\bibitem[Bo2]{Bo2} B. H. Bowditch, {\sl Convergence groups and configuration
spaces},  in ``Group theory down under'' (ed.\ J.Cossey,
C.F.Miller, W.D.Neumann, M.Shapiro), de Gruyter (1999) 23--54.

\bibitem[Bo3]{Bo3}
 B. H. Bowditch, {\sl A topological characterisation of hyperbolic
groups}, J. Amer. Math. Soc.  11 (1998), no. 3  no. 3 643--667.

\bibitem[Bo4]{Bo4}
 B. H. Bowditch, {\sl A short proof that a subquadratic isoperimetric inequality implies
a linear one}, Michigan Math. J. Volume 42, Issue 1 (1995),
103-107.



\bibitem[BH]{BH}
M.~Bridson, A.~Haefliger, {\sl Metric spaces of non-positive
curvature}, Grundlehren der Mathematischen Wissenschaften, 319.
Springer-Verlag, Berlin, 1999.

\bibitem[Bourb]{Bourb}
 N.~Bourbaki, {Topologie G\'en\'erale}
Hermann, Paris, 1965.

\bibitem[BBI]{BBI}
D.~Burago, Y.~Burago, S.~Ivanov, {A Course in Metric Geometry}
Graduate Studies in Math, AMS, Vol 13, 2001.

\bibitem[Da]{Da}
F.~Dahmani, {\sl PhD thesis}, cds.cern.ch/record/851759/files/cer-
002534906.ps.gz?version=1

\bibitem[DS]{DS} C.~Drutu, M.~Sapir, Tree-graded spaces and
asymptotic cones of groups, Topology 44 (2005), no. 5, 959-1058.


\bibitem[Fa]{Fa} B.~Farb, {\sl Relatively hyperbolic groups,}
GAFA \bf 8 \rm no. 5 (1998) 810-840.

\bibitem[F]{F}
W. J. Floyd, {\sl Group completions and limit sets of Kleinian
groups}, Inventiones Math. \bf 57 \rm(1980), 205--218.


\bibitem[Fr]{Fr}
 E. M. Freden, {\sl Properties of convergence groups and spaces}, Conformal
Geometry and Dynamics, \bf 1 \rm(1997) 13--23.

\bibitem[Fu]{Fu}

H. Furstenberg{\sl. Poisson boundaries and envelopes of discrete
groups,} Bull. Amer. Math. Soc. {\bf 73} (1967) 350--356.


\bibitem[Ge1]{Ge1}
V.~Gerasimov, {\sl Expansive Convergence Groups are Relatively
Hyperbolic}, GAFA {\bf 19} (2009) 137--169.

\bibitem[Ge2]{Ge2}
V.~Gerasimov\sl, Floyd maps to the boundaries of relatively
hyperbolic groups\rm, GAFA {\bf 22} (2012) 1361--1399.

\bibitem[GePo1]{GePo1}
V.~Gerasimov, L.~Potyagailo, {\sl Quasi-isometric maps and Floyd
boundaries of relatively hyperbolic groups}, arXiv:0908.0705
[math.GR], 2009,  to appear in Journal of EMS.

\bibitem[GePo3]{GePo3} V.~Gerasimov, L.~Potyagailo, {\sl
Quasiconvexity in the relatively hyperbolic groups},
arXiv:1103.1211 [math.GR], 2011.

\bibitem[GM]{GM}
F. W. Gehring and G. J. Martin\sl, Discrete quasiconformal groups
I\rm. Proc.\ London Math.\ Soc.\ \bf 55 \rm(1987) 331--358.

\bibitem[Gr]{Gr} M. Gromov, {\sl Hyperbolic groups}, in: ``Essays in Group
Theory'' (ed. S.~M.~Gersten) M.S.R.I. Publications No.~8,
Springer-Verlag (1987) 75--263.

\bibitem[Gr1]{Gr1} M. Gromov, {\sl Asymptotic Invariants of Infinite Groups},
``Geometric Group Theory II'' LMS Lecture notes \bf 182\rm,
Cambridge University Press (1993)

\bibitem[GrMa]{GrMa} D. Groves and J.F. Manning, Dehn filling in relatively hyperbolic
groups, Israel J. Math. 168 (2008) 317429.

\bibitem[Hr]{Hr}
C.~Hruska, {\sl Relative hyperbolicity and relative quasiconvexity
for countable groups}, Algebr. Geom. Topol., 10 (2010) 1807--1856.



\bibitem[Ka]{Ka}
 A. Karlsson {\sl Free subgroups of groups with non-trivial Floyd boundary}, Comm.
Algebra, 31, (2003), 5361--5376.

\bibitem[My]{My}
P.~J. Myrberg, {\sl Untersuchungen ueber die automorphen
Funktionen beliebiger vieler Variabelen}, Acta Math. \bf 46
\rm(1925).


\bibitem[Os]{Os}
D.~Osin, {\sl Relatively hyperbolic groups: intrinsic geometry,
algebraic properties and algorithmic problems}, \rm, Mem. AMS \bf
179 \rm(2006) no. 843 vi+100pp.


\bibitem[Tu1]{Tu1}
P. Tukia, {\sl A remark on the paper by Floyd}, Holomorphic
functions and moduli, vol.II(Berkeley CA,1986),165--172, MSRI
Publ.11, Springer New York 1988.

\bibitem[Tu2]{Tu2}
  P. Tukia, {\sl Convergence groups and Gromov's metric
hyperbolic spaces} : New Zealand J.\ Math.\ \bf 23 \rm (1994)
157--187.

\bibitem[Tu3]{Tu3}
 P. Tukia, {\sl Conical limit points and uniform convergence
groups}, J.\ Reine.\ Angew.\ Math.\ {\bf 501} (1998) 71--98.

\bibitem[W]{W}
A.~Weil, {\rm Sur les espace à structure uniforme et sur la
topologie générale}, Paris, 1937.

\bibitem[Ya]{Ya}
A. Yaman, {\sl A topological characterisation of relatively
hyperbolic groups}, J.\ reine ang.\ Math. \bf 566 \rm(2004),
41--89.


\bibitem[Yang]{Yang}
W. Yang, {\sl Peripherical structures of  relatively hyperbolic
groups}, J.\ reine angew. Math.,\ to appear, DOI 10.1515/
crelle-2012-0060.

\end{thebibliography}
\end{document}